\documentclass[11pt,a4paper,righttag]{amsart}
 \usepackage{amsfonts, amssymb}
 \usepackage{overpic}
\newtheorem{teo}{\sc Theorem}[section]

\newtheorem{cor}[teo]{\sc Corollary}
\newtheorem{lemma}[teo]{\sc Lemma}
\newtheorem{proposition}[teo]{\sc Proposition}

\theoremstyle{definition}
\newtheorem{defi}[teo]{\sc Definition}

\theoremstyle{remark}
\newtheorem{rem}[teo]{\sc Remark}

\newtheorem{notation}{\sc Notation}

\newcommand{\bte}{\begin{teo}}
\newcommand{\ete}{\end{teo}}
\newcommand{\bc}{\begin{cor}}
\newcommand{\ec}{\end{cor}}
\newcommand{\bp}{\begin{pro}}
\newcommand{\ep}{\end{pro}}
\newcommand{\bl}{\begin{lemma}}
\newcommand{\el}{\end{lemma}}
\newcommand{\bd}{\begin{defi}}
\newcommand{\ed}{\end{defi}}
\newcommand{\bno}{\begin{notation}}
\newcommand{\eno}{\end{notation}}

\newcommand{\bca}{\begin{cases}}
\newcommand{\eca}{\end{cases}}
\newcommand{\la}{\langle}
\newcommand{\ra}{\rangle}

\newcommand{\bq}{\begin{equation}}
\newcommand{\eq}{\end{equation}}
\newcommand{\btabu}{\begin{table}}
\newcommand{\etabu}{\end{table}}
\newcommand{\bt}{\begin{tabular}}
\newcommand{\et}{\end{tabular}}
\newcommand{\ba}{\begin{array}}
\newcommand{\ea}{\end{array}}
\newcommand{\br}{\begin{eqnarray}}
\newcommand{\er}{\end{eqnarray}}
\newcommand{\brn}{\begin{eqnarray*}}
\newcommand{\ern}{\end{eqnarray*}}
\newcommand{\benu}{\begin{enumerate}}
\newcommand{\eenu}{\end{enumerate}}
\newcommand{\bite}{\begin{itemize}}
\newcommand{\eite}{\end{itemize}}

\newcommand{\supp}{\operatorname{supp }}

\newcommand{\field}[1]{\mathbb{#1}}

\newcommand{\R}{\field{R}}

\newcommand{\C}{\field{C}}

\renewcommand{\Re}{\mathop{\rm Re}}
\renewcommand{\Im}{\mathop{\rm Im}}

\title[Weak asymptotics of multiple orthogonal polynomials] {On Nikishin systems with discrete components and weak asymptotics of multiple orthogonal polynomials}

\thanks{The work of the first author  was supported by a grant from the  Russian Science Foundation project
142100025. The second and the third authors were supported by MICINN of Spain under grants MTM2015-65888-C4-2-P and MTM2011-28952-C02-01, respectively, and by the European Regional Development Fund (ERDF). Additionally, the third author was supported by Junta de Andaluc\'ia (the Excellence Grant P11-FQM-7276 and the research group FQM-229) and by Campus de Excelencia Internacional del Mar (CEIMAR) of the University of Almer\'ia.}

\dedicatory{Dedicated to our teachers and friends Andrei Alexandrovich Gonchar, Eugene Mikhailovich  Nikishin, and Herbert Stahl}

\author[Aptekarev]{A. I. Aptekarev}

\address[Aptekarev]{Keldysh Institute of Applied Mathematics, Moscow, Russia}

\email[Aptekarev]{aptekaa@keldysh.ru}

\author[L\'opez-Lagomasino]{G. L\'opez Lagomasino}
\address[L\'opez-Lagomasino]{Department of Mathematics,
Universidad Carlos III de Madrid, Legan\'es, Spain}

\email[L\'opez]{lago@math.uc3m.es}

\author[Mart\'inez-Finkelshtein]{A. Mart\'inez-Finkelshtein}

\address[Mart\'inez-Finkelshtein]{Department of Mathematics, University of Almer\'ia, Almer\'ia, Spain}

\email[Mart\'inez]{andrei@ual.es}

\begin{document}

\maketitle


\begin{abstract}
We consider multiple orthogonal polynomials with respect to Nikishin systems generated by two measures $(\sigma_1, \sigma_2)$ with unbounded supports ($\mbox{supp} \, \sigma_1 \subseteq \mathbb{R}_+$, $\mbox{supp} \, \sigma_2 \subseteq  \,(-\infty,0)$) and $\sigma_2$ discrete. A Nikishin type equilibrium problem in the presence of an external field acting on $\mathbb{R}_+$ and a constraint on $ \mathbb{R}_-$ is stated and solved. The solution is used for deriving the contracted zero distribution of the associated multiple orthogonal polynomials.

Bibliography: 56 titles.
 \end{abstract} \vspace{1cm}

\vspace{.3cm}

{\it A.M.S. Subject Classification.} Primary: 30E10, 42C05;
Secondary: 41A20.
 
\vspace{1cm}

\section{Introduction}

In a celebrated paper published in 1980, E.~M.~Nikishin \cite{Nik1} introduced a general class of systems of measures, now called Nikishin systems. Let $\Delta_{\alpha}, \Delta_{\beta}$
be two non-intersecting bounded intervals of the real line ${\mathbb{R}}$, measures  $\sigma_{\alpha} \in {\mathcal{M}}(\Delta_{\alpha})$ and $\sigma_{\beta} \in {\mathcal{M}}(\Delta_{\beta})$, where $\mathcal{M}(\Delta)$ denotes the set of all finite Borel measures on the interval $\Delta$ with constant sign. With $\sigma_{\alpha}$ and $\sigma_{\beta}$ we construct a third measure $\langle \sigma_{\alpha}, \sigma_{\beta} \rangle$, which using the differential notation is given by
\begin{equation} \label{CT} d \langle \sigma_{\alpha}, \sigma_{\beta} \rangle(x) := \widehat{\sigma}_{\beta}(x) d\sigma_{\alpha}(x), \qquad \widehat{\sigma}_{\beta}(x) = \int (x-t)^{-1} d\sigma_{\beta}(t).
\end{equation}
\begin{defi} \label{Nikishin} Take a collection  $\Delta_j$, $j=1,\ldots,m,$ of intervals such that
\[ \Delta_j \cap \Delta_{j+1} = \emptyset, \qquad j=1,\dots,m-1,
\]
and a system of measures $(\sigma_1,\dots,\sigma_m)$ with $\sigma_j \in {\mathcal{M}}(\Delta_j)$, $j=1,\dots, m$; we assume additionally that for each $j$, the convex hull of the support $ \supp (\sigma_j)$ of $\sigma_j$ coincides with $\Delta_j$. Let
\[ s_1 = \sigma_1, \quad s_2 = \la \sigma_1,\sigma_2 \ra, \quad \dots , \quad  s_m = \la \sigma_1,\la \sigma_2,\ldots,\sigma_m \ra \ra.
\]
We say that $(s_1,\dots,s_m)$ is the \emph{Nikishin system of measures} generated by $(\sigma_1,\dots,\sigma_m)$, and denote it by $(s_1,\dots,s_m)  = {\mathcal{N}}(\sigma_1,\dots,\sigma_m)$.
\end{defi}

This model system was introduced in order to study general properties of multiple orthogonal polynomials and Hermite-Pad\'e approximants.

Fix ${\bf n}  := (n_1,\ldots,n_m) \in {\mathbb{Z}}_+^m \setminus \{\bf 0\}$, where ${\bf 0}$ is the $m$ dimensional zero vector. Define $P_{\bf n}$ as a non-zero polynomial of degree $\deg (P_{\bf n}) \leq |{\bf n}| := n_1+\cdots +n_m$ such that
\[ \int x^{\nu} P_{\bf n}(x) d s_j(x) = 0,   \qquad \nu=0,\ldots,n_j -1, \qquad j=1,\ldots,m.
\]
The existence of $P_{\bf n}$ reduces to solving a homogeneous linear system of $|{\bf n}|$ equations on the  $|{\bf n}| +1$ coefficients of $P_{\bf n}$; therefore, a non-trivial solution is guaranteed. However, in contrast with the scalar case $(m=1)$ of standard orthogonal polynomials (OP), uniqueness up to a constant factor is not a trivial matter (and, in general, not true for systems of arbitrary measures $(s_1,\ldots,s_m)$). In connection with this question in \cite{Nik1} it  was shown that in presence of a Nikishin system uniqueness holds, with $\deg P_{\bf n} = |{\bf n}|$, for multi-indices of the form $(n+1,\ldots,n+1,n,\ldots,n)$, and stated without proof that it is also true whenever $n_1\geq\cdots\geq n_m$. In the sequel we assume that $P_{\bf n}$ is monic.

Motivated by the structure of Nikishin systems, Herbert Stahl studied their analytic and algebraic properties (see \cite{AptSta}). In a series of papers \cite{DrSt0}, \cite{DrSt1}, \cite{DrSt2}, among other results, K.~Driver and H.~Stahl showed that uniqueness remains valid whenever $n_j \leq n_k+1$, $1\leq k < j\leq m$. The problem for arbitrary multi-indices was definitely solved in \cite{FL1} (and \cite{FL} when the generating measures have unbounded and/or touching supports).

A remarkable property of Nikishin orthogonal polynomials is that they not only share orthogonality relations with respect to several measures but they also satisfy full orthogonality relations with respect to a single (varying with respect to $\bf n$) measure. For $m=2$ and $n_2 \leq n_1 +1$ this was first observed by Andrei Aleksandrovich Gonchar\footnote{On one of the regular Monday seminars at the Steklov Institute A.~A.~Gonchar was reporting on the results contained in \cite{Nik1} but after a short while he had to leave to attend an important meeting. After  an hour or so he returned and started anew his presentation proving \eqref{eq:a} and \eqref{eq:b} and from there deduced the convergence of the corresponding Hermite-Pad\'e approximants.}  by showing that the  function of the second kind
\[ R_{{\bf n},1}(z) = \int \frac{P_{\bf n}(x)}{z-x}\,d\sigma_1(x)
\]
satisfies the orthogonality relations
\begin{equation} \label{eq:a}
   \int x^{\nu}R_{{\bf n},1}(x) d\sigma_2(x) = 0, \qquad \nu =0,\dots,n_2 -1.
\end{equation}
From here it follows that $R_{{\bf n},1}$ has exactly $n_2$   zeros in ${\mathbb{C}}\setminus \Delta_1$, they are all simple, and lie in the interior of $\Delta_2$. If $P_{{\bf n},2}$ denotes the monic polynomial of degree $n_2$ vanishing at these points, then
\begin{equation}
  \label{eq:b}
  \int x^{\nu}P_{{\bf n}}(x) \frac{d\sigma_1(x)}{P_{{\bf n},2}(x)} = 0, \qquad \nu =0,\dots,n_1 + n_2 -1.
\end{equation}

The study of the asymptotic behavior of multiple orthogonal polynomials is greatly indebted to A.~A.~Gonchar. In joint papers with E.~A.~Rakhmanov \cite{GR1}, \cite{GR2}, \cite{GR3}, they introduced the notion of vector equilibrium problem to describe the asymptotic zero distribution of such polynomials.
For a Nikishin system of two measures and $n_1=n_2=n$ the result may be stated as follows. Define the normalized zero counting measure $\nu_P$ of a polynomial $P$ as
\[ \nu_P = \frac{1}{\deg P}\sum_{P(x) =0} \delta_x,
\]
where $\delta_x$ denotes the Dirac measure with mass $1$ at the point $x$, and each zero of~$P$ is taken
with account of its multiplicity, so that the total variation $|\nu_P|$ of $\nu_P$ is 1.
Assume that $\sigma_j \in \mbox{\bf Reg}$, $j=1,2$ (for the definition of the class $\mbox{\bf Reg}$
of measures, see \cite[Chapter 3]{stto}). Then there exist positive measures $\lambda_j \in {\mathcal{M}}(\Delta_j)$,
$j=1,2$, $|\lambda_1| = 2$, $|\lambda_2| = 1$,
such that
\begin{equation} \label{limZer}
\lim_n \nu_{P_{\bf n}} = \lambda_1/2, \qquad \lim_n \nu_{P_{{\bf n},2}}  = \lambda_2,
\end{equation}
in the weak-* topology of measures, where $\lambda_1$ and $\lambda_2$ are uniquely determined by the solution of the vector equilibrium problem
\begin{equation} \label{equiN}
\begin{array}{c}  2 U^{\lambda_1}(x) - U^{\lambda_2}(x) \,\,
\left\{
\begin{array}{cl}
= w_1, & x \in \supp(\lambda_1), \\
\geq w_1, &   x \in \Delta_1 \setminus \supp(\lambda_1),
\end{array}
\right.
\\\\
    2 U^{\lambda_2}(x) - U^{\lambda_1}(x)  \,\,
\left\{
\begin{array}{cl}
= w_2, & x \in \supp(\lambda_2), \\
\geq w_2, &   x \in \Delta_2 \setminus \supp(\lambda_2),
\end{array}
\right.
\end{array}
\end{equation}
where $w_1, w_2$ are certain constants, and $U^{\lambda}$ denotes the logarithmic potential of $\lambda$ (see the definition below). At the time, this result and its extensions were well known within a small circle of specialists. With some variations, for general Nikishin systems it appeared in papers by H. Stahl \cite{St}, and with the highest degree of generality by A.~A.~Gonchar, E.~A.~Rakhmanov, and V.~N.~Sorokin \cite{GRS}. For other extensions and generalizations see  \cite{AptKal}, \cite{AL}, \cite{BLM}, \cite{Bus},  \cite{LF2},  \cite{NikTypeI}, \cite{Rak}, \cite{RS}.

In recent years, Nikishin systems have attracted new attention because this construction has been identified in different models of random matrix theory and multiple orthogonal polynomial ensembles,  see  \cite{aptkuij}, \cite{Kuij}, and \cite{MR2470930}. In some of these models new ingredients appear in which some of the generating measures turn out to be discrete and/or have unbounded support. V.~N.~Sorokin has studied the asymptotic distribution of the zeros for several multiple orthogonal polynomials of this type, see \cite{Sorokin:2009fk}--\cite{Sor2}.

Orthogonal polynomials with respect to discrete measures have the characteristic that between two consecutive mass points there may be at most one zero of the  polynomial. This fact induces a constraint on the equilibrium problem whose solution describes the asymptotic zero distribution of the orthogonal polynomials. This effect was first pointed out by E.~A.~Rakhmanov in \cite{Rak2} (see also \cite{DrSa} and \cite{KW}). A similar situation occurs in the case of multiple orthogonal polynomials.

The present paper is devoted to the study of multiple orthogonal polynomials with respect to Nikishin systems generated by two measures $(\sigma_1, \sigma_2)$ with unbounded supports
\[ \mbox{supp} ( \sigma_1) \subseteq \mathbb{R}_+ := [0,+\infty) , \qquad \mbox{supp} ( \sigma_2 ) \subset (-\infty,0).
\]
The second measure $\sigma_2$ is  discrete. To obtain the limiting zero distribution \eqref{limZer} of such multiple OP we state and solve a Nikishin type equilibrium problem which generalizes  \eqref{equiN} by having an external field acting on $\mathbb{R}_+$ and a constraint on $\mathbb{R}_- := (-\infty,0]$.

The main results are stated in Section 2. In Section 3 we review some examples of explicit solutions of the type of equilibrium problems that we consider. Section 4 contains new results related with potentials with unbounded support and scalar equilibrium problems. The last two sections include the proofs of the main results.

\bigskip

This paper has a long story. It was started in 2011 while the first author visited Spain in the framework of  the Excellence Chair Program sponsored by Universidad Carlos III
de Madrid and the Bank of Santander. Then, an essential progress on this project was achieved in 2014 when the Editorial Boards  of Sbornik Mathematics and Journal of Approximation Theory were preparing the special issues \cite{Gon-St1} and \cite{Gon-St2} of their journals in memory of A.A. Gonchar  (1931 - 2013) and H. Stahl  (1945 - 2012). However, it was impossible for us to complete the task in due form.   Finally, the 70th anniversary of E. M. Nikishin's birthday in 2015 and the 30th anniversary of his death in 2016 motivated the authors to conclude the work, which is dedicated to the memory of these outstanding analysts.

\section{Statement of the main results} \label{sec2}

Let $d\sigma_1(x) = \sigma_1'(x) dx$ be a positive, absolutely continuous measure on $\R_+$, and $\sigma_2$ a purely discrete measure whose support is contained in $(-\infty,0)$ given by
\begin{equation}\label{sig2}
\sigma_2 = \sum_{k \geq 1} \beta_k \, \delta_{t_k}, \quad 0> t_k \searrow - \infty, \qquad \beta_k > 0, \qquad \sum_{k\geq 1} \frac{\beta_k}{|t_k|}  < +\infty.
\end{equation}
All the moments of $\sigma_1$ are assumed to be finite. Notice that $\widehat{\sigma}_2$ is integrable with respect to $\sigma_1$.  Let $(s_1,s_2) = {\mathcal{N}}(\sigma_1,\sigma_2)$ be the Nikishin system generated by these measures. For ${\bf n} = (n_1,n_2) \in {\mathbb{Z}}_+^2 \setminus \{0\}$ we define $P_{\bf n}$ as the monic polynomial of degree $|{\bf n}|$ which satisfies
\begin{equation}\label{eq:NikPs}
\int x^{\nu} P_{\bf n}(x) ds_j(x) = 0, \qquad \nu=0,\dots,n_j-1, \quad j=1,2.
\end{equation}
The zeros of $P_{\bf n}$ are simple and lie in the interior of $\mathbb{R}_+$.
We will restrict our attention to sequences of multi-indices of the form ${\bf n} = (n,n)$. In order to simplify the notation we write $P_n$ instead of $P_{\bf n} $. Thus, $\deg P_n = 2n$. Our goal is to describe the (rescaled) asymptotic zero distribution of the polynomials $\left(P_{n}\right), n \in \mathbb{N}$, under appropriate assumptions on the generating measures $\sigma_j, j=1,2$.

Using the properties of Nikishin systems (see \cite{FL} and \cite{GRS}) it is easy to deduce that there exists a monic polynomial $P_{n,2}$, $\deg P_{n,2} = n$, whose zeros are simple and contained in the interior of the convex hull of $\supp(\sigma_{2})$, such that
\begin{equation}\label{eq:3^*P}
\int  x^{\nu} \frac{P_{n}(x)}{P_{n,2}(x)} {d\sigma_1(x)}    =0 , \qquad \nu =0,\ldots, 2n -1,
\end{equation}
and
\begin{equation}\label{eq:4^*P}
\int  t^{\nu} \frac{P_{n,2}(t)}{P_{n}(t)}\int  \frac{P_{n}^2(x)}{P_{n,2}(x)} \frac{d\sigma_1(x)}{x-t }d \sigma_{2}(t)  =0 , \qquad \nu =0,\ldots,  n  -1.
\end{equation}
In other words, $P_n$ and $P_{n,2}$ satisfy full orthogonality relations with respect to varying measures.

Let $(d_n)_{n \in {\mathbb{Z}}_+}$, $d_n \geq 1$, $\lim_n d_n^{1/n} = 1,$ be an increasing  sequence of numbers, and let
\begin{equation}
\label{scaling}
Q_n(x) = P_n(d_n x)/d_n^{2n}, \qquad Q_{n,2}(t) = P_{n,2}(d_n t)/d_n^{n}.
\end{equation}
Making the change of variables $x \to d_nx$, $t \to d_n t$  it follows that the monic polynomials $Q_n$, $Q_{n,2}$ verify the orthogonality relations
\begin{equation} \label{eq:3^*}
\int  x^{\nu} \frac{Q_{n}(x)}{Q_{n,2}(x)} {\sigma_1'(d_n x)d x}    =0 , \qquad \nu =0,\ldots, 2n -1,
\end{equation}
and
\begin{equation} \label{eq:4^*}
\int  t^{\nu} \frac{Q_{n,2}(t)}{Q_{n}(t)}\int  \frac{Q_{n}^2(x)}{Q_{n,2}(x)} \frac{\sigma_1'(d_n x)d x}{x-t }\, d \sigma_{2,n}(t)  =0 , \qquad \nu =0,\ldots,  n  -1,
\end{equation}
where
\begin{equation}
\label{sigma2scaled}
\sigma_{2,n} = \sum_{k \geq 1} \beta_k \delta_{\xi_{k,n}}, \qquad \xi_{k,n} = t_k/d_n.
\end{equation}

\medskip

The asymptotic zero distribution of the multiple orthogonal polynomials $Q_n$, $Q_{n,2}$ is described in terms of an associated vector equilibrium problem that we now present.

For a closed subset $\Delta \subset {\mathbb{R}}$ we denote by
${\mathcal{M}}^+(\Delta)$ the class of all finite positive Borel measures $\mu$ such that $\supp(\mu) \subset \Delta$. We write $\mu \in {\mathcal{M}}_c^+(\Delta)$ if, additionally, $|\mu| =c$.
Let $\mu \in {\mathcal{M}}^+({\mathbb{R}})$.
Its logarithmic potential and energy are given by
\begin{equation}
\label{def_log_pot}
U^{\mu}(x) := \int \log \frac{1}{| x- y|}  d\mu (y), \quad I(\mu) := \int \int \log \frac{1}{| x- y|} d\mu(x) d\mu (y),
\end{equation}
respectively, whenever these integrals are well defined.

Assume that $\mu_1,\mu_2 \in \mathcal{M}^+({\mathbb{R}})$ verify
\begin{equation}\label{C1} I(\mu) < +\infty, \qquad \int \log(1+|x|^2) d\mu(x) < +\infty.
\end{equation}
Their mutual energy may be defined as
\[ I(\mu_1,\mu_2) := \int\int \log \frac{1}{| x- y|} d\mu_1(x) d\mu_2 (y).
\]
Analogously, one can define the potential, energy, and mutual energy of signed measures. In particular, if $\mu_1$ and $\mu_2$ verify \eqref{C1},  then
\[ I(\mu_1 - \mu_2) = I(\mu_1) + I(\mu_2) - 2I(\mu_1,\mu_2).
\]
Moreover, if $\mu_1,\mu_2 \in {\mathcal{M}}_c^+({\mathbb{R}})$ (only finite energy is required), we have
\begin{equation} \label{desig} I(\mu_1 - \mu_2) \geq 0,
\end{equation}
with equality if and only if $\mu_1 = \mu_2$ (see \cite[Theorem 2.5]{CKL},  \cite[Theorem 4.1]{Sim}, and also \cite[Lemma 1.1.8]{ST} if the measures have compact support).

Let $\sigma$ be a positive Borel measure, $\supp (\sigma) = {\R}_-$, $|\sigma| > 1$,  such that for every compact subset $K \subset {\R}_-$ we have that $U^{\sigma|_K}$ is continuous on $\C$, where $\sigma|_K$ denotes the restriction of $\sigma$ to $K$. We define
\begin{equation} \label{star}
{\mathfrak{M}}(\sigma) := \{\vec{\mu} = (\mu_1,\mu_2)^t \in {\mathcal{M}}_2^+({\mathbb{R}}_+) \times {\mathcal{M}}_1^+({\mathbb{R}}_-):\,  \mu_2 \leq \sigma\},
\end{equation}
where the superscript $t$ stands for transpose. By $\mu_2 \leq \sigma$ we mean that $\sigma - \mu_2$ is a positive measure. Since we have assumed that $U^{\sigma|_K}$ is continuous on $\C$ for every compact $K$, it readily follows that $U^{\mu_2}$ is continuous on $\C$.
Eventually we will require that a measure $\mu$ on $\R$ (in particular $\sigma$) satisfies the condition that for every $\varepsilon > 0$ there exists $0 < \delta < 1/2$ and $R_0 >0$ such that
\begin{equation}
\label{cond5*}
\sup_{|R| \geq R_0} \int_{R-\delta}^{R +\delta} \log \frac{1}{|R - y|} d\mu(y) < \varepsilon.
\end{equation}

Let $\varphi$ be a real valued continuous function on $\R_+$ satisfying
\begin{equation} \label{cond1*} \lim_{x \to \infty} (\varphi(x) - 4 \log x) = +\infty .
\end{equation}
Set
\[{\mathfrak{M}}^*(\sigma) := \{\vec{\mu} \in {\mathfrak{M}}(\sigma): \mu_1,\,  \, \mu_2\,\,\, \mbox{verify}\,\,\, \eqref{C1} \},
    \]
\[ J_{\varphi} := \inf\{J_{\varphi}(\vec{\mu}):\vec{\mu} \in {\mathfrak{M}}^*(\sigma) \}, \qquad J_{\varphi}(\vec{\mu}) := 2 \left( I(\mu_1) - I(\mu_1,\mu_2) + I(\mu_2) +\int \varphi\,  d\mu_1\right),
\]
and
\[  W_1^{\vec{\mu}}(x) := 2 {U}^{ \mu_1}(x) - {U}^{ \mu_2}(x) + \varphi(x), \qquad  W_2^{\vec{\lambda}}(x) := 2 {U}^{ \lambda_2}(x) - {U}^{ \lambda_1}(x).
\]

\medskip

\begin{teo} \label{teo:equil} Let $\sigma, \supp (\sigma) = {\mathbb{R}}_-, |\sigma| > 1$,  be a positive Borel measure such that for every compact subset $K \subset {\R}_-$ we have that $U^{\sigma|_K}$ is continuous on $\C$. Let $\varphi$  be a continuous function on $\R_+$ which verifies $\eqref{cond1*}$.
Then, the following statements are equivalent and have the same unique solution:
\begin{itemize}
\item[$(A)$] There  exists $\vec{\lambda} \in {\mathfrak{M}}^*(\sigma)$ such that ${{J}}_{\varphi}(\vec{\lambda}) = {{J}}_{\varphi} > -\infty$.
\item[$(B)$] There  exists $\vec{\lambda} \in {\mathfrak{M}}^*(\sigma)$ such that  for all $\vec{\nu} \in \mathfrak{M}^*(\sigma)$,
    \[\int {{W}}_1^{\vec{\lambda}} d({\nu_1}-{\lambda_1}) + \int {{W}}_2^{\vec{\lambda}} d({\nu_2}-{\lambda_2})\geq 0.
    \]
\item[$(C)$] There  exist  $\vec{\lambda} = ( \lambda_1,  \lambda_2) \in {\mathfrak{M}}^*(\sigma)$ and constants $w_1 = w_1(\sigma,\varphi),w_2 = w_2(\sigma,\varphi)$  such that
\begin{equation}\label{EqCo1p}   2 {U}^{ \lambda_1}(x) - {U}^{ \lambda_2}(x) + \varphi(x)
\left\{
\begin{array}{ll}
= w_1, & x \in \supp( \lambda_1) , \\
\geq w_1, & x \in {\mathbb{R}}_+,
\end{array}
\right.
\end{equation}
\begin{equation}\label{EqCo2p} 2 {U}^{ \lambda_2}(x) - {U}^{ \lambda_1}(x)
\left\{
\begin{array}{ll}
\leq w_2, & x \in \supp( \lambda_2) = \R_-, \\
=  w_2, & x \in \supp(\sigma -  \lambda_2).
\end{array}
\right.
\end{equation}
\end{itemize}
The constants $w_1,w_2$ are uniquely determined as well. We also have that $U^{\lambda_1},U^{\lambda_2}$ are continuous on $\mathbb{C}$, $\supp(\lambda_1)$ is compact.  If  $x\varphi'(x)>0$ is increasing on ${\R}_+$ then $\supp(\lambda_1)$ is also connected. If $\varphi$ is increasing on $\R_+$ then $0 \in \supp(\lambda_1)$. If $\int \log(1+y^2) d\sigma(y) = +\infty$ and $\sigma$ verifies \eqref{cond5*} then $w_2(\sigma,\varphi) = 0$.
\end{teo}

Results of this nature (in a more general setting regarding the dimension of the vector equilibrium problem and the supports of the corresponding measures) may be seen in \cite{BKMW}. There,  the action of constraints on the measures is not considered and the external fields, which verify restrictions of the form \eqref{cond1}, act on all the components of the vector measures. This implies in turn that all the components of the equilibrium vector measure have compact support. However, taking into consideration certain applications, we are especially interested in allowing the second component of the equilibrium measure to be unbounded. For this reason, in the proof of Theorem \ref{teo:equil} (see also Lemma \ref{lem1}) we follow the approach presented in \cite{HK} where results similar to Theorem \ref{teo:equil}, except for part $(C)$, also appear. It's worth mentioning that when dealing with vector potentials involving measures with overlapping supports, in general, there is no reason for Euler-Lagrange variational conditions to hold everywhere, even in the presence of positive definite interaction matrices\footnote{see Section 5 for the definition of the interaction matrix relevant in our case} and strongly confining external fields (see the interesting examples contained in \cite{BKMW}). In our case, the solution is due to the Nikishin type structure of the problem and the action of the constraint $\sigma$ satisfying adequate conditions.

\medskip

In order to study the contracted zero distribution of the polynomials $Q_{n}, Q_{n,2}$, we must impose some restrictions on  the points $\xi_{k,n}$ and the numbers $\beta_k, d_n$. These conditions are inspired by similar ones introduced for the study of the contracted zero distribution of discrete orthogonal polynomials in the scalar case as you can see in \cite[Theorem 2]{Rak2}, \cite[Definition 3.1]{DrSa}, \cite[Section 6]{KW}, and  \cite[Theorem 7.1]{KR} whose model we follow closely. In the sequel we assume that:
\begin{itemize}
\item[(i)] There exists a positive continuous function $\rho$ on $\R_-$   such that
\[ |\xi_{k+1,n} - \xi_{k,n}| > \rho(\xi_{k,n})/n, \qquad k \geq 0 \qquad (\xi_{0,n} = 0).
\]
\item[(ii)] \[\lim_{n\to \infty}\left(\min \{\beta_k: \xi_{k,n} \in [-n,0]\}\right)^{1/n} = 1.\]
\item[(iii)] There exists a positive Borel measure $\sigma, \supp (\sigma) = {\mathbb{R}}_-, |\sigma| > 1$,  such that:
\begin{itemize}
\item for every compact subset $K \subset {\R}_-$, the logarithmic potential $U^{\sigma|_K}$   of the restriction of $\sigma$ to $K$ is continuous on $\C$,
\item \[\int \log (1+y^2) d\sigma(y) = +\infty,\]
\item  for every $\varepsilon > 0$ there exists $0 < \delta < 1/2$ and $R_0 < 0$ verifying
\[
\sup_{R \leq R_0} \int_{R-\delta}^{R +\delta} \log \frac{1}{|R - y|} d\sigma(y) < \varepsilon,
\]
\end{itemize}
and
\begin{equation} \label{limsigma}
\lim_{n \to \infty} \frac{1}{n} \int  f(x) d\,(\sum_{k \geq 1} \delta_{\xi_{k,n}}) (x)= \lim_{n \to \infty} \frac{1}{n} \sum_{k \geq 1} f(\xi_{k,n}) = \int f(x) d\sigma(x)
\end{equation}
for every continuous function $f$ with compact support in $\mathbb{R}_-$.
\item[(iv)] There exists a continuous function $\varphi$ on $\R_+$ satisfying
\begin{equation} \label{cond1}
\liminf_{x\to +\infty}  \varphi(x)/( 4\log x ) > 1.
\end{equation}
such that for a certain $\alpha < 1$
\begin{equation} \label{b} \lim_{n \to \infty} \frac{1}{n}\log (x^{\alpha} \sigma_1'(d_n x)) = -{\varphi}(x)
\end{equation}
uniformly on each compact subset of $\mathbb{R}_+$, and
\begin{equation} \label{cond3}
\liminf_{n \to \infty, x \to +\infty} \frac{- \log (x^{\alpha} \sigma_1'(d_n x))}{4n \log x} > 1.
\end{equation}
\end{itemize}

\medskip

Now we are ready to formulate the main result about the zero asymptotics of Nikishin orthogonal polynomial.

\begin{teo} \label{teo:weak}
Let the  assumptions ${\rm (i)-(iv)}$ formulated above hold, and let $\vec{\lambda}=(\lambda_1, \lambda_2) \in {\mathfrak{M}}^*(\sigma) $ be the solution of the extremal problem in Theorem~\ref{teo:equil}. Assume that
\begin{equation}
\label{int}
 0 \not\in \supp(\sigma-\lambda_2), \qquad \int |y|^{\alpha} d\lambda_2(y) < \infty, \qquad \lambda > 1/2.
\end{equation}
Then
\begin{equation} \label{lessweak} \lim_n \nu_{ Q_{n}} = \lambda_1/2, \qquad \lim_n \nu_{ Q_{n,2}} = \lambda_2,
\end{equation}
in the weak topology of measures. That is for every bounded continuous functions $f$ and $g$ on $\R_+$ and $\R_-$, respectively, we have
\[ \lim_n \int f d \nu_{ Q_{n}} = \frac{1}{2} \int f d \lambda_1, \qquad \lim_n \int g d \nu_{ Q_{n,2}} =   \int g d \lambda_2.
\]
\end{teo}

Although the assumptions of this theorem may seem too restrictive, it encompasses many interesting examples. Some of them will be discussed in the next section. In particular, we will analyze briefly the case of the modified Bessel weights (appearing in the analysis of the non-intersecting squared Bessel paths), the multiple Hermite polynomials (useful when studying ensembles of random matrices with an external source), and finally, the multiple Pollaczek polynomials, studied previously in  \cite{Sorokin:2009fk}, which will be discussed in more detail, and for which an alternative method for solving the equilibrium problem of Theorem~\ref{teo:equil} is presented. These examples verify all the assumptions of Theorem \ref{teo:weak} except the integral condition in \eqref{int}. It remains a difficult unsolved problem to eliminate this condition from a general theorem like Theorem \ref{teo:weak}.

\medskip

Let us finish this section noting that we can easily translate the results of  Theorems~\ref{teo:equil} and \ref{teo:weak} to the equivalent setting of the whole real axis ${\mathbb{R}}$ (with symmetric measures with respect to the origin). Indeed, let $\{P_m\}$ be a sequence of multiple orthogonal polynomials satisfying \eqref{eq:NikPs} with respect to a Nikishin system \eqref{eq:3^*P}--\eqref{eq:4^*P} on the semiaxis  ${\mathbb{R}}_+$, and define the polynomial sequence  $\{\tilde{P}_n\}$  with polynomials of even degrees by
\begin{equation}\label{MapR}
\tilde{P}_n(x):= P_m(x^2), \qquad m=\displaystyle\frac{n}{2}, \quad n\in 2{\mathbb{N}}.
\end{equation}
Then $ \tilde{P}_n $ is a  multiple orthogonal polynomials satisfying conditions of the form \eqref{eq:NikPs} with respect to what can be seen as a natural generalization of a  Nikishin system: now the first generating measure $\sigma_1$ is supported on  the whole real axis ${\mathbb{R}}$, while the second generating measure $\sigma_2$ is a discrete measure on the imaginary axis. Then for the rescaled polynomials 
$
\tilde{Q}_n(x):= P_n(d_n x^2)/d_n^{2n}
$
we have straightforward analogues of  Theorems~\ref{teo:equil} and \ref{teo:weak}, but now in terms of the solution of the following equilibrium problem:
there exists a unique pair of measures $(\lambda_1, \lambda_2)$, $|\lambda_1|=2$, $|\lambda_2| = 1, \lambda_2(x) \leq \tilde{\sigma},$  and unique constants $w_1,w_2,$ such that
\begin{equation}\label{EqCo1}
 2 U^{\lambda_1}(x) - U^{\lambda_2}(x) + \tilde{\varphi}(x)\,
\left\{
\begin{array}{ll}
= w_1, & x \in \supp(\lambda_1) \subset {\mathbb{R}}, \\
\geq w_1, & x \in {\mathbb{R}},
\end{array}
\right.
\end{equation}
\begin{equation}\label{EqCo2}
 2 U^{\lambda_2}(x) - U^{\lambda_1}(x)
\left\{
\begin{array}{ll}
\leq w_2, & x \in \supp (\lambda_2) = i{\mathbb{R}}, \\
\geq  w_2, & x \in \supp(\sigma - \lambda_2).
\end{array}
\right.
\end{equation}
The external field and the constraint are related to their
analogues in  \eqref{EqCo1p}--\eqref{EqCo2p} by $\tilde{\varphi}(x)
= {\varphi} (x^2)$, $ \tilde{\sigma}^{\prime}(x) = 2x
{\sigma}^{\prime}(x^2).$ We note, that the polynomials $ \tilde{Q}_n(x)
$ are multiple orthogonal with respect to the varying
weights $s'_{j,n}(x):=s'_j(d_{n}x)$
\begin{equation} \label{VarOrth}
\int_{\mathbb{R}} x^{k} \tilde {Q}_{n}(x)\, s'_{j,n}(x)\, dx\, =0 , \qquad k =0,\ldots, m-1,\quad j=1,2.
\end{equation}

\section{Examples of explicit solutions of the equilibrium problem}\label{sec:5}
As we already mentioned in the introduction, in recent years various models from random matrix theory have been reformulated in  terms of multiple orthogonal polynomials corresponding to Nikishin systems of type~\eqref{eq:NikPs}--\eqref{eq:4^*P}. In all of them, the  generated weights are given by entire functions whose  ratio is a meromorphic  function, which can be considered as the Cauchy transform of a discrete measure
$\sigma_2$ as in \eqref{sig2}.

In this section we discuss three  examples of this type of   Nikishin systems for which  explicit  solutions of the associated equilibrium problems stated  in Theorem  \ref{teo:equil} are available. One of them (see subsection~\ref{ssec:5.4} below) is analyzed in more detail, along with a new approach for expressing the density of the equilibrium measure as a jump of the logarithm of an algebraic function. In this representation, the component of the equilibrium measure constrained by the Lebesgue measure is modeled as the jump of  the logarithm of a negative function. In contrast with the standard approach, where either the underlying differential equations or the recurrence relations of the corresponding multiple orthogonal polynomials are used, we derive this representation directly from the equilibrium conditions.

\subsection{Modified Bessel weights (and non-intersecting squared Bessel paths)}\label{ssec:5.1}
In \cite{ElsWal}, \cite{ElsWal2} multiple orthogonal polynomials
$\{P_n\}$ satisfying \eqref{eq:NikPs} for the system of weights
\begin{equation} \label{BesselWeight}
\begin{aligned}
    s'_1(x) & = x^{\nu/2} e^{- \frac{x}{2} } I_{\nu} \left( \sqrt{x} \right),  \\
    s'_2(x) & = x^{(\nu+1)/2}e^{- \frac{x}{2} } I_{\nu+1} \left( \sqrt{x} \right),
    \end{aligned} \qquad x\in {\mathbb{R}}_+\,,
    \end{equation}
where $I_\nu$ is the modified Bessel function, $\nu > -1 $, were introduced and studied. This system  has found
applications (see \cite{MR2470930}, \cite{KMW2}, and  \cite{HK2}) in the description of ensembles of particles following non-intersecting squared Bessel
paths  (i.e. the radial component of the multidimensional Brownian motion \cite{Shir}).  Since this system of multiple orthogonal polynomials is
quite well studied we just briefly notice, that the  polynomials $\{P_n\}$,
rescaled as in \eqref{scaling} have the asymptotic zero distribution given in \eqref{lessweak}.

The ratio of two weights from \eqref{BesselWeight} is a
meromorphic function which has its poles at the squares of the
zeros of the modified Bessel functions, i.e. $t_k$ in \eqref{sig2}
equals
\begin{equation*}
t_k\,:=\, - \, j_{k,\nu+1}^2\,,\quad k\in \mathbb{Z}_{+}\,,
 \end{equation*}
where $j_{k,\nu}$ is the $k$-th zero of the Bessel function
$J_{\nu}$. To apply  Theorem~\ref{teo:weak} we do not need to have
explicit expressions of the mass points $t_k$ and the values of
the masses $\beta_k$ for the measure $\sigma_2$, but we will need
the asymptotics of the zeros of the Bessel function, see
\cite[p.192]{AS}
\begin{equation}\label{BesselZ}
j_{k,\nu} \,=\, \pi\,(k+\frac{\nu}{2}-\frac{1}{4}) + O\left(\frac{1}{k}\right), \qquad
k\rightarrow \infty,
\end{equation}
and for estimating the values of the  masses $\beta_{k}$ we
can use the asymptotics of the modulus $M_\nu$ of the amplitude of
the Bessel function $J_{\nu}=:M_\nu\,\cos \,\theta_\nu$, see
\cite[p.186]{AS}
\begin{equation}\label{BesselMod}
M_{\nu} (x)\,=\, \sqrt{\frac{2}{\pi\,x}}\,\left(1 + O\left(\frac{1}{x^2}\right)\right), \qquad
x\rightarrow +\infty.
\end{equation}
Choosing the scaling coefficient in \eqref{scaling} as   $d_n =
n^2$ for the measure $\sigma_{2,n}$, see \eqref{sigma2scaled}, we
have $\xi_{k,n} = -(j_{k,\nu}/n)^2$. Using  \eqref{BesselZ},
\eqref{BesselMod} and the asymptotic of the modified Bessel
function on the right half plane, see \cite[p.199]{AS}
$$
I_{\nu}(z)\,=\, \frac{e^z}{\sqrt{2\pi z}}\,\left(1+O\left(\frac{1}{|z|}\right)\right)\,,\qquad |\mbox{arg}z|<\frac{\pi}{2},
$$
it is possible to verify that conditions (i)--(iv) of
Section~\ref{sec2} are satisfied with
$$
\rho(x)\sim     \sqrt{|x|},   \qquad \alpha = 1/2,
$$
(here $f \sim g$ means $0<C_{1}<|f/g|<C_{2}<\infty$ where $C_{1}$,
$C_{2}$ do not depend on $x$), and
\begin{equation} \label{FiSi1}
   \varphi (x)= \frac{x}{2} - \sqrt{x} , \quad x > 0,  \qquad \qquad \frac{d\sigma}{dx} = \frac{1}{\pi \sqrt{|x|}}, \quad x < 0 .
    \end{equation}
(regarding (i), it follows from the fact that \eqref{BesselZ} implies
$\lim_{n \rightarrow \infty} (\xi_{k+1,n} - \xi_{k,n})=\pi$, see
proof of \cite[Lemma 4.4]{HK2}). In subsection~\ref{ssec:5.4} below we give more details
verifying in a similar situation some of the limits in conditions
(iii) and (iv).

\smallskip

The rescaled weak asymptotics of the polynomial sequence $\{P_n\}$ is described by means of
the extremal problem solved in  Theorem~\ref{teo:equil}, with the
particular choice of the external field $\varphi$ and the upper
constraint $\sigma$ indicated in \eqref{FiSi1}. We note, that the
example of this subsection and some other relevant examples were
also discussed in \cite{HK2}, providing an insight of why such
vector equilibrium problem should appear.

\smallskip

 An explicit solution of the equilibrium problem \eqref{EqCo1p}--\eqref{EqCo2p} and \eqref{FiSi1} is known (see \cite{MR2470930}, or \cite[p. 1188]{aptkuij}). The measures $\lambda_j$, $j=1,2$, are absolutely continuous with respect to the Lebesgue measure with
densities that can be expressed in terms of solutions of the cubic equation (a.k.a.~the spectral curve)
\begin{equation} \label{eq-5.9}
    H^3\, -\, 2 H^2 \,+\,
    H \,- \,\frac{2}{z}\,=\,0.
\end{equation}
Equation \eqref{eq-5.9} has three solutions, enumerated in  such a way  that
\begin{equation*} 
\begin{aligned}
    H_0(z) & = \frac{2}{z}+ O(z^{-2}), \\
    H_1(z) & = 1-\frac{\sqrt{2}}{z^{1/2}}-\frac{1}{z} + O(z^{-3/2}), \\
    H_2 (z) & = 1+\frac{\sqrt{2}}{z^{1/2}}-\frac{1}{z} + O(z^{-3/2}),
\end{aligned}
\end{equation*}
as $z \to \infty$.
Then, as it was shown in \cite{MR2470930},  $\lambda_1$ and $\lambda_2$ can be written as
\begin{equation} \label{eq-5.11}
    \begin{aligned}
      \lambda_1'(x) & =  \frac{1}{\pi} \Im H_{0,+}(x), \quad x > 0, \\
  \lambda_2'(x) & = \frac{d \sigma}{dx} - \frac{1}{\pi} \Im H_{1,+}(x), \quad x < 0,
    \end{aligned}
\end{equation}
where the $+$ subindices indicate the boundary values from the upper half plane.

\subsection{Multiple Hermite polynomials (and random matrices with an external source)}\label{ssec:5.3}

Another set of multiple orthogonal polynomials was described in
\cite{ABVA}. It turns out that it is more convenient to deal with
the polynomials $\{\tilde {Q}_{n}\}$, defined by \eqref{VarOrth},
with respect to the system of varying weights
\begin{equation*} 
    s'_{j,n}(x) = e^{- n(\frac{1}{2} x^2 - a_j x)}, \qquad x\in {\mathbb{R}}, \quad j =1, \ldots, p.
        \end{equation*}
This system has found applications in the description of ensembles
of non-intersecting Brownian bridges or random matrices with
external source \cite{ABK}, \cite{BK2}. There, for the case   $p=2$
and $a_1=-a_2=a$,  it was proved that the zero counting measures of
the rescaled polynomials $\{\tilde{Q}_n\}$ (corresponding to
$\{\tilde{P}_n\}$) have a weak limit $\lambda$ which can be
described by means of the spectral curve
\begin{equation} \label{HermiteSpCurve}
H^3\, -\, z H^2 \,+\,(2-a^2)
    H \,+ \,{z}a^2\,=\,0.
\end{equation}
This equation is due to Pastur \cite{Pas}. If we enumerate the branches in \eqref{HermiteSpCurve}
so that, as $z \to \infty$,
\begin{equation*} 
\begin{aligned}
    H_0(z) & = z - \frac{2}{z}+O(z^{-2}), \\
    H_1(z) & = a+\frac{1}{z} + O(z^{-2}), \\
    H_2 (z) & = -a+\frac{1}{z} + O(z^{-2}),
\end{aligned}
\end{equation*}
then  $\lambda$ is given by
\begin{equation} \label{EqMeHerm}
    \lambda'(x)  =  \frac{1}{\pi} \Im H_{0,+}(x), \quad x \in {\mathbb{R}}.
\end{equation}
A generalization of Pastur's curve for arbitrary $p$ can be
seen in  \cite{H}.

It was noticed in  \cite{BDK} (see also \cite{ALT}) that the
measure $\lambda$ in \eqref{EqMeHerm} coincides with the component
$\lambda_1$ in the solution of the equilibrium problem
\eqref{EqCo1}--\eqref{EqCo2} corresponding to the external field
$\tilde{\varphi}$ and the constraint $\tilde{\sigma}$ as follows:
 \begin{equation*} 
   \tilde{\varphi}(x)=  \frac{x^2}{2}-a|x| , \quad x \in {\mathbb{R}},  \qquad
   \qquad d\widetilde{\sigma}(z)=\frac{a}{\pi}\,|dz|, \quad z\in i\mathbb{R}.
 \end{equation*}
 Indeed, multiple Hermite polynomials are  orthogonal as well
 with respect to the weights
 $$\begin{aligned}
 \tilde{s}'_{1,n}(x)  := s'_{1,n}(x) + s'_{2,n}(x) = e^{-n \, x^{2}/2} \cosh(nax), \qquad x\in
 {\mathbb{R}},\\
 \tilde{s}'_{2,n}(x)  := s'_{1,n}(x) - s'_{2,n}(x) =  \tanh(nax)\, \tilde{s}'_1(x) ,\quad \qquad x\in
 {\mathbb{R}}.
 \end{aligned}$$
 Since
 $$
\tanh(nax)=\lim_{N \to \infty}\sum_{k=-N}^{N}\left(\frac{1}{n a}\,\,\frac{1}{x+\frac{i\pi}{n a}(k-\frac{1}{2})}\,\,-\,\,
\frac{1}{i\pi\,(k-\frac{1}{2})}\right)\,,
 $$
 then $(\tilde{s}_{1,n},\,\tilde{s}_{2,n})$ is a Nikishin system
 generated by $\tilde{\sigma}_{1,n}:=\tilde{s}_{1,n}$ and the discrete
 measure
 $$
d\tilde{\sigma}_{2,n}\,:=\,\lim_{N \to \infty}\sum_{k=-N}^{N} \, \frac{1}{n a}\,\delta_{\xi_{k,n}},
\qquad\xi_{k,n}:=\frac{i\pi}{n a}(k-\frac{1}{2}).
 $$
It is clear that
$$
\#\{k: \xi_{k,n} \in [-ix,\,ix]\}\,\sim \,\left[\frac{2nax}{\pi}\right]\,,
$$
thus
$$
\frac{1}{n }\, \lim_{N \to \infty}\sum_{k=-N}^{N} \delta_{\xi_{k,n}}\quad{\stackrel{*}{\longrightarrow}}_n\quad d\,\widetilde{\sigma}(z)=\frac{a}{\pi}\,|dz|, \quad z\in i\mathbb{R},
$$
and conditions (ii) -- (iii) of (an analogue on the real line and the imaginary axis of) Theorem 2.2 
are fulfilled. Regarding (iv) one can use that
$$
-\frac{1}{n}\,\log  \tilde{s}'_{1,n}(x)\,=\,\frac{x^{2}}{2}\,
\begin{cases}  -ax\,-\,\frac{1}{n} \, \log \left(1\,+\,e^{-2 n a x} \right),\quad x\geqslant 0\\[12pt]
+ax\,-\,\frac{1}{n} \, \log \left(1\,+\,e^{+2 n a x} \right),\quad x\leqslant 0
  \end{cases}\,,
  $$
  which leads, in particular, to the uniform convergence when $n \rightarrow \infty$
$$
-\frac{1}{n}\,\log  \tilde{s}'_{1,n}(x)\,\,\rightrightarrows\,\,
\tilde{\varphi}(x):=  \frac{x^2}{2}-a|x|,
$$
on compact subsets of $\mathbb{R}$. As to (i)  it can be derived as in the previous example.

 \medskip

 Actually,   \cite{BDK}
 contains a more  general result for the multiple orthogonal polynomials $\{\tilde{Q}_n\}$ given by \eqref{VarOrth},
 corresponding to the system of  varying weights
\begin{equation*} 
    s'_{j,n}(x)  = e^{-n (V(x) - a_j x)}, \qquad x\in {\mathbb{R}}, \quad j =1,  2,\qquad
            \end{equation*}
where $V(x) = \sum_{j=1}^d v_j x^{2j}$ is an even polynomial potential with $v_d > 0$; it was shown that the  zero counting measures of the  scaled polynomials
$\{\tilde{Q}_n\}$  converge (in a weak-* sense) to the first component $\lambda=\lambda_1$  of the solution to the equilibrium problem (26)-(27)
, with the constraint $\tilde{\sigma}$ and the external field $\tilde{\varphi}$  given by
 \begin{equation} \label{FiSi2}
   \tilde{\varphi}(x)= V(x)-a|x| , \quad x \in {\mathbb{R}}, \qquad d\widetilde{\sigma}(z)=\frac{a}{\pi}\,|dz|, \quad z\in i\mathbb{R}.
 \end{equation}
For a detailed proof of the existence and uniqueness of the
solution of this equilibrium problem see \cite{HK}.

\smallskip

Moreover, it was also proved in \cite{BDK} that the equilibrium problem \eqref{EqCo1}--\eqref{EqCo2} with
input data \eqref{FiSi2}  has always a unique solution $(\lambda_1, \lambda_2), |\lambda_1|=2,|\lambda_2| = 1,$ and that the functions
\begin{equation*} 
    \begin{aligned}
      H_0(z) &  = V'(z) - \int \frac{d\lambda_1(s)}{z-s}, \,\,\quad \quad \qquad \qquad z \in {\mathbb{C}} \setminus S(\lambda_1), \\
     H_1(z) &  = \pm a + \int \frac{d\lambda_1 (s)}{z-s} - \int \frac{d\lambda_2(s)}{z-s}, \qquad z \in {\mathbb{C}}\setminus \left(S(\lambda_1)\cup S(\sigma- \lambda_2)\right),\quad \pm \Re z > 0, \\
      H_2(z) &  = \mp a + \int \frac{d\lambda_2(s)}{z-s}, \quad\,\,\qquad\qquad\qquad z \in
      {\mathbb{C}}\setminus  S(\sigma- \lambda_2), \qquad \qquad \qquad\pm \Re z > 0,
      \end{aligned}
      \end{equation*}
are  the three solutions of the equation
\begin{equation} \label{ExtSourSpCurve}
H^3 + p_2(z) H^2 + p_1(z) H + p_0(z) = 0 \end{equation}
with polynomial coefficients,  
whose degrees can be easily determined from the degree of the potential $V$. However, finding the coefficients of these polynomials explicitly in the most general situation is a very difficult problem. In \cite{ALT} (see also \cite{BDK}) this was done for a general even quartic potential,
\[ V(x) = \frac{1}{4} x^4 - \frac{b}{2} x^2 \]
in the cases when the  Riemann surface of \eqref{ExtSourSpCurve} is of genus either 0 or 1. For instance, when the genus is 1 we have from \cite{ALT} that
$$
H^3-(z^3+bz)H^2+z^2H+a^2z^3=0,
$$
where $a$ and $b$ belong to the triangular domain on the $(a, b)$-plane, bounded  by the curves
 \begin{alignat*}{2}
a_m(b)&:=\frac{\sqrt{6b^3-27b-6(b^2-3)^{3/2}}}{9}>0,
&\qquad
b&\in(-2,-\sqrt{3}),
\\
a_M(b)&:=\frac{\sqrt{6b^3-27b+6(b^2-3)^{3/2}}}{9}>0,
&\qquad
b&\in(-\infty,-\sqrt{3}).
\end{alignat*}
and by the $b$-axis ($a=0$).

\subsection{Multiple Pollaczek polynomials }\label{ssec:5.4}

We have come to the main example as discussed at the end of Section~\ref{sec2}.

The sequence of polynomials, studied in  \cite{Sorokin:2009fk}, is defined by the multiple orthogonality conditions \eqref{eq:NikPs} on $\R_+$ with
\begin{equation}
\label{def_ss}
d s_1(x)=   \frac{dx}{\sinh \frac{\pi \sqrt{x}}{2}}, \quad ds_2(x)=  \frac{1}{\cosh \frac{\pi \sqrt{x}}{2}}\, \frac{dx}{\sqrt{x}} = \frac{\tanh \frac{\pi \sqrt{x}}{2}}{\sqrt{x}}ds_1(x).
\end{equation}
Decomposing $\tanh (\pi z/2)/z$ into simple fractions, it is easy to check that
\[ \frac{ \tanh \frac{\pi \sqrt{z}}{2}}{\sqrt{z}} = \frac{4}{\pi}  \sum_{k \geq 0} \frac{1}{z+ (2k+1)^2} = \int\frac{d\sigma_2(x)}{z-x}
\]
where
\[ \sigma_2 =  \frac{4}{\pi} \sum_{k \in {\mathbb{Z}}_+} \delta_{-(2k + 1)^2}
\]
(cf.~\eqref{sig2}). Hence, $(s_1,s_2) = {\mathcal{N}}(\sigma_1,\sigma_2)$ is a Nikishin system generated by $\sigma_1 = s_1$, supported on $\mathbb{R}_+$, and the discrete measure $\sigma_2$ made of equal masses of size $4/\pi$,  whose support is contained in $(-\infty,0)$. In this case, the re-scaling  \eqref{scaling} is done taking $d_n = 4 n^2$. This yields the measure $\sigma_{2,n}$, see \eqref{sigma2scaled}, with $\xi_{k,n} = -((2k + 1)/2n)^2$ and $\beta_k=4/\pi$. It is easy to verify that conditions (i)--(iv) of Section~\ref{sec2} are satisfied with
\begin{equation}
\label{constraints1}
\rho(x) =  \sqrt{|x|}, \qquad d\sigma(x) = {dx}/{2\sqrt{|x|}}, \qquad \varphi(x) = \pi \sqrt{x}, \qquad \alpha = 1/2.
\end{equation}

For example, to derive the expression of $\sigma$, let $T \in (-\infty,0)$, then
\[ \lim_n \frac{1}{n} \int_{[T,0]} d\left(\sum_{k\geq 1} \delta_{k,n}(t)\right) = \lim_{n } \frac{\sharp \{k: (2k + 1)^2 \leq 4n^2 |T|\}}{n} = \sqrt{|T|} = \int_{[T,0]} \frac{|dt|}{2\sqrt{|t|}}.
\]
Since $d\sigma$ has no mass point this is sufficient to prove convergence in the vague topology (for continuous functions with compact support). We wish to underline that the constraint comes purely from the fact that in between two mass points of $\sigma_{2,n}$ there is at most one zero of $Q_{n,2}$. In this property only the positions of the masses of $\sigma_{2,n}$ intervene not their weights; therefore, the constant $4/\pi$ must be discarded.

Regarding \eqref{b} we have
\[ \frac{1}{n}\log (x^{1/2} s_1'(4 n^2 x))^{-1} =   \frac{1}{n} \log \left( \frac{\sinh (\pi n \sqrt{x})}{\sqrt{x}}\right).
\]
At $x=0$ we give this function its limiting value $\log(\pi n)/n$
to make it continuous. For the proof of the uniform convergence we
make the change of variables $\sqrt{x} = y$. Notice that
\[ \frac{1}{n} \log \left( \frac{\sinh (\pi n y)}{y}\right) = \pi y + \frac{1}{n} \log \frac{1 - e^{-2n \pi y}}{2 y}
\]
Obviously, for $y > 0$ the pointwise limit is $\pi y$. On the other hand,
\[ \left(\frac{1 - e^{-2n \pi y}}{2 y}\right)' = \frac{(4n\pi y +2)e^{-2n\pi y} -2}{4y^2} <0, \qquad y >0,
\]
since the numerator equals $0$ at $y=0$ and
\[ ((4n\pi y +2)e^{-2n\pi y} -2)' = -8n^2\pi^2 y e^{-2n\pi y} < 0,\qquad y >0.
\]
Consequently, on any interval $[0,T], T > 0,$ the function
\[ h_n(x) := \frac{1}{n} \log \left( \frac{\sinh (\pi n y)}{y}\right) - \pi y
\]
attains it's maximum and minimum at the extreme points. We have
\[ \lim_{n\to \infty} h_n(0) = \lim_{n\to \infty} \frac{\log (\pi n)}{n} = 0
\]
and from the pointwise limit
\[ \lim_{n\to \infty} h_n(T) = 0.
\]
Therefore, the uniform convergence follows.

Obviously, a pair of measures $(fds_1,fds_2),$ where $f$ is any
continuous function  such that $0 < c_1 \leq f(x) \leq c_2 < +
\infty, x \in \R_+$, has associated the same vector equilibrium
problem. Thus, the corresponding multiple orthogonal polynomials
exhibit the same rescaled normalized zero distribution as those
corresponding to \eqref{def_ss}. Other examples may be constructed
replacing the discrete component of the Nikishin system by a
Meixner or a Charlier type measure (see, for example, \cite{KW},
\cite{Sor2} or \cite{AA}). A large class, depending on two
parameters, of Meixner-Pollaczek type multiple orthogonal
polynomials was studied in \cite{BenDK} and\cite{BDK} for which the
rescaled logarithmic and ratio asymptotic were given. Our example is
a confluent case of those analyzed in \cite{BenDK}, \cite{BDK}.

We will also consider the corresponding polynomials transplanted to the whole real axis,   for multi-indices of the form $(n,n)$. Using the transformation \eqref{MapR} we obtain a sequence of monic polynomials $\tilde P_{n}$ of degree $2n$, satisfying  the orthogonality relations
\begin{equation} \label{eq:1}
\int_{\mathbb{R}} x^{\nu} \tilde P_{n}(x) \frac{x dx}{ \sinh \pi x} =0 , \qquad \nu =0,\ldots, n-1,
\end{equation}
\begin{equation} \label{eq:2}
\int_{\mathbb{R}} x^{\nu} \tilde P_{n}(x) \frac{dx}{\cosh \pi x}  =0 , \qquad \nu =0,\ldots, n-1,
\end{equation}
that are known as multiple (or generalized) Pollaczek
polynomials (see  \cite{Sorokin:2009fk}). In order to guarantee normality,  we
will assume additionally that the $n$ are even. In this case, the zeros of $\tilde P_{n}$
are real and simple.

In a similar fashion as it is done for Nikishin systems (on the real line) it
can be deduced that there exists a monic polynomial $\tilde P_{n,2}$,
$\deg \tilde P_{n,2} = n,$ whose zeros are also simple and contained in
$i{\mathbb{R}} \setminus \{0\}$, such that
\begin{equation} \label{eq:3}
\int_{ \mathbb{R}} x^{\nu} \frac{\tilde P_{n}(x)}{\tilde P_{n,2}(x)}\frac{xdx}{\sinh(\pi x)}   =0 , \qquad \nu =0,\ldots, 2n -1,
\end{equation}
and
\begin{equation} \label{eq:4}
\int_{\mathbb{R}} t^{\nu} \frac{\tilde P_{n,2}(t)}{\tilde P_{n}(t)}\int_{i\mathbb{R}} \frac{\tilde P_{n}^2(x)}{\tilde P_{n,2}(x)} \frac{x d x}{(x-t)\sinh(\pi x)}d \beta(t)  =0 , \qquad \nu =0,\ldots,  n  -1,
\end{equation}
where $\beta$ is a discrete measure supported on the imaginary line. Set
\[ \tilde{Q}_{n} (z) =  \tilde P_{n}(n z)/n^{2n}, \qquad \tilde{Q}_{n,2} (z) = \tilde P_{n,2}(n z)/n^n.
\]
The logarithmic (weak) asymptotic behavior of these polynomials was
studied by V.~N.~Sorokin in \cite{Sorokin:2009fk}. Sorokin's
approach is based on the existence of an explicit expression of the generating
function for the polynomials $\tilde Q_n(x)$, to which a weak form of the Darboux method can be applied. On this path,  the weak asymptotics of the polynomials can be deduced from the singularities of the
generating function.

By \eqref{constraints1}, the zero counting measures of the  scaled polynomials
$\{\tilde{Q}_n\}$ (corresponding to $\{\tilde{P}_n\}$) have a weak
limit $\lambda$, which  is  the first component ($\lambda=\lambda_1$) of the solution to
the equilibrium problem \eqref{EqCo1}--\eqref{EqCo2}, with
\begin{equation} \label{FiSi3}
   \tilde{\varphi}(x)=\pi |x|, \quad x \in {\mathbb{R}},  \qquad \qquad d\widetilde{\sigma}(z)= |dz|\quad\text{on}\,\,\, i\mathbb{R}.
 \end{equation}
One of the goals of this section is to obtain $\lambda$  by  a
direct solution of this equilibrium problem.


From electrostatic considerations we expect that $\supp(\lambda_2) =
i{\mathbb{R}}$, because the external field created by
${U}^{\lambda_1}$ on $i{\mathbb{R}}$  is too weak to make $\supp(\lambda_2)$
compact. An alternative argument is that, if there were no
restrictions on $\lambda_2$, the measure $2\lambda_2$ in \eqref{EqCo2} would
coincide with the balayage of $\lambda_1$
onto $i{\mathbb{R}}$. Hence, the upper constraint forces the balayage
measure to redistribute its mass precisely where it exceeds $\sigma$ in order to attain equilibrium on the rest of $i{\mathbb{R}}$. This
consideration makes us look for a solution $\lambda_2$ for which there is an equality on $\supp(\sigma -\lambda_2)$ in
the equilibrium conditions \eqref{EqCo2}.

We shall try to find the Cauchy transform of the equilibrium measure
$\lambda_1$,
\begin{equation} \label{eq:25} H(z) := - \widehat{\lambda}_1(z) = \int_{\mathbb{R}} \frac{d\lambda_1(x)}{x-z}.
\end{equation}
If we ``complexify'' the equilibrium relations \eqref{EqCo1}--\eqref{EqCo2} and \eqref{FiSi3}, differentiate them and
take the real parts, we obtain
\[ \Re\left( 2 \widehat{\lambda}_1(x) - \widehat{\lambda}_2(x) \right)  = \left\{
\begin{array}{ccc}
-\pi, & \mbox{on} & {\mathbb{R}}_- \cap \supp(\lambda_1), \\
\pi, & \mbox{on} & {\mathbb{R}}_+ \cap \supp(\lambda_1),
\end{array}
\right.
\]
and
\[ \Re\left( 2 \widehat{\lambda}_2(x) - \widehat{\lambda}_1(x)\right) =
\begin{array}{ccc}
0, & \mbox{on} &  \supp(\sigma - \lambda_2).
\end{array}
\]
Using the Riemann--Schwartz symmetry principle, from the first
relation we deduce that the function $H$ can be continued
analytically from both sides of the cut along ${\mathbb{R}}_- \cap
\supp(\lambda_1)$. Thus, $H$ can be lifted to a Riemann surface, where
\begin{equation} \label{eq:27}
H(z) = \pi + \widehat{\lambda}_1(z) - \widehat{\lambda}_2(z) := H_1(z)
\end{equation}
is considered on the next sheet.
Analogously, $H$ can be continued analytically from both sides of
the cut along ${\mathbb{R}}_+ \cap \supp(\lambda_1)$, so that
\begin{equation} \label{eq:28}
H(z) = -\pi + \widehat{\lambda}_1(z) - \widehat{\lambda}_2(z) := H_2(z)
\end{equation}
is defined on another sheet of the same surface. Let us assume that the complete Riemann surface
${\mathcal{R}}= \{\overline{{\mathcal{R}}^{(j)}}\}_{j=0}^2$,
$\overline{{\mathcal{R}}^{(j)}} =\overline{\mathbb{C}},$ has three
sheets. With appropriate cuts we will have three branches
of $H = \{H_j\}_{j=0}^2$, where $ H_0(z) = - \widehat{\lambda}_1(z)$ is holomorphic in $\overline{\mathbb{C}} \setminus \supp(\lambda_1)
$,
and (\ref{eq:25})--(\ref{eq:28}) give us that, as $z \to \infty$,
\begin{equation} \label{eq:29}
\begin{split}
H_0(z) & =  -\frac{2}{z} + \dots  \\
H_1(z) & =  \pi + \frac{1}{z} + \dots  \\
H_2(z) & =  - \pi + \frac{1}{z} + \dots.
\end{split}
\end{equation}

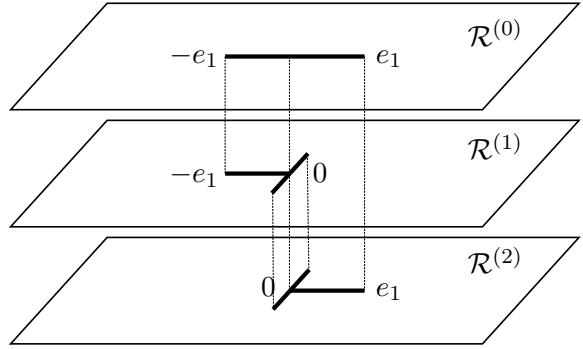
\begin{figure}[htb]
\centering \begin{overpic}[scale=1]%
{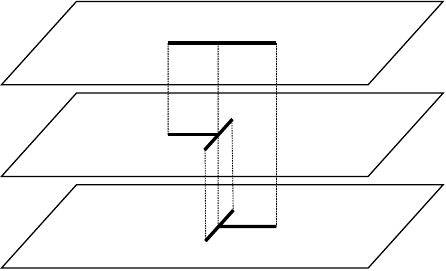}%
 \put(28,50){$-e_1 $}
 \put(64,50){$e_1 $}
 \put(28,29){$-e_1 $}
  \put(64,9){$e_1 $}
   \put(53,29){$0 $}
      \put(44,9){$0 $}
 \put(80,53){$\mathcal R^{(0)} $}
 \put(80,33){$\mathcal R^{(1)} $}
  \put(80,13){$\mathcal R^{(2)} $}
\end{overpic}
\caption{Sheet structure of the Riemann surface ${\mathcal{R}}$.}
\label{fig:RS}
\end{figure}

We make an ansatz that the function $H$ can be found in the form
\begin{equation} \label{eq:30}
H(\zeta) = \frac{2}{i} \log \psi(\zeta) \qquad \mbox{on} \qquad {\mathcal{R}} \setminus \{\zeta \in {\mathcal{R}}: \psi(\zeta) \in {\mathbb{R}}_-\},
\end{equation}
where $\psi$ is a meromorphic function on the compact three sheeted
Riemann surface ${\mathcal{R}}$. At this moment, ${\mathcal{R}}$ is still
unknown (should it exist);  however, representation (\ref{eq:30}) and relations (\ref{eq:29}) yield that
\begin{equation} \label{eq:31}
\psi(\zeta) = \begin{cases}
1 - \frac{i}{\zeta} + \ldots,  & \zeta \to \infty^{(0)},\\[2mm]
i - \frac{1}{2\zeta} + \ldots,  & \zeta \to \infty^{(1)},\\[2mm]
-i + \frac{1}{2\zeta} + \ldots, & \zeta  \to \infty^{(2)},
\end{cases}
\end{equation}
where $q^{(j)}$ denotes the point on $\mathcal R^{(j)}$ whose canonical projection on the plane is $q\in \overline \C$.
We try to take $\psi$ as the simplest meromorphic function which
maps  ${\mathcal{R}}$ conformally   onto $\overline{\mathbb{C}}$.
The inverse of this function is a rational function $\zeta =
r(\psi)$. From the main term in the asymptotic expansion
(\ref{eq:31}) we have that
\[ \zeta = \frac{A}{\psi -1} + \frac{B}{\psi - i}+ \frac{C}{\psi + i},
\]
and the second term gives us that
\[ A= -i, \qquad B = \frac{-1}{2} ,\qquad C = \frac{1}{2}.\]
Thus,
\begin{equation} \label{eq:40} \zeta = -i \frac{\psi(\psi +1)}{(\psi^2 +1)(\psi -1)}
\end{equation}
or, what is the same,
\begin{equation} \label{eq:32} \psi^3 + \frac{i-\zeta}{\zeta} \psi^2 + \frac{i+\zeta}{\zeta} \psi -1 = 0.
\end{equation}
The discriminant of (\ref{eq:32}) is equal to
\[ 16 \zeta^4 - 44\zeta^2 -1.
\]
Therefore, the algebraic function has four branch points $\pm e_1$
and $\pm e_2$, where
\[ e_1 = \frac{1}{4} \sqrt{22- 10\sqrt{5}}, \qquad e_2 = \frac{i}{4}\sqrt{-22 + 10\sqrt{5}}.
\]
Taking into account (\ref{eq:31}) we fix the following sheet structure of ${\mathcal{R}}$ (see Figure~\ref{fig:RS})
\begin{equation} \label{eq:33}
\begin{split}
{\mathcal{R}}^{(0)} := \overline{\mathbb{C}} \setminus & [-e_1,e_1], \qquad {\mathcal{R}}^{(1)} := \overline{\mathbb{C}} \setminus ([-e_1,0] \cup [-e_2,e_2]),\\ & {\mathcal{R}}^{(2)} := \overline{\mathbb{C}} \setminus ([0,e_1]\cup[-e_2,e_2]).
\end{split}
\end{equation}

Therefore, the algebraic function $\psi$ has the following
single-valued meromorphic branches (in fact holomorphic, since
$\psi(0) = \{0,-1,\infty\}$):
\[ \psi_0(\zeta) \in {\mathcal{H}}(\overline{\mathbb{C}} \setminus [-e_1,e_1]), \qquad \psi_1(\zeta) \in {\mathcal{H}}(\overline{\mathbb{C}} \setminus ([-e_1,0] \cup [-e_2,e_2])),
\]
\[ \psi_2(\zeta) \in {\mathcal{H}}(\overline{\mathbb{C}} \setminus ([0,e_1]\cup[-e_2.e_2])),
\]
where $\mathcal H(\Omega)$ stands for the class of functions holomorphic (and single-valued) in a domain $\Omega$.
From the analysis of the roots of (\ref{eq:32}) it follows that
\begin{equation} \label{eq:34}
 \begin{split}
 &\{i{\mathbb{R}}\}^{(0)} = \{\zeta \in {\mathcal{R}}: \psi(\zeta) \in {\mathbb{R}}_+\}, \\ &\{[-e_2,e_2]\}^{(1)} \cup \{[-e_2,e_2]\}^{(2)} = \{\zeta \in {\mathcal{R}: \psi(\zeta) \in {\mathbb{R}}_-}\}.
\end{split}
\end{equation}
Thus, if we cut our compact Riemann surface ${\mathcal{R}}$ along
the second set in (\ref{eq:34}) and denote
\begin{equation} \label{eq:35}
\widetilde{\mathcal{R}}:={\mathcal{R}} \setminus (\{[-e_2,e_2]\}^{(1)} \cup \{[-e_2,e_2]\}^{(2)}),
\end{equation}
we get that the function $H$ in (\ref{eq:30}) is single-valued and
holomorphic in the open Riemann surface $\widetilde{\mathcal{R}}$.
Now, we can formulate our result about the solution of the equilibrium
problem \eqref{EqCo1}--\eqref{EqCo2}:
\begin{proposition} \label{pro:1}  Let
$$
H_j(\zeta) = \frac{2}{i} \log \psi_j(\zeta), \quad \zeta\in \mathcal R^{(j)}, \quad j=0, 1,
$$
where the $\psi_j$ are the solutions of \eqref{eq:32} satisfying \eqref{eq:31}.
Define the absolutely continuous measures
\[ d\lambda_1(x) = \lambda_1'(x)dx,  \qquad d\lambda_2(x) = \lambda_2'(x)|dx|,
\]
by
\begin{equation} \label{eq:36}
\begin{array}{clll}
\lambda_1'(x) & = & \displaystyle{\frac{1}{\pi} \lim_{\varepsilon \to 0+} |\Im H_0(x + i\varepsilon)}|, & x \in {\mathbb{R}}, \\
\lambda_2'(x) & = & \displaystyle{-1 + \frac{1}{\pi} \lim_{\varepsilon \to 0+} \Re H_1(x - \varepsilon)}, & x \in i{\mathbb{R}} = \supp(\lambda_2).
\end{array}
\end{equation}
The pair  $(\lambda_1,\lambda_2)$
is the solution of the equilibrium problem \eqref{EqCo1}--\eqref{EqCo2} and \eqref{FiSi3}. More precisely, $|\lambda_1| = 2$, $|\lambda_2|= 1$,
and these measures  verify
\begin{equation} \label{eq:37}
d\sigma(x) = |dx|, \quad \lambda_2 \leq \sigma,  \quad \text{and} \quad \lambda_2'(x) = 1 \text{ for } x \in [-e_2, e_2];
\end{equation}
\begin{equation} \label{eq:46}
2 {U}^{\lambda_1}(x) - {U}^{\lambda_2}(x) + \pi |x|
\left\{
\begin{array}{ll}
= w_1, & x \in [-e_1, e_1] = \supp(\lambda_1) \subset {\mathbb{R}}, \\
> w_1, & x \in {\mathbb{R}} \setminus [-e_1,e_1],
\end{array}
\right.
\end{equation}
and
\begin{equation} \label{eq:47}
2 {U}^{\lambda_2}(x) - {U}^{\lambda_1}(x)
\left\{
\begin{array}{ll}
=  w_2, & x \in \supp(\sigma - \lambda_2) = i{\mathbb{R}} \setminus (-e_2,e_2), \\
< w_2, & x \in (-e_2,e_2).
\end{array}
\right.
\end{equation}
\end{proposition}

Before proving Proposition \ref{pro:1} we discuss some properties
of the primitive function $G$ defined by
\begin{equation} \label{eq:26} G^{\prime} = H,
\end{equation}
which we  now consider on the open Riemann surface,
$\widetilde{\mathcal{R}}$. That is
\begin{equation} \label{eq:39}
G(\zeta) = \int_{\zeta_0}^{\zeta} H(t) dt, \qquad \zeta_0, \zeta, t \in \widetilde{\mathcal{R}}.
\end{equation}
The uniformization of $\mathcal{R}$ defined in (\ref{eq:40})  allows
us to integrate by parts obtaining
\begin{equation} \label{eq:41}
 \begin{split}
G(\zeta) &= -2 \int_{\psi(\zeta_0)}^{\psi(\zeta)} \log (\psi)\,\, d \frac{\psi(\psi +1)}{(\psi^2 +1)(\psi -1)}\\& = C + \zeta H(\zeta)
+ 2\log(\psi(\zeta) -1) - \log(\psi^2(\zeta) +1),
\end{split}
\end{equation}
where $C$ is a constant which depends on $\zeta_0$. According to
(\ref{eq:41}), $G$ is multivalued on $\widetilde{\mathcal{R}}$ and
has local analytic extension to the whole ${\mathcal{R}}$ (and beyond),
with possible singular points at $\zeta =0$ and $\zeta =\infty$
(notice that by \eqref{eq:40}, $\psi(\infty) = \{1, i, -i\}$).
However, its periods are purely imaginary. Therefore, its real part
is a single valued harmonic function on ${\mathcal{R}}\setminus
\{0,\infty\}$,
\[ g := \{g_j = \Re G_j\}_{j=0}^2,
\]
which is defined up to an additive constant. We fix the constant so
that
\[ g_0(\infty) + g_1(\infty) + g_3(\infty) = 0.
\]
This normalization  in turn implies that
\begin{equation} \label{eq:42} g_0(\zeta) +g_1(\zeta) +g_2(\zeta) \equiv 0,\qquad \zeta \in {\mathbb{C}}.
\end{equation}
Indeed, $g_0 + g_1 +g_2$ is a symmetric function of $g$ which is
harmonic on $\overline{\mathbb{C}} \setminus \{0,\infty\}$. From
(\ref{eq:40}) and (\ref{eq:41}), one sees that the singularity it
has at  $\zeta =0$ is removable. On the other hand, from
(\ref{eq:29}) and (\ref{eq:39}), we have that the branches of $g$
at infinity have the following behavior
\begin{equation} \label{eq:43}
g(\zeta)  \simeq \begin{cases}
-2 \log|\zeta|, & \zeta \to \infty^{(0)}, \\
\pi \Re z + \log|\zeta|, & \zeta \to \infty^{(1)}, \\
-\pi \Re z + \log|\zeta|, & \zeta \to \infty^{(2)}.
\end{cases}
\end{equation}
So, $\zeta = \infty$ is also a removable singularity of $g_0 + g_1
+g_2$. Since $g_0 + g_1 +g_2$ is harmonic in
$\overline{\mathbb{C}}$ and equal to zero at $\infty$, it is
identically equal to zero.

\vspace{0,2cm}
{\bf Proof of Proposition \ref{pro:1}.} We must verify that the measures defined by their densities  in (\ref{eq:36}) verify (\ref{eq:37})--(\ref{eq:47}). In order to identify the potentials of the measures $\lambda_1,\lambda_2$, let us change the sheet structure of $\mathcal{R}$. Define 
\begin{equation} \label{eq:48}
g_0^*:= g_0, \quad g_1^*:=
\begin{cases}
g_1(z), & \Re z < 0, \\
g_2(z), & \Re z > 0,
\end{cases}
\quad g_2^*:=
\begin{cases}
g_2(z), & \Re z < 0, \\
g_1(z), & \Re z > 0.
\end{cases}
\end{equation}
On $i{\mathbb{R}}$, $g^*$ is defined by continuity. Notice that now $g_1^*$, $g_2^*$ have a harmonic continuation through the interval $[-e_2,e_2]$. 

Now, we see that the function $g_0^*$ is superharmonic,  and that
$g_2^*$ is subharmonic (being  the maximum of two harmonic functions).
Therefore, taking into account the behavior at $\infty$ (see
(\ref{eq:43})), from the Riesz decomposition theorem for
superharmonic functions we obtain a global representation of the
branches of $g^*$ in ${\mathbb{C}}$ in the form
\begin{equation} \label{eq:49}
\begin{array}{lll}
g_0^*(z) & = & {U}^{\lambda_1}(z) + \kappa_1, \\
g_2^*(z) & = & -{U}^{\lambda_2}(z) -  v(z) + \kappa_2,
\end{array}
\end{equation}
where $\lambda_1, \lambda_2$ are measures supported on $[-e_1,e_1]$
and $i{\mathbb{R}}$, respectively, and $v(z)$ is the superharmonic
function
\begin{equation} \label{eq:50}
v(z) = \left\{
\begin{array}{rr}
\pi \Re z, &  \Re z  \leq 0, \\
-\pi \Re z,&  \Re z  > 0.
\end{array}
\right.
\end{equation}
As a consequence of (\ref{eq:42}), we also have that
\begin{equation} \label{eq:51}
g_1^*(z) = -{U}^{\lambda_1}(z) + {U}^{\lambda_2}(z) + v(z) - \kappa_1 -\kappa_2.
\end{equation}
Using (\ref{eq:29}) and (\ref{eq:49}), it is easy to verify that
\[ |\lambda_1| = 2, \qquad    |\lambda_2| =1,
\]
and taking into consideration the definition of $g$, the
Stieltjes-Perron formula applied to the calculation of the measures
yields (\ref{eq:36}).

Since $g_0^*(x) = g_1^*(x)$ for $ x \in [-e_1,e_1]$, using (\ref{eq:49})
and (\ref{eq:51}) we obtain the equality in (\ref{eq:46}) with
$w_1 := -2\kappa_1 -\kappa_2$. The fact that $g_0^*(x) >
g_1^*(x)$ on $ {\mathbb{R}} \setminus [-e_1,e_1]$ allows us to verify the
inequality in (\ref{eq:46}). Analogously, comparing $g_1^*$ and
$g_2^*$ on $i{\mathbb{R}}$, and using (\ref{eq:49}), (\ref{eq:51}) and
the fact that $v(z) \equiv 0, z \in i{\mathbb{R}}$ (see
(\ref{eq:50})), we obtain (\ref{eq:47}) with $w_2 := 2\kappa_1
+ \kappa_2$.

Finally, notice that the functions $\psi_1,\psi_2$ have negative
limiting values on $[-e_2,e_2]$ (see the second  relation in
(\ref{eq:34})). Therefore, taking into consideration (\ref{eq:30}),
it follows that
$$
\lim_{\varepsilon \to 0+} \Re H_1(x - \varepsilon)
= 2\pi, \quad x \in [-e_2.e_2],
$$
and $\lambda_2'(x) \equiv 1$, $x \in
[-e_2, e_2]$. On the rest of the imaginary axis,
$$
\pi <
\lim_{\varepsilon \to 0+} \Re H_1(x - \varepsilon) < 2\pi
$$ (see
also (\ref{eq:29})). Thus, we obtain (\ref{eq:37}). We wish to
remark that when applying the Stieltjes-Perron formula in the
second half of (\ref{eq:36}) we take the imaginary part because
$|dx| = -idx, x \in  i{\mathbb{R}}$.  This concludes the proof. \hfill $\Box$

\section{Scalar case}

\subsection{Potentials of measures with unbounded support}

In all that follows, finite positive Borel measures $\mu$ supported in $\R$ which verify
\begin{equation} \label{cond4}
\int \log (1 +y^2) d\mu(y) < +\infty.
\end{equation}
play  a central role. It is easy to see that \eqref{cond4} is equivalent to $\int \log (1 +|y|) d\mu(y) < +\infty$ or $\int_{|y|\geq 1} \log |y| d\mu(y) < +\infty$.

Another  important assumption on a measure $\mu$ which we will  use is that for every $\varepsilon > 0$ there exists $0 < \delta < 1/2$ and $R_0 >0$ such that
\begin{equation}
\label{cond5}
\sup_{|R| \geq R_0} \int_{R-\delta}^{R +\delta} \log \frac{1}{|R - y|} d\mu(y) < \varepsilon.
\end{equation}
Obviously, if $\mu \leq \mu^*$ and $\mu^*$ verifies \eqref{cond5} then  $\mu$ verifies \eqref{cond5}. In particular, a  sufficient condition is that there exists $R_0 > 0$ such that $d \mu|_{\R \setminus (-R_0,R_0)} \leq |f| dm $, where $f \in L_{\infty}(m)$ and $m$ is the Lebesgue measure.

We have
\begin{lemma}
\label{convcero}
Let $\mu$ be a finite positive Borel measure on $\R_+$ such that $U^{\mu}$ is continuous at some point $x_0 \in \supp(\mu)$, then for every compact $K\subset \C$ and every $\varepsilon > 0$ there exists $\delta > 0$ such that
\begin{equation} \label{unif1}
\sup_{x \in K} \int_{x_0-\delta}^{x_0 +\delta}  \left| \log  |x-y|\right| d\mu(y) < \varepsilon.
\end{equation}
Suppose that \eqref{cond4}-\eqref{cond5} take place. Then, for every $\varepsilon >0$ there exists $R_0$   such that
\begin{equation} \label{unif}
\sup_{R \geq R_0} \sup_{x \in [0,R]} \int_R^{+\infty}\left| \log  |x-y|\right| d\mu(y) < \varepsilon
\end{equation}
and
\begin{equation} \label{unif2}
 \lim_{x \to \infty}   \int \left| \log  {\left|1-\frac{y}{x}\right|}  \right| d\mu(y)= 0,
\end{equation}
where $x \to \infty$ in any direction in $\C$.
\end{lemma}

{\bf Proof.} Let us prove \eqref{unif1}.
Consider the closed disk $B = \{x: |x-x_0| \leq 1/2\}$. For all $x\in B$
\[ 0 < \int_B |\log |x-y|| d\mu(y) = \int_B \log \frac{1}{|x-y|} d\mu(y).
\]
Obviously, $U^{\mu|_B}$ is  continuous at $x_0$. Therefore, $\log (1/|x-x_0|) \in L_1(\mu|_B)$ and $x_0$ is not a mass point of $\mu|_B$. Consequently, for every $\varepsilon > 0$ there exists $0 < \delta_1 < 1/2$ such
\[0 <  \int_{x_0 - \delta_1}^{x_0+\delta_1} \log \frac{1}{|x_0-y|} d\mu(y)  < \varepsilon/2. \]
The potential of the measure $\mu|_{[x_0 - \delta_1,x_0 +\delta_1]}$ is also continuous at $x_0$, so there exists $0 < \delta_2 < 1/2$ such that
\[\left|\int_{x_0 - \delta_1}^{x_0+\delta_1} \log \frac{1}{|x-y|} d\mu(y)-  \int_{x_0 - \delta_1}^{x_0+\delta_1} \log \frac{1}{|x_0-y|} d\mu(y)\right| < \varepsilon/2, \qquad |x - x_0| < \delta_2.\]
Using these two inequalities we obtain
\[0 < \int_{x_0 - \delta_1}^{x_0+\delta_1} \log \frac{1}{|x-y|} d\mu(y)  < \varepsilon , \qquad |x - x_0| < \delta_2.\]

Fix a compact set $K \subset \C$ and take $K_1 = K \setminus \{x :  |x - x_0| < \delta_2\}$. Since the distance from $K_1$ to $x_0$ is positive and $x_0$ is not a mass point of $\mu|_{[x_0 - \delta_1,x_0 +\delta_1]}$ there exists $0 < \delta_3 < \delta_1$ such that
\[ \int_{x_0 - \delta_3}^{x_0+\delta_3} \left| \log  {|x-y|}\right| d\mu(y)  < \varepsilon , \qquad x \in K_1.\]
On the other hand,
\[0 < \int_{x_0 - \delta_3}^{x_0+\delta_3} \log \frac{1}{|x-y|} d\mu(y) \leq \int_{x_0 - \delta_1}^{x_0+\delta_1} \log \frac{1}{|x-y|} d\mu(y)  < \varepsilon, \qquad |x - x_0| < \delta_2.\]
The last two relations imply \eqref{unif1}.

If $\mu$ has compact support assertions \eqref{unif} and \eqref{unif2} are trivial so in their proof we restrict our attention to measures with unbounded support in $\R_+$. We will analyze \eqref{unif} by sections.  Take $R >1$.

Assume that $x \in [0,R-1]$, then $y -x \geq 1$  for all  $y \in [R,+\infty)$. Using the monotonicity of the logarithm and \eqref{cond4}, we obtain
\[0 \leq \lim_{R\to +\infty}\sup_{x \in [0,R-1]} \int_{R}^{+\infty}\left| \log \frac{1}{|x-y|}\right|\, d\mu(y) = \lim_{R\to +\infty} \sup_{x \in [0,R-1]} \int_{R}^{+\infty} \log (y -x)\,  d\mu(y)
\]
\[ \leq \, \lim_{R\to +\infty} \int_{R}^{+\infty} \log (y) \, d\mu(y) = 0.
\]
By the same token
\[ \lim_{R\to +\infty}\sup_{x \in [0,R]}\int_{R+1}^{+\infty}  \left|\log \frac{1}{|x-y|}\right|\, d\mu(y) =  0
\]

Choose a constant $\delta, 0 < \delta <1/2$. For $x \in [R-1,R-\delta]$ and $y \in [R,R+1]$
\[\log \frac{1}{2} \leq \log \frac{1}{|x-y|} \leq \log \frac{1}{\delta},\]
which implies that
\[\left|\log \frac{1}{|x-y|}\right| \leq \log \frac{1}{\delta}.\]
Consequently
\[ 0 \leq \lim_{R\to +\infty}\sup_{x \in [R-1,R-\delta]} \int_{R}^{R+1} \left| \log \frac{1}{|x-y|}\right|\, d\mu(y) \leq \log \frac{1}{\delta} \lim_{R\to +\infty}\mu([R,R+1])  = 0
\]
since $\mu$ is finite. Analogously,
\[ \lim_{R\to +\infty}\sup_{x \in [R-1,R]} \int_{R+\delta}^{R+1}\left| \log \frac{1}{|x-y|}\right|\, d\mu(y)   = 0.\]

Fix $\varepsilon >0$, from \eqref{cond5} there existe $0<\delta <1/2$ and $R_0 > 0$ such that
\[ \sup_{R \geq R_0} \sup_{x \in [R-\delta,R]} \int_{R}^{R+\delta} \left|\log \frac{1}{|x-y|}\right|\, d\mu(y) \leq  \sup_{R \geq R_0} \sup_{x \in [R-\delta,R]}  \int_{R}^{R+\delta} \log \frac{1}{y-R}\, d\mu(y) < \varepsilon.
\]
Putting everything together \eqref{unif} follows immediately.

To prove \eqref{unif2} first let us restrict to the limiting case when $x \in \R_+$ and without loss of generality we can assume that $x > 2$. For the moment fix $x$. As a function of $y$ on $\R_+$, the non negative function $\left| \log  {\left|1-\frac{y}{x}\right|}  \right|$  has a vertical asymptote at $y = x$ and zeros at $y \in \{0,2x\}$. It is convex in $[0,x)$ and $(x,2x]$ and concave in $[2x,+\infty)$. The functions $\log(1+y)$ and $\log (y -1)$ are concave in their domain of definition. On the interval $[0,x]$, it is easy to verify that $\left| \log  {\left|1-\frac{y}{x}\right|}  \right| = \log (1+y)$  if and only if $y = 0$ or $y= x-1$. On the interval $[x,2x]$, $\left| \log  {\left|1-\frac{y}{x}\right|}  \right| = \log (y-1)$  if and only if $y = x+1$. Taking account of the concavity properties of the functions in the specified intervals and the monotonicity of the logarithm it follows that
\begin{equation}
\label{desig3}
 \left| \log  {\left|1-\frac{y}{x}\right|}  \right|
\left\{
\begin{array}{ll}
\leq  \log(1+y), & 0 \leq y \leq x-1, \\
\leq \log(y-1) \leq \log(1+y), & x+1 \leq y \leq 2x ,\\
= \log\left(\frac{y}{x} -1 \right) \leq \log(1+y), & 2x \leq y < +\infty.
\end{array}
\right.
\end{equation}

Denote $E_x = [x^{\alpha}/2,+\infty) \setminus (x-1,x+1), 0 < \alpha < 1$. Fix $\varepsilon > 0$ and take $0 <\delta < 1/2$, such that \eqref{cond5} takes place. We have
\[ 0 \leq \int \left| \log  {\left|1-\frac{y}{x}\right|} \right|d\mu(y) \leq \int_0^{x^{\alpha}/2} + \int_{E_x} + \int_{x-1}^{x+1} \left| \log  {\left|1-\frac{y}{x}\right|} \right|d\mu(y).\]
Let us analyze these integrals separately.

First
\[0 \leq \int_0^{x^{\alpha}/2}\left| \log  {\left|1-\frac{y}{x}\right|} \right|d\mu(y) = \int_0^{x^\alpha/2}\left| \log  {\left(1-\frac{y}{x}\right)} \right|d\mu(y) \leq
\]
\[|\mu| \left|\log \left(1-\frac{1}{2x^{1-\alpha}}\right) \right|\leq C \left| \frac{1}{x^{1-\alpha}}\right|\to 0, \qquad x \to +\infty.\]
On $E_x$, taking \eqref{desig3} and \eqref{cond4} into account,
\[ 0 \leq \int_{E_x}\left| \log  {\left|1-\frac{y}{x}\right|} \right|d\mu(y) \leq \int_{y \geq x^{\alpha}/2} \log (1+y) d\mu(y) \to 0, \qquad x \to +\infty.\]
Finally, on $[x-1,x+1]$
\[  0 \leq \int_{x-1}^{x+1}\left| \log  {\left|1-\frac{y}{x}\right|} \right|d\mu(y) \leq \mu([x-1,x+1]) \log(x) + \int_{x-1}^{x+1} \log  \frac{1}{|y-x|} d\mu(y) \leq \]
\[\mu([x-1,x+1])  \left( \log(x) + \log(1/\delta) \right) + \int_{x-\delta}^{x+\delta} \log  \frac{1}{|y-x|} d\mu(y),\]
where the first term tends to zero as $x\to +\infty$, on account of \eqref{cond4}, and the second term is bounded by $\varepsilon$ for all sufficiently large $x$ due to \eqref{cond5}.

Summarizing,  we have
\[ 0 \leq \liminf_{x \to +\infty} \int \left| \log  {\left|1-\frac{y}{x}\right|} \right|d\mu(y) \leq \limsup_{x \to +\infty} \int \left| \log  {\left|1-\frac{y}{x}\right|} \right|d\mu(y) \leq \varepsilon,\]
for each $\varepsilon > 0$. Letting $\varepsilon \to 0$ we obtain \eqref{unif2} for the case when $x \in \R_+$.

Now, take $\theta, 0 < \theta < \pi/2$ and define the region $F_{\theta} =  \C \setminus \{x: |\arg (x)|  \leq \theta\}$. Assume that $x \to \infty, x \in F_{\theta}$, In this case, for all $y \geq 0$ and $x \in F_{\theta}$, we have
$y/x \in F_{\theta}$. Consequently $|1 - (y/x)| \geq |\sin (\theta)|  > 0.$ Therefore, if $|x| \geq 1$,
\[ 0 < |\sin (\theta)| \leq  |1 - (y/x)| \leq 1 + |y/x| \leq 1+y.\]
Thus
\[|\log|1 - (y/x)|| \leq \max\{-\log|\sin (\theta)|, \log (1 +y)\}, \quad x \geq1, y \in \R_+.\]
The function defined by the maximum is in $L_1(\mu)$. By  Lebesgue's dominated convergence theorem it follows that
\[ \lim_{x \to \infty, x \in F_{\theta}} \int \left|\log\left|1 - \frac{y}{x}\right|\right| d\mu (y) = 0.
\]

Denote $a_x = \arg(x)$. Now assume that $x \to \infty, a_x \not\to 0$ and
\[\limsup_{x \to \infty }\int |\log|1 - (y/x)|| d\mu (y) > 0.\]
Then we can find $\theta, 0 < \theta < \pi/2$ sufficiently small and a sequence $x_n \in F_\theta, x_n \to \infty$ such that
\[\limsup_{n \to \infty }\int |\log|1 - (y/x_n)|| d\mu (y) > 0.\]
against what was proved before. Consequently, to prove \eqref{unif2} it remains to show that the assertion is true when $x \to \infty$ and $a_x \to 0$. This case is similar to the one when $x \to \infty, x \in \R_+$ so we focus on the main ingredients.

Without loss of generality we can assume that $|x| \geq 1$ and $\mbox{Re}(x) > 2$ where $\mbox{Re}(x)$ denotes the real part of $x$. Let $ |1 - (y/x)|  \geq 1$. This implies that $y \geq 2\mbox{Re}(x)$. Then
\begin{equation}
\label{desig5} |\log|1 - (y/x)|| = \log|1 - (y/x)| \leq \log(1+y), \qquad y \geq 2\mbox{Re}(x).
\end{equation}
Notice that
\[|1 - (y/x)|^2 = |e^{i a_x} - (y/|x|)|^2 = (\cos(a_x) - (y/|x|))^2 +  \sin^2(a_x)  \geq
\]
\[  (\cos(a_x) - (y/|x|))^2 = \cos^2(a_x)(1 - (y/\mbox{Re}(x)))^2.\]
Consequently, when $0 < |1 - (y/x)|  \leq 1$, that is $0 \leq y \leq 2 \mbox{Re}(x)$,
\[0 \geq \log |1 - (y/x)| \geq \log |\cos(a_x)(1 - (y/\mbox{Re}(x)))| \]
and
\[
|\log |1 - (y/x)|| \leq |\log |\cos(a_x)(1 - (y/\mbox{Re}(x)))|| \leq
\]
\begin{equation}
\label{desig6}
|\log |\cos(a_x)| + |\log|1 - (y/\mbox{Re}(x))||, \qquad 0 \leq y \leq 2 \mbox{Re}(x).
\end{equation}
Analyzing separately $y \!\in \![0,\mbox{Re}(x)-\!1]$, $y \!\in\! [\mbox{Re}(x)+\!1, 2\mbox{Re}(x)]$, and  $y\! \in \! [2\mbox{Re}(x), +\infty]$, reasoning as in the deduction of \eqref{desig3} (with $x$ replaced by $\mbox{Re}(x)$), with the help of  \eqref{desig5}-\eqref{desig6} one obtains
\begin{equation}
\label{desig7}
 \left| \log  {\left|1-\frac{y}{x}\right|}  \right| \leq \log(1+y) + \log(\sec{a_x}), \qquad y \in \R_+ \setminus (\mbox{Re}(x)-1,\mbox{Re}(x) +1).
\end{equation}
In the final part of the proof we take $E_x = [\mbox{Re}(x)^{\alpha}/2,+\infty) \setminus (\mbox{Re}(x)-1,\mbox{Re}(x)+1), 0 < \alpha < 1,$ and proceed as in the case when $x \in \R_+$ observing that
\[ \lim_{x \to \infty, a_x \to 0}\int \log(\sec{a_x}) d\mu(y) = \lim_{x \to \infty, a_x \to 0} \log(\sec{a_x}) = 0.
\]
With this we conclude the proof. \hfill $\Box$

With the aid of \eqref{unif} we prove a version of the principle of domination for measures with unbounded support.

\begin{lemma} \label{prindesc} Suppose that $\mu, \nu$ are finite positive Borel measures supported in $\R_+$ such that $|\mu| = |\nu|, { {I}}(\mu) < \infty$,  and verify \eqref{cond4}\! -- \eqref{cond5}. If $\supp(\mu)$ is unbounded and $\supp(\nu)$ is compact we also suppose that $U^{\nu}$ is continuous at some point $x_0 \in \supp(\nu)$. Assume that for some constant $c \in \R$
\begin{equation}
\label{desigpot}
 U^{\mu}(x) \leq U^{\nu}(x) + c, \qquad \mu \quad \mbox{almost everywhere}.
\end{equation}
Then
\begin{equation}
\label{domination2}
U^{\mu}(x) \leq U^{\nu}(x) + c, \qquad x \in \C.
\end{equation}
\end{lemma}

{\bf Proof.} If the supports of $\mu$ and $\nu$ are compact sets the lemma gives the standard statement of the principle of domination (see, for example,  \cite[Theorem II.3.2]{ST}), so this result is new when at least one of the two measures has unbounded support. We will reduce the proof to the case of measures with compact support. We will analyze in detail the case when the supports of $\mu$ and $\nu$ are both unbounded and then mention how to proceed when one of them is bounded and the other unbounded.

Assume that $\supp(\mu)$ and $\supp(\nu)$ are unbounded. Fix $\varepsilon > 0$. According to \eqref{unif} there exist $R_1(\varepsilon), R_2(\varepsilon)$ such that  $\mu([0,R_1]) = \nu([0,R_2])$ and
\begin{equation}
\label{varep3}
\max \left(\sup_{x \in [0,R_1]} \left|\int_{R_1}^{+\infty} \log \frac{1}{|x-y|} d\mu(y)\right|, \sup_{x \in [0,R_2]} \left|\int_{R_2}^{+\infty} \log \frac{1}{|x-y|} d\nu(y)\right|\right) < \varepsilon.
\end{equation}
One can take $R_1(\varepsilon)$ and $R_2(\varepsilon)$ so that $\lim_{\varepsilon \to 0} R_1(\varepsilon) = +\infty, \lim_{\varepsilon \to 0} R_2(\varepsilon) = +\infty$.

Denote $\mu_1 = \mu_1(\varepsilon) = \mu|_{[0,R_1(\varepsilon)]}$ and $\nu_1 = \nu_1(\varepsilon)= \nu|_{[0,R_2(\varepsilon)]}$. We have $|\mu_1| = |\nu_1|$. Since $\mu_1 \leq \mu$ from \eqref{desigpot} and \eqref{varep3} it follows that
\[
U^{\mu_1}(x) \leq U^{\nu_1}(x) + c +2\varepsilon, \qquad \mu_1 \quad \mbox{almost everywhere}.
\]
Notice that $I(\mu_1) < +\infty$. Using \cite[Theorem II.3.2]{ST} we have
\begin{equation}
\label{domination1} U^{\mu_1}(x) \leq U^{\nu_1}(x) + c +2\varepsilon, \qquad x \in \C.
\end{equation}

Fix an arbitrary compact set $K \subset \C$ and let $M = \sup_{x\in K} |x|$. For all sufficiently large $R$
\[ |\log|x-y|| = \log|x -y| \leq \log(M + y), \qquad y \geq R, \qquad x \in K,
\]
and using \eqref{cond4} it follows that
\[\lim_{\varepsilon \to 0} U^{\mu_1(\varepsilon)} = U^{\mu}, \qquad \lim_{\varepsilon \to 0} U^{\nu_1(\varepsilon)} = U^{\nu}\]
uniformly on $K$. Letting $\varepsilon$ tend to zero, \eqref{domination2} follows from \eqref{domination1}  and we are done.

When only $\supp(\nu)$ is unbounded, we proceed as before to reduce $\nu$ to a measure $\nu_1$ with compact support but we can maintain $\mu$ as it is because the principle of domination for compact sets allows $|\nu_1| \leq |\mu|$ to deduce \eqref{domination1}. If $\supp(\mu)$ is unbounded we take $\mu_1$ as before, but we must reduce $\nu$ so that $|\nu_1| \leq |\mu_1| (< |\mu|)$. In order to achieve this, since $\supp(\mu)$ is a compact set we take away mass from a neighborhood of a point $x_0 \in \supp(\nu)$ where $U^{\nu}$ is continuous and use \eqref{unif1} instead of \eqref{unif}.
\hfill $\Box$

\begin{rem}
\label{remdesc}
Lemmas \ref{convcero} and \ref{prindesc} are valid for measures supported on all $\R$. In fact, Lemma  \ref{prindesc}  will be used in the next section for measures supported on $\R_-$.
\end{rem}

\subsection{Equilibrium measure with constraint and external field} \label{4.2}

This question has been considered by several authors (see, for example, \cite{BKMW}, \cite{DrSa}, \cite{GR1}, \cite{HK}, \cite{KR}, and \cite{Rak2}). Our contribution consists in studying the corresponding variational problem in cases when the equilibrium measure does not have compact support.  We will state the corresponding results for measures supported on $\R_-$ because this is the setting in which they will be needed for the proof of Theorem \ref{teo:equil} but they may be restated for measures supported on $\R$.

In order to deal with measures with unbounded support it is convenient to follow the approach used in \cite{HK}. For arbitrary $\mu_1, \mu_2 \in {\mathcal{M}}^+({\mathbb{R}})$, we define a modified logarithmic potential  and mutual energy as follows
\begin{equation}
\label{modpot}
{\mathcal{U}}^{\mu_1}(x) := \int \log \frac{\sqrt{1 + y^2}}{|x-y|} d\mu_1(y),
\end{equation}
\begin{equation}
\label{modener} {\mathcal{I}}(\mu_1,\mu_2) := \int\int \log \frac{\sqrt{1+x^2}\sqrt{1 + y^2}}{|x-y|} d\mu_1(y) d\mu_2(x).
\end{equation}
The modified energy of $\mu$ is then given by ${\mathcal{I}}(\mu) := {\mathcal{I}}(\mu,\mu)$.
The new kernel is connected with the inverse stereographic projection from the ball in ${\mathbb{R}}^3$ centered at $(0,0,1/2)$ and radius $1/2$ onto the extended complex plane. Therefore,
\begin{equation} \label{desig2} \frac{\sqrt{1+x^2}\sqrt{1 + y^2}}{|x-y|} \geq 1
\end{equation}
(for more details see (2.9)--(2.11) in \cite{HK}). Consequently, the modified potential and the mutual energy are uniformly bounded from below for all $\mu_1, \mu_2 \in {\mathcal{M}}^+({\mathbb{R}})$. When $\mu_1,\mu_2$ have finite energy and verify \eqref{cond4} then
\[ {\mathcal{I}}(\mu_1,\mu_2) = I(\mu_1,\mu_2) + \frac{|\mu_2|}{2} \int \log(1+x^2) d\mu_1(x) + \frac{|\mu_1|}{2} \int \log(1+x^2) d\mu_2(x).
\]

In the sequel, $\sigma$ denotes a positive Borel measure, $\supp (\sigma) = {\R}_-$, $|\sigma| > 1$,  such that $U^{\sigma|_K}$ is continuous on $\C$ for every compact subset $K \subset {\R}_-$. Set
\[ {\mathcal{M}}(\sigma) := \{\mu \in {\mathcal{M}}_1^+(\R_-): \mu \leq \sigma\}, \qquad  \widetilde{\mathcal{M}}(\sigma) := \{ {\mu} \in {\mathcal{M}}(\sigma): \mathcal{I}(\mu) < \infty\}.\]

\begin{lemma}
\label{cont}
For any $\mu \in {\mathcal{M}}(\sigma)$, $\mathcal{U}^{\mu}$ is continuous on $\C$.
\end{lemma}
{\bf Proof.} Take $\mu \in {\mathcal{M}}(\sigma)$. Obviously, $\mathcal{U}^{\mu}$ is continuous on $\C \setminus \supp(\mu)$, so we only have to check the continuity on $\R_-$. Choose $x_0 \in \R_-$. Take a compact set $K \subset \R_-$ that contains $x_0$ in its interior.
Since
\[\mathcal{U}^{\mu} = \mathcal{U}^{\mu|_K} + \mathcal{U}^{\mu - \mu|_K}\]
and $x_0 \not \in \supp(\mu - \mu|_K)$ then $\mathcal{U}^{\mu - \mu|_K}$ is continuous at $x_0$. However, $\mu|_K \leq \sigma|_K$ and $U^{\sigma|_K}$ is continuous on $\C$, so (see \cite[Lemma 5.2]{DrSa}) $U^{\mu|_K}$ and $\mathcal{U}^{\mu|_K}$ are continuous on $\C$, in particular at $x_0$. Thus, $\mathcal{U}^{\mu}$ is continuous at any $x_0 \in \R_-$. \hfill $\Box$

\medskip

Let  $\phi$ be a real valued continuous function on $\R_-$ such that
\begin{equation} \label{admis}
\liminf_{x \to -\infty} \phi^*(x) > -\infty. \qquad \phi^*(x):=\phi(x) - \log(1 +x^2).
\end{equation}
For $\mu \in {\mathcal{M}}^+_1(\R_-)$ define
\[ {\mathcal{W}}^{\mu}(x):= 2 \int \log \frac{\sqrt{1 + x^2} \sqrt{1 + y^2}}{|x-y|}d \mu(y) + \phi^*(x) = 2 \mathcal{U}^{\mu}(x)  + \phi(x),
\]
and
\[ {\mathcal{J}}_{\phi^*}(\mu) := 2 \int \left(\int \log \frac{\sqrt{1 + x^2} \sqrt{1 + y^2}}{|x-y|}d \mu(y) + \phi^*(x)\right)d\mu(x) =
\]
\[ 2 {\mathcal{I}}(\mu) +2 \int \phi^*(x) d\mu(x).
\]
If ${\mathcal{I}}(\mu) = +\infty$ we take ${\mathcal{J}}_{\phi^*}(\mu) = +\infty$.

\medskip

Condition \eqref{admis} guarantees that the energy problem for the functional ${\mathcal{J}}_{\phi^*}(\mu)$ is weakly admissible as defined in \cite[Section 2.1]{HK} and according to \cite[Corollary 2.7]{HK} there exists a unique $\lambda \in \widetilde{\mathcal{M}}(\sigma) $ such that
\begin{equation}
\label{minprob}
{\mathcal{J}}_{\phi^*}(\lambda) = \inf \{{\mathcal{J}}_{\phi^*}(\mu): \mu \in {\mathcal{M}}(\sigma)\}.
\end{equation}
The measure $\lambda$ is said to be extremal.

For $\mu \in \widetilde{\mathcal{M}}(\sigma)$ we also introduce the following characteristic value
\[{\mathcal{F}}_\mu := \max \{C \in \R: \mathcal{W}^{\mu}(x) \geq C\,\, \mbox{holds}\,\,(\sigma -\mu)\,\, \mbox{a.e.}  \}. \]

We have

\begin{teo} \label{lem1scalar} Let $\phi$ satisfy \eqref{admis} and let $\sigma$, $\supp (\sigma) = {\R}_-$, $|\sigma| > 1$, be a positive Borel measure such that $U^{\sigma|_K}$ is continuous on $\C$ for every compact subset $K \subset {\R}_-$. The following statements are equivalent and have the same unique solution:
\begin{itemize}
\item[$(A')$] There  exists $ {\lambda} \in \widetilde{\mathcal{M}}(\sigma)$ which is extremal.
\item[$(B')$] There  exists $ {\lambda} \in \widetilde{\mathcal{M}}(\sigma)$ such that for all $ {\nu} \in \widetilde{\mathcal{M}}(\sigma)$
    \[\int {\mathcal{W}}^{ {\lambda}}   d( {\nu}- {\lambda})  \geq 0 .\]
\item[$(C')$] There  exist  $ {\lambda}  \in \widetilde{\mathcal{M}}(\sigma)$ and a constant ${\mathfrak{w}}  = {\mathfrak{w}} (\sigma,\phi) $  such that
 \[ {\mathcal{W}}^{ {\lambda}}(x) = 2 {{\mathcal{U}}}^{ \lambda }(x) + \phi(x)
\left\{
\begin{array}{ll}
\leq {\mathfrak{w}} , & x \in \supp( \lambda ), \\
\geq  {\mathfrak{w}} , & x \in \supp(\sigma -  \lambda ).
\end{array}
\right.
\]
\end{itemize}
The constant ${\mathfrak{w}}$ is uniquely determined and equals  ${\mathcal{F}}_\lambda$. The extremal measure verifies \eqref{cond4}.
\end{teo}

{\bf Proof.} As mentioned above the existence of a unique extremal measure follows from \cite[Corollary 2.7]{HK}. The equivalence of $(A')$ and $(B')$ follows from the identity
\[ {\mathcal{J}}_{\phi^*}( {\nu}_{\varepsilon}) - {\mathcal{J}}_{\phi^*}( {\lambda}) = \varepsilon^2 {\mathcal{J}}_0( {\nu} -  {\lambda}) + 2\varepsilon \int {{\mathcal{W}}}^{ {\lambda}}  d( {\nu} -  {\lambda}),
\]
valid for all $\lambda, \nu \in \widetilde{\mathcal{M}}(\sigma)$, where $\nu_{\varepsilon} = \varepsilon \nu + (1-\varepsilon)\lambda, 0\leq \varepsilon \leq 1$ and ${\mathcal{J}}_0( {\nu} -  {\lambda})$ is the functional applied to $ {\nu} -  {\lambda}$ with $\phi^* \equiv 0$.

Assume that $ {\lambda}$ is extremal. From the identity it follows that
\[ \varepsilon^2 {\mathcal{J}}_0( {\nu} -  {\lambda}) + 2\varepsilon \int {{\mathcal{W}}}^{ {\lambda}}   d( {\nu} -  {\lambda}) \geq 0.
\]
Dividing by $\varepsilon$ and letting $\varepsilon \to 0$, we have
\begin{equation} \label{eq:24} \int {{\mathcal{W}}}^{ {\lambda}}   d( {\nu} -  {\lambda})\geq 0, \qquad  {\nu} \in \widetilde{\mathcal{M}}(\sigma),
\end{equation}
so $(A')$ implies $(B')$. Taking $\varepsilon = 1$, we get
\[ {\mathcal{J}}_{\phi^*}( {\nu}) - {\mathcal{J}}_{\phi^*}( {\lambda}) =  {\mathcal{J}}_0( {\nu} -  {\lambda}) + 2 \int {{\mathcal{W}}}^{ {\lambda}}   d( {\nu} -  {\lambda}).
\]
From  \cite[Theorem 2.5]{CKL} we have ${\mathcal{J}}_0( {\nu} -  {\lambda}) \geq 0$ with equality if and only if $ {\nu} = {\lambda}$. Therefore, $(B')$ implies $(A')$ and the solution to $(B')$ is unique.

Now, let us prove that any solution to $(C')$ solves $(B')$.  Let $ {\lambda} $  verify $(C')$ and take $ {\nu}  \in \widetilde{\mathcal{M}}(\sigma) $.
Since $ |\lambda|= |\nu| = 1$
\[  \int {{\mathcal{W}}}^{ {\lambda}} \, d(\nu  -  \lambda ) = \int ({{\mathcal{W}}} ^{ {\lambda}} - {\mathfrak{w}} ) \, d(\nu  -  \lambda ).
\]
Define
\[ E_+ = \{ t \in \R_- : {{\mathcal{W}}} ^{ {\lambda}}(t) - {\mathfrak{w}}  > 0\}, \qquad E_- = \{ t \in R_-: {{\mathcal{W}}} ^{ {\lambda}}(t) - {\mathfrak{w}}  < 0\}.
\]
According to   $(C')$, $ \lambda (E_+) =0$, so
\[ \int_{E_+} ({{\mathcal{W}}} ^{ {\lambda}} - {\mathfrak{w}} ) d(\nu  -  \lambda ) = \int_{E_+} ({{\mathcal{W}}} ^{ {\lambda}} - {\mathfrak{w}} ) d\nu  \geq 0.
\]
Additionally, $(\sigma -  \lambda )(E_-) =0$. Take an increasing sequence of compact sets $K_n \subset E_-$ such that
$\lim_{n \to \infty} (\sigma -  \lambda )(K_n) = (\sigma -  \lambda )(E_-)$. By Lemma \ref{cont}, ${{\mathcal{W}}}^{ {\lambda}}$ is continuous on all $\C$, in particular on $K_n$, and therefore ${{\mathcal{W}}} ^{ {\lambda}} - {\mathfrak{w}}$ is bounded on $K_n$. Using Lebesgue's monotone convergence theorem it follows that
\[ \int_{E_-} |{{\mathcal{W}}} ^{ {\lambda}} - {\mathfrak{w}} | d(\sigma - \lambda ) = \lim_{n \to \infty} \int_{E_-} 1_{K_n}|{{\mathcal{W}}} ^{ {\lambda}} - {\mathfrak{w}} | d(\sigma - \lambda ) = 0,
\]
where $1_{K_n}$ is the function which equals $1$ on $K_n$ and $0$ elsewhere. Consequently, taking into account that
$\nu \leq \sigma$, we obtain
\[ \int_{E_-} ({{\mathcal{W}}}^{ {\lambda}} - {\mathfrak{w}}) d(\nu -  \lambda) = \int_{E_-} ({{\mathcal{W}}}^{ {\lambda}} - {\mathfrak{w}}) d(\nu - \sigma) + \int_{E_-} ({{\mathcal{W}}}^{ {\lambda}} - {\mathfrak{w}}) d(\sigma -  \lambda) \geq 0.
\]
Putting these relations together, we obtain
\[ \int {{\mathcal{W}}}^{ {\lambda}}  d( {\nu} -  {\lambda}) \geq 0, \qquad \nu \in \widetilde{\mathcal{M}}(\sigma),
\]
as claimed. Therefore, $(C')$ has a unique solution. Let's see that $(B')$ implies $(C')$.

Suppose that $\lambda$ solves $(B')$ and consider the value
\[  \mathcal{F}_{\lambda} =  \max\{C  \in {\R}: {\mathcal{W}} ^{ {\lambda}} \geq C, \quad (\sigma - \lambda)\,\,\mbox{a.e.}\}.
\]
Suppose that there exists $x_0 \in \supp(\lambda)$ such that ${\mathcal{W}}^{ {\lambda}}(x_0) > \gamma >  \mathcal{F}_{\lambda}$. By the definition of $\mathcal{F}_{\lambda}$, there exists a compact $K_1 \subset \supp(\sigma - \lambda),$ such that ${\mathcal{W}}^{ {\lambda}}(x) < \gamma$, $x \in K_1,$ and $(\sigma - \lambda)(K_1) > 0$. On the other hand, ${\mathcal{W}}^{ {\lambda}}(x)$ is continuous on $\R_-$, so there exists $\delta > 0$ sufficiently small such that ${\mathcal{W}}^{ {\lambda}}(x) > \gamma$ for $ |x -x_0| < \delta,$ and by the same token there exists a compact set $K_2$ with $\lambda(K_2) > 0,$ such that ${\mathcal{W}}^{ {\lambda}}(x) > \gamma$ for $ x\in K_2$. Obviously, $K_1 \cap K_2 = \emptyset$. Choose $\alpha, \beta  \in (0,1)$ such that $\beta(\sigma - \lambda)(K_1) = \alpha\lambda(K_2)$. Define a signed measure $\eta$ equal to $-\alpha\lambda$ on $K_2$, equal to $\beta(\sigma - \lambda)$ on $K_1$, and zero otherwise.

Let us prove that $\nu:=\lambda + \eta \in \widetilde{\mathcal{M}}(\sigma)$.
In fact,
\[0 \leq \nu|_{K_2} = (1 - \alpha)\lambda|_{K_2} \leq \sigma|_{K_2}, \]
\[0 \leq \nu|_{K_1} = \beta \sigma|_{K_1} + (1-\beta)\lambda|_{K_1} \leq \sigma|_{K_1},\]
and since $\supp(\nu) = \supp(\lambda)$, we have
\[\nu(\supp(\nu)) = \nu(\supp(\lambda)) = \lambda(\supp(\lambda)) - \alpha\mu(K_2) + \beta(\sigma - \lambda)(K_1) = 1.\]
The energy of $\nu$ is finite since $\lambda$ and $(\sigma - \lambda)|_{K_1}$ have finite energy.
Then
\[ \int {\mathcal{W}}^{ {\lambda}}    d( {\nu}- {\lambda}) =  \int {\mathcal{W}}^{ {\lambda}} \, d\eta < \gamma \beta(\sigma - \lambda)(K_1) - \gamma \alpha \lambda(K_2)= 0,
\]
in contradiction with $(B')$. So, $\mathcal{W}^\lambda(x) \leq \mathcal{F}_{\lambda}, x \in \supp(\lambda).$ By definition, $\mathcal{W}^\lambda(x) \geq \mathcal{F}_{\lambda}, (\sigma-\lambda)$ almost everywhere. Since  $ {\mathcal{W}}^{ {\lambda}}$ is continuous on $\C$, we have $\mathcal{W}^\lambda(x) \geq \mathcal{F}_{\lambda}, x \in \supp(\sigma-\lambda)$. Thus, $\lambda$  solves $(C')$ with $\mathfrak{w} = \mathcal{F}_{\lambda}$.

The uniqueness of $ {\lambda}$ and the fact that $\supp(\sigma - \lambda) \cap \supp(\lambda)   \neq \emptyset$ imply that  ${\mathfrak{w}}$ is uniquely determined.

If the extremal measure $\lambda$ has compact support, then obviously it satisfies \eqref{cond4}. Now suppose that $\supp(\lambda)$ is unbounded.
Using \eqref{desig2}, we get
\[ \frac{\sqrt{1 + y^2}}{|1 - (y/x)|} \geq \frac{|x|}{\sqrt{1 + x^2}}.
\]
Therefore, for all $x \leq -1$
\[\log \frac{\sqrt{1 + y^2}}{|1 - (y/x)|} \geq - (\log 2)/2.
\]
Using Fatou's Lemma \cite[p. 22]{WR}, $(C')$, and \eqref{admis}, we get
\[ \int \log (1+y^2) d\lambda(y) = 2 \int \liminf_{x\to -\infty} \log\frac{\sqrt{1+y^2}}{|1-(y/x)|} d{\lambda}(y) \leq
 \]
\[\liminf_{x\to -\infty } 2 \int \log\frac{\sqrt{1+y^2}}{|1-(y/x)|} d{\lambda}(y) \leq \liminf_{x\to -\infty, x\in \supp(\lambda) } 2 \int \log\frac{\sqrt{1+y^2}}{|1-(y/x)|} d{\lambda}(y) \leq \]
\[ \leq {\mathfrak{w}} + \limsup_{x \to -\infty} (2 \log |x| - \phi(x)) < +\infty.
\]
Thus, in this case \eqref{cond4} is also fulfilled by $\lambda$.
\hfill $\Box$

\medskip

We are ready to return to the standard potential. Define
\[ {\mathcal{M}}^*(\sigma) := \{\mu \in {\mathcal{M}}(\sigma): I(\mu) < +\infty, \int \log (1 + y^2) d\mu(y) < +\infty\}.
\]
Notice that
\[ {\mathcal{M}}^*(\sigma)   \subset \widetilde{\mathcal{M}}(\sigma) \subset {\mathcal{M}}(\sigma).
\]
According to the last assertion of Theorem \ref{lem1scalar},  $\lambda \in {\mathcal{M}}^*(\sigma)$. Therefore, under the present assumptions, \eqref{minprob} admits the same solution when we minimize the functional over ${\mathcal{M}}^*(\sigma)$.

Set
\[ J_\phi = \inf\{J_{\phi}(\mu): \mu \in {\mathcal{M}}^*(\sigma)\}, \qquad J_{\phi}(\mu) := 2 \left(I(\mu) + \int \phi(x) d\mu(x)\right).
\]
We take $J_{\phi}(\mu) = +\infty$ when $I(\mu) = +\infty$. It is easy to verify that
\[ {\mathcal{J}}_{\phi^*}(\mu)  = J_{\phi}(\mu), \qquad \mu \in {\mathcal{M}}^*(\sigma).
\]
Likewise
\[W^{\mu}(x) := 2U^{\lambda}(x) + \phi(x) = 2\mathcal{U}^{\lambda}(x) + \phi(x) - \int \log (1+y^2)\, d\mu(y), \qquad \mu \in {\mathcal{M}}^*(\sigma).\]
Let
\begin{equation}
\label{F}
F_\mu := \max \{C \in \R: 2 U^{\mu}(x) + \phi(x) \geq C\,\, \mbox{holds}\,\,(\sigma -\mu)\,\, \mbox{a.e.} \}, \qquad \mu \in {\mathcal{M}}^*(\sigma).
\end{equation}
Notice that
\[ F_\mu = \mathcal{F}_\mu - \int \log (1+y^2)\, d\mu(y), \qquad \mu \in {\mathcal{M}}^*(\sigma) .\]

The following result follows from  Theorem \ref{lem1scalar}.

\begin{cor}
\label{potescalar} Under the assumptions of Theorem \ref{lem1scalar}, the following statements are equivalent and have the same unique solution:
\begin{itemize}
\item[$(A'')$] There  exists $ {\lambda} \in {\mathcal{M}}^*(\sigma)$ which is extremal.
\item[$(B'')$] There  exists $ {\lambda} \in {\mathcal{M}}^*(\sigma) $ such that for all $ {\nu} \in  {\mathcal{M}}^*(\sigma)$
    \[\int { {W}}^{ {\lambda}}   d( {\nu}- {\lambda})  \geq 0.\]
\item[$(C'')$] There  exist  $ {\lambda}  \in  {\mathcal{M}}^*(\sigma)$ and a constant ${ {w}}  = { {w}} (\sigma,\phi) $  such that
 \[  2 {{ {U}}}^{ \lambda }(x) + \phi(x)
\left\{
\begin{array}{ll}
\leq { {w}} , & x \in \supp( \lambda ), \\
\geq  { {w}} , & x \in \supp(\sigma -  \lambda ).
\end{array}
\right.
\]
\item[$(D'')$] If $\sigma$ also satisfies \eqref{cond5} then the solution $\lambda$ of $(A'')-(C'')$ verifies
\[{{F}}_\lambda = \max \{ {{F}}_\mu : \mu \in {\mathcal{M}}^*(\sigma)\}.\]
In addition, should
\begin{equation}
\label{condlambda}
\lim_{x\to +\infty} \sqrt{x} \int \log(1 - y/x) d\lambda(y) =0,
\end{equation}
then $\lambda$ is the unique measure which verifies $(D'')$. A sufficient condition for  \eqref{condlambda} is
\begin{equation}
\label{condlambda2}
\int (-y)^{\alpha} d\lambda(y) < \infty, \qquad \alpha > 1/2.
\end{equation}
\end{itemize}
The constant $w(\sigma,\phi) = F_{\lambda}$ is uniquely determined.
\end{cor}

{\bf Proof.} The equivalence of the statements $(A''), (B'')$ and $(C'')$ and the uniqueness of the extremal measure for the functional $J_\phi(\cdot)$   is immediate from Theorem \ref{lem1scalar} and the connections established above.
For  $(D'')$ we have assumed that $\sigma$ also verifies \eqref{cond5}. Then all measures in ${\mathcal{M}}^*(\sigma)$ satisfy \eqref{cond4} and \eqref{cond5} (see sentence right after the introduction of \eqref{cond5}).

\medskip

Notice that  $(C'')$ implies that ${ {F}}_\lambda = w(\sigma,\phi)$.
We must show that ${ {F}}_\mu \leq { {F}}_\lambda$ for all $\mu \in {\mathcal{M}}^*(\sigma)$.
Assume that  ${ {F}}_\mu > { {F}}_\lambda$ for some $\mu \in { \mathcal{M}}^*(\sigma)$. Following the proof of \cite[Theorem 2.1.e]{DrSa}, but replacing the use of the standard principle of domination by Lemma \ref{prindesc}, one obtains that there exists $c > 0$ such that
\[ U^\lambda(x) \leq U^{\mu}(x) - c, \qquad x \in \C. \]
Deleting $\log(1/|x|)$ from both sides and letting $x\to + \infty$ one obtains the contradiction $0 \leq -c$. Therefore,
\begin{equation}
\label{max}
 \max\{{ {F}}_\mu: \mu \in {\mathcal{M}}^*(\sigma)\}= { {F}}_\lambda.
\end{equation}

\medskip

If ${ {F}}_\lambda = { {F}}_\mu$, repeating the scheme used in \cite[Theorem 2.1.e]{DrSa} we arrive to
\[ U^\lambda(x) \leq U^{\mu}(x), \qquad x \in \C. \]
In other words
\[ U^{\mu - \lambda}(x) \geq 0, \qquad x \in \C. \]
If $\supp(\lambda)$ and $\supp(\mu)$ were compact sets, considering that $\lim_{x \to \infty} U^{\mu -\lambda}(x) =0$, this inequality immediately implies, using the minimum principle for harmonic functions, that $U^{\mu -\lambda}(x) \equiv 0, x \in \C \setminus (\supp(\lambda) \cup \supp(\mu))$ which in turn implies that $\mu = \lambda$. Observe that \eqref{condlambda2} is verified when $\lambda$ has compact support.

\medskip

Suppose there exists $x_0 \in \C \setminus \R_-$ where $U^{\mu - \lambda}(x_0) =0$. Then, by the minimum principle $U^{\lambda-\mu}(x) \equiv 0, x \in \C \setminus (\supp(\lambda) \cup \supp(\mu))$ since on the whole boundary (including $\infty$) this harmonic functions has limiting values $\geq 0$. In this case, as in the compact one, we conclude that $\mu  = \lambda$.

\medskip

Assume that $U^{\mu - \lambda}(x) > 0, x \in \C \setminus \R_-$. Define
\[ G^{\mu - \lambda}(x) = \int \log \frac{1}{x-y} d(\mu - \lambda)(y)
\]
the associated complex potential. This function is analytic and never equals zero in $\C \setminus \R_-$. Set
\[ \widetilde{G}^{\mu - \lambda}(z) :=  i G^{\mu - \lambda}(-z^2).
\]
$\widetilde{G}^{\mu -\lambda}$ is analytic and different from zero in $\mbox{Im}(z) >0$, where $\mbox{Im}(\cdot)$ denotes the imaginary part of $(\cdot)$. Moreover,
\[\mbox{Im}(\widetilde{G}^{\mu -\lambda}(z)) = \mbox{Re}(G^{\mu - \lambda}(-z^2)) = U^{\mu -\lambda}(-z^2) > 0, \qquad \mbox{Im}(z) > 0.\]
Therefore, $\widetilde{G}^{\mu -\lambda}$ transforms the upper half plane into the upper half plane. From here we have an integral representations for $\widetilde{G}^{\mu -\lambda}(z)$.

\medskip

Indeed, from \cite[Theorem A.2]{KN}, we know that
\begin{equation}
\label{repint1}
 \widetilde{G}^{\mu - \lambda}(z) = \kappa + \beta z +  \int_\R  \left(\frac{1}{t-z} - \frac{t}{1+t^2}
  \right) d \rho(t),
\end{equation}
where $\kappa \in \R, \beta \geq 0,$ and $\rho$ is a positive Borel measure on $\R$ such that $\int (1+t^2)^{-1} d\rho(t) < \infty$.  Similarly,  from \cite[Theorem A.3]{KN}, it follows that
\begin{equation}
\label{repint2}
 \log \left(\widetilde{G}^{\mu - \lambda}(z)\right) = \gamma + \int_\R  \left(\frac{1}{t-z} - \frac{t}{1+t^2}
  \right) f(t) dt,
\end{equation}
where $\gamma \in \R$ and $f$ is an integrable function on $\R$ such that $0 \leq f(t) \leq 1$ almost everywhere. Let us simplify these representations a bit.

\medskip

If $z = iu, u > 0$, using the definition of $\widetilde{G}^{\mu -\lambda}$, it follows that
\begin{equation}
\label{symmetric}
\widetilde{G}^{\mu -\lambda}(iu) = i \int \log \frac{1}{|u^2 -y|} d(\mu -\lambda)(y) - \int \arg \frac{1}{u^2 -y} d(\mu -\lambda)(y) =
\end{equation}
\[i \int \log \frac{1}{|u^2 -y|} d(\mu -\lambda)(y) = i U^{\mu - \lambda}(u^2),\]
is purely imaginary. By the symmetry principle, $\widetilde{G}^{\mu -\lambda}$ is symmetric with respect to the imaginary axis. That is for $\mbox{Im} z >0$,
\begin{equation}
\label{symmetric2}
\mbox{Im}(\widetilde{G}^{\mu -\lambda}(z))= \mbox{Im}(\widetilde{G}^{\mu -\lambda}(-\overline{z})), \qquad \mbox{Re}(\widetilde{G}^{\mu -\lambda}(z))= -\mbox{Re}(\widetilde{G}^{\mu -\lambda}(-\overline{z})).
\end{equation}
In particular,
\begin{equation}
 \label{rel1}
\arg\left(\widetilde{G}^{\mu - \lambda}(z) \right) = \pi - \arg\left(\widetilde{G}^{\mu - \lambda}(-\overline{z}) \right), \qquad \mbox{Im} z > 0.
 \end{equation}
Actually, $\widetilde{G}^{\mu - \lambda}$ can be extended continuously to $\R$ from the upper half plane; therefore, the last relation implies that
\begin{equation}
\label{rel2}
\arg\left(\widetilde{G}^{\mu - \lambda}(t) \right)_+ = \pi - \arg\left(\widetilde{G}^{\mu - \lambda}(-t) \right)_+, \qquad t \in \R.
\end{equation}

\medskip

Due to the Stieltjes inversion formula, the first relation in \eqref{symmetric2} implies that the measure $\rho$ is symmetric with respect to the origin ($d\rho(t) = d\rho(-t)$). Therefore, \eqref{repint1} can be transformed as follows
\begin{equation}
\label{repint4}
 \widetilde{G}^{\mu - \lambda}(z) = \kappa + \beta z +  \int_{-\infty}^0 \left(\frac{1}{t-z} - \frac{t}{1+t^2}
  \right) d \rho(t) + \int_0^{\infty} \left(\frac{1}{t-z} - \frac{t}{1+t^2}
  \right) d \rho(t) =
\end{equation}
\[ \kappa + \beta z + \int_{-\infty}^0 \frac{2z}{t^2-z^2} d\rho(t).\]
Evaluating \eqref{repint4} at  $iu$ we obtain a purely imaginary number (see \eqref{symmetric}) so comparing both sides we see that $\kappa = 0$. Now, dividing by $u$ and letting $u$ tend to $\infty$, we get that $\beta=0$. Consequently,
\[ \widetilde{G}^{\mu - \lambda}(z) =  i G^{\mu - \lambda}(-z^2) = \int_{-\infty}^0 \frac{2z}{t^2-z^2} d\rho(t).
\]
Changing variables $-z^2 = x, -t^2 = y$, we obtain
\begin{equation}
\label{repint5}
G^{\mu - \lambda}(x) = \int_{-\infty}^0 \frac{2\sqrt{x}}{x-y} d\widetilde{\rho}(y), \qquad x \in \C \setminus \R_-,
\end{equation}
with $\sqrt{1} = 1$ and $d\widetilde{\rho}(y) = d\rho(\sqrt{-y})$. Notice that $\int (1+|y|)^{-1} d\widetilde{\rho}(y) < \infty$.

\medskip

Take $x > 0$ and $N > 0$. From \eqref{repint5}, we have
\[\sqrt{x} G^{\mu - \lambda}(x) \geq \int_{-N}^0 \frac{2 {x}}{x-y} d\widetilde\rho(y), \qquad x \in \C \setminus \R_-.\]
Assume that $\sqrt{x} G^{\mu - \lambda}(x) \leq M$ for all sufficiently large $x$. Taking $\limsup$ as $x \to \infty$, it follows that
$\widetilde{\rho}[-N,0] \leq M/2$. Since this would take place for all $N >0$ we would conclude that $\widetilde{\rho}$ is finite with total mass $\leq M/2$. Now, if $\lim_{x \to \infty} \sqrt{x} G^{\mu - \lambda}(x) = 0$, we would have that $\widetilde{\rho}$ is the null measure and \eqref{repint5} would render that $G^{\mu -\lambda}(x) \equiv 0, x\in \C \setminus \R_-$, implying  $\mu = \lambda$ as we wish.

\medskip

Notice that for $x > 0$, we have
\[ 0 \leq \sqrt{x} G^{\mu - \lambda}(x) = \sqrt{x} U^{\mu - \lambda}(x) = \]
\[\sqrt{x}\left(\int \log(1 - y/x) d(\lambda - \mu)(y)\right) \leq \sqrt{x} \int \log(1 - y/x) d\lambda(y).
\]
Therefore, $\lim_{x\to +\infty} \sqrt{x} G^{\mu - \lambda}(x)= 0$ under  \eqref{condlambda}.

\medskip

On the other hand, if \eqref{condlambda2} takes place we can assume that $1/2 < \alpha \leq 1$
\[ 0 \leq  \sqrt{x} \int \log(1 - y/x) d\lambda(y) =  \frac{\sqrt{x}}{\alpha} \int \log(1 - y/x)^{\alpha} d\lambda(y) \leq  \]
\[\frac{\sqrt{x}}{\alpha} \int \log(1 + (-y/x)^{\alpha}) d\lambda(y) \leq \frac{\sqrt{x}}{\alpha} \int  (-y/x)^{\alpha} d\lambda(y). \]
Consequently, \eqref{condlambda2} is sufficient to have \eqref{condlambda}.The proof of Corollary \ref{potescalar} is complete. \hfill $\Box$

\medskip

Alternatively, we could have concluded the proof of Corollary \ref{potescalar} with the following arguments which lead to a different integral representation. Let $\tau$ denote the distribution function of the measure $f(t) dt$. By the Stieltjes inversion formula
\[ \tau(t_2) - \tau(t_1) = \lim_{\varepsilon \to 0} \frac{1}{\pi} \int_{t_1}^{t_2} \arg \left(\widetilde{G}^{\mu - \lambda}(t + i\varepsilon)\right) dt, \qquad t_1 < t_2.
\]
Using \eqref{rel1}-\eqref{rel2} it follows that for $\infty < t_1 < t_2 \leq 0$
\[\tau(t_2) - \tau(t_1) =   \lim_{\varepsilon \to 0} \frac{1}{\pi} \int_{t_1}^{t_2} \arg \left(\widetilde{G}^{\mu - \lambda}(t + i\varepsilon)\right) dt =
\]
\[ \frac{1}{\pi} \int_{t_1}^{t_2} \arg \left(\widetilde{G}^{\mu - \lambda}(t)\right)_+ dt =  t_2 - t_1 - \frac{1}{\pi} \int_{t_1}^{t_2} \arg \left(\widetilde{G}^{\mu - \lambda}(-t)\right)_+ dt  = \]
\[t_2 - t_1 - \frac{1}{\pi} \int_{-t_2}^{-t_1} \arg \left(\widetilde{G}^{\mu - \lambda}(t)\right)_+ dt = t_2 - t_1 - (\tau(-t_1) - \tau(-t_2)), \]
Consequently, almost everywhere on $\R$, we have
\begin{equation}
\label{ft}
f(t) = \frac{d\tau(t)}{dt} = 1 - f(-t).
\end{equation}

\medskip

From \eqref{repint2} and \eqref{ft}, we obtain
\[ \log \left(\widetilde{G}^{\mu - \lambda}(z)\right) = \gamma + \int_{-\infty}^0  \left(\frac{1}{t-z} - \frac{t}{1+t^2}
  \right) f(t) dt +
  \]
  \[\int_{0}^{\infty}  \left(\frac{1}{t-z} - \frac{t}{1+t^2}
  \right) (1-f(-t)) dt =\]
  \[\gamma + \int_{-\infty}^0  \left(\frac{1}{t-z} + \frac{1}{t+z}- \frac{2t}{1+t^2}
  \right) f(t) dt  + \int_0^{\infty} \frac{1 + tz}{t-z} \frac{dt}{1+t^2} = .
\]
\[\gamma + 2(1 + z^2) \int_{-\infty}^0   \frac{t f(t)  }{t^2-z^2} \frac{dt}{1 + t^2}
     + \int_0^{\infty} \frac{1 + tz}{t-z} \frac{dt}{1+t^2} .
\]

\medskip

Integrating with respect to $t$ the function $(1 + tz)\log (t)/(t-z)$,   over the closed contour consisting of the circles $\{t: |t| = R\}, \{t:|t| = \varepsilon\}$, and the segment $[\varepsilon,R]$ oriented positively, where the branch of the logarithm in $\C \setminus \R_+$ is taken so that $\log (-1) = i\pi$, and using the residue theorem one obtains
\[ \int_0^{\infty} \frac{1 + tz}{t-z} \frac{dt}{1+t^2} = i\pi - \log (z), \qquad z \in \C \setminus \R_+ .
\]
Therefore,
\begin{equation}
\label{repint3}
\log \left(\widetilde{G}^{\mu - \lambda}(z)\right) = \gamma + i\pi - \log(z) +  2(1 + z^2) \int_{-\infty}^0   \frac{t f(t)  }{t^2-z^2} \frac{dt}{1 + t^2},\quad  \mbox{Im}(z) > 0,
\end{equation}
or what is the same,
\[ \log\left({G}^{\mu - \lambda}(-z^2)\right) = \gamma + \frac{i\pi}{2} - \log(z) +  2(1 + z^2) \int_{-\infty}^0   \frac{t f(t)  }{t^2-z^2} \frac{dt}{1 + t^2},\quad  \mbox{Im}(z) > 0,\]
Making the change of variables $-z^2 = x$ and $-t^2 = y$, this relation becomes
\[ \log\left({G}^{\mu - \lambda}(x)\right) = \gamma   -  \log(\sqrt{x}) +  \left(1 - \frac{1}{x}\right) \int_{-\infty}^0   \frac{x}{x-y} \frac{f(-\sqrt{|y|})dy}{1 + |y|},\qquad  x \in \C \setminus \R_-,\]
where $\sqrt{1} = 1, \log(1) = 0$. Evaluating at $x=1$, it follows that $\gamma = \log\left({G}^{\mu - \lambda}(1)\right)$. Therefore
\[ \log\left(\frac{\sqrt{x}{G}^{\mu - \lambda}(x)}{{G}^{\mu - \lambda}(1)}\right) =   \left(1 - \frac{1}{x}\right) \int_{-\infty}^0   \frac{x}{x-y} \frac{f(-\sqrt{|y|})dy}{1 + |y|},\qquad  x \in \C \setminus \R_-.\]
Notice that for $x > 1$ the right hand is positive. So ${\sqrt{x}{G}^{\mu - \lambda}(x)} > {G}^{\mu - \lambda}(1)$ for all $x > 1$, which is not possible under \eqref{condlambda} unless ${G}^{\mu - \lambda}(1) = 0$ which implies, as we know, that $\mu = \lambda$.

\medskip

Let us see some other properties of the extremal measure.

\begin{cor}
\label{miscelanea} Suppose that the assumptions of Theorem \ref{lem1scalar} are verified.
\begin{itemize}
\item[(a)] If $\liminf_{x \to -\infty} \phi^*(x) = + \infty$, then $\supp(\lambda)$ is compact.
\item[(b)] If $\supp(\lambda)$ is unbounded and $\lambda$ verifies \eqref{cond5},
then

$\liminf_{x\to -\infty }   \phi^*(x) \leq w(\sigma,\phi) $.
\item[(c)] Should $\int \log(1+y^2) d\sigma(y) = +\infty$, then $\supp(\sigma -\lambda)$ is unbounded.
\item[(d)] If $\supp (\sigma - \lambda)$  is unbounded and $\lambda$ verifies \eqref{cond5}, then

$\liminf_{x\to -\infty, x \in \supp(\sigma -\lambda)} \phi^*(x) \geq  {w}(\sigma,\phi)$.
\item[(e)] If $\supp (\sigma - \lambda)$ and $\supp(\lambda)$ are unbounded, $\lambda$ verifies \eqref{cond5}, and $\lim_{x\to -\infty} \phi^*(x)$ exists, the limit is $w(\sigma,\phi)$.
\item[(f)] Assume that  $\phi(x)$ is decreasing on $\R_-$, then $0 \in \supp(\lambda)$.
\item[(g)] Should $x\phi'(x)$ be decreasing on $\R_-$, then $\supp(\lambda)$ is connected.
\item[(h)] Let $\phi(x) = -U^{\tau}(x)$, where $\tau \in \mathcal{M}_2^+(\R_+)$ has compact support and $U^{\tau}(x)$ is continuous at $x=0$, then $\supp(\lambda) = \R_-$. If $\supp(\sigma -\lambda)$ is unbounded and $\lambda$ satisfies \eqref{cond5} then $w(\sigma,\phi) = 0$.
\end{itemize}
\end{cor}

{\bf Proof.} According to  $(C')$
\[\phi^*(x) \leq  2 \int \log \frac{\sqrt{1+x^2}\sqrt{1+y^2}}{|x-y|} d\lambda(y) + \phi^*(x) \leq \mathfrak{w}(\sigma,\phi), \qquad x \in \supp(\lambda).
\]
If $\supp (\lambda)$ is unbounded, it follows that
\[ \limsup_{x\to -\infty, x \in \supp(\lambda)} \phi^*(x) \leq \mathfrak{w}.\]
Therefore, if $\liminf_{x\to -\infty} \phi^*(x) = +\infty$ we get a contradiction. Thus {\rm (a)} takes place.

\medskip

According to $(C'')$ we have
\[   W^{\lambda}(x) =  2 \int \log \frac{1}{|1-(y/x)|} d\lambda(y) + \phi^*(x) + \log \frac{1+x^2}{x^2} \leq  {w}, \qquad x \in \supp(\lambda).\]
If $\supp(\lambda)$ is unbounded and $\lambda$ verifies \eqref{cond5}, due to \eqref{unif2}  it follows that
\[\liminf_{x\to -\infty } \phi^*(x) \leq \liminf_{x\to -\infty, x \in \supp(\lambda)} \phi^*(x) \leq w.\]
Therefore, {\rm (b)} is valid.

\medskip

Suppose that $\supp(\sigma - \lambda)$ is a compact set $K$. We have $\lambda|_{\R_- \setminus K} = \sigma|_{\R_- \setminus K}$. However, $\int \log(1+y^2) d\lambda(y) < + \infty$. Consequently, $\int \log(1+y^2) d\sigma(y) < + \infty$. We conclude that {\rm (c)} holds.

\medskip

From $(C'')$ we know that
\[  W^{\lambda}(x) =  2 \int \log \frac{1}{|1-(y/x)|} d\lambda(y) + \phi^*(x) + \log \frac{1+x^2}{x^2} \geq  {w}, \qquad x \in \supp(\sigma - \lambda).\]
Thus, if $\supp(\sigma \! - \! \lambda)$ is unbounded and $\lambda$ verifies \eqref{cond5}, for $x\to -\infty, x \in \supp(\sigma - \lambda),$ from  \eqref{unif2} we obtain {\rm (d)}. Now, {\rm (e)} is a direct consequence of {\rm (b)} and {\rm (d)}.

\medskip

For $x \in \R \setminus \supp(\lambda)$, we have
\[ \left(U^{\lambda}(x)\right)' = -\int \frac{d\lambda(y)}{x-y}, \qquad \left(x\left(U^{\lambda}(x)\right)'\right)'= \int \frac{yd\lambda(y)}{(x-y)^2}.
\]
If $\phi$ decreases on $\R_-$ and $0 \not\in \supp(\lambda)$ the first of these formulas implies that $W^{\lambda}(x)$ decreases immediately to the right of $\supp(\lambda)$ but this contradicts $(C'')$; therefore, {\rm (f)} follows. Should $x\phi'(x)$ be decreasing, the second formula implies that $x \left(W^{\lambda}(x)\right)'$ is decreasing on any connected component of $\R_- \setminus \supp(\lambda)$. From here it follows that  $\left(W^{\lambda}(x)\right)'$ cannot change sign from plus to minus on any such connected component. Suppose that $\supp(\lambda)$ is not connected, then there exist $x_1,x_2 \in \supp(\lambda_2)$, $x_2 < 0,$ such that $(x_1,x_2) \cap  \supp(\lambda_2)=\emptyset$. According to $(C'')$, $\left(W^{\lambda}(x)\right)'$ changes sign from  plus to minus on $(x_1,x_2)$; thus $\supp(\lambda)$ must be connected and we obtain {\rm (g)}.

\medskip

Finally, it is easy to check that $\phi = - \mathcal{U}^{\tau}$, as indicated in part {\rm (g)}, is decreasing on $\R_-$ and $x\phi'(x)$ is decreasing in $\R_-$; therefore, according to $(f)$ and $(g)$, $\supp(\lambda)$ is a closed interval in $\R_-$ which touches $x=0$. Suppose that $\supp(\lambda)$ is bounded. Then, $W^{\lambda}(x)$ is subharmonic in $\overline{\C} \setminus \supp(\lambda)$, continuous on $\supp(\lambda)$, and using the second part of $(C'')$
\[W^{\lambda}(\infty) = \lim_{x \to -\infty} W^{\lambda}(x) = 0 \geq w.\]
However, $W^{\lambda}(x)  \leq w, x \in \supp(\lambda),$ as the first part of $(C'')$ states. Using the maximum principle for subharmonic functions it follows that $2 U^{\lambda}(x) \equiv U^{\tau}(x), x \in \C \setminus \supp(\lambda)$ which is impossible. Therefore, $\supp(\lambda) = \R_-$. Now, if $\lambda$ satisfies \eqref{cond5} and $\supp(\sigma-\lambda)$ is unbounded from {\rm (e)} we get $w(\sigma,\phi) =0$. \hfill $\Box$

\medskip

\begin{rem}
In this corollary we have assumed on several occasions that $\lambda$ satisfies \eqref{cond5}. One way to ensure this  is requiring in the initial data that $\sigma$ fulfills this condition. However, it is possible that $\lambda$  satisfies \eqref{cond5} but not necessarily $\sigma$ (for example, when $\supp(\lambda)$ is compact, see part {\rm (a)} of the corollary). In connection with {\rm (h)} notice that the unboundedness of $\supp(\sigma-\lambda)$ is ensured when $\int \log(1+ y^2) d\sigma(y) = +\infty$ (see {\rm (c)}).
\end{rem}

\medskip

\begin{rem} \label{unconstrained}
We wish to call attention to the case when $\sigma \equiv +\infty$ which corresponds to an equilibrium problem with no constraint. This case is considered in \cite{GR2}. In this situation, one cannot rely on $\sigma$ to guarantee that $\lambda$ verifies \eqref{cond5} or deduce the continuity of   $\mathcal{U}^\lambda$. Nevertheless, if $\liminf \phi^* = +\infty$, one can assert that $\lambda$ has compact support which in turn trivially implies \eqref{cond5} on $\lambda$ and the continuity of $\mathcal{U}^{\lambda}$ follows from $(C')$ since $2\mathcal{U}^{\lambda}$ is equal on $\supp(\lambda)$ to the continuous function $\mathfrak{w} - \phi$.

\end{rem}

\section{Proof of Theorem \ref{teo:equil}}

In this section, we use again the notion of modified potential \eqref{modpot} and modified energy \eqref{modener} introduced in Section 4.2.

\medskip

Let $\varphi$ be a continuous function on ${\mathbb{R}}_+$ which verifies
\begin{equation} \label{C2} \liminf_{x\to +\infty} \left( 2\varphi(x) - 3 \log(1 + x^2)\right) > -\infty.
\end{equation}
This assumption is much weaker than \eqref{cond1}.
Set $\varphi^*(x) := \varphi(x) - \frac{3}{2}\log(1+ x^2)$, and define
\[ \mathcal{A} =
\begin{pmatrix}
2 & -1 \\
-1 & 2
\end{pmatrix},
\qquad
 f =  \begin{pmatrix}
\varphi^* \\
0
\end{pmatrix}.
\]
For $ \vec{\mu}=(\mu_1,\mu_2)^t \in \mathfrak{M}(\sigma)$ (see the definition in \eqref{star}), we introduce the vector function
\[ {\mathcal{W}}^{\vec{\mu}}(x) = ({\mathcal{W}}_1^{\vec{\mu}}(x),{\mathcal{W}}_2^{\vec{\mu}}(x))^t := \int \log \frac{\sqrt{1 + x^2}\sqrt{1 + y^2}}{|x-y|} d \mathcal{A} \vec{\mu}(y) + f(x)
\]
and the functional
\begin{equation} \label{func}
{\mathcal{J}}_{\varphi^*}(\vec{\mu}) := \int ({\mathcal{W}}^{\vec{\mu}} + f) \cdot d \vec{\mu} =  \int ({\mathcal{W}}_1^{\vec{\mu}} + \varphi^*) d \mu_1 + \int {\mathcal{W}}_2^{\vec{\mu}} d \mu_2
\end{equation}
(when either ${\mathcal{I}}(\mu_1) = +\infty$ or ${\mathcal{I}}(\mu_2) = +\infty$, we take ${\mathcal{J}}_{\varphi^*}(\vec{\mu})= +\infty$). That is,
\[ {\mathcal{J}}_{\varphi^*}(\vec{\mu}) = 2({\mathcal{I}}(\mu_1) - {\mathcal{I}}(\mu_1,\mu_2) + {\mathcal{I}}(\mu_2)) + \int (2\varphi - 3 \log(1+x^2)) d\mu_1.
\]
Condition \eqref{C2} and the fact that $\mathcal{A}$ is positive definite guarantee that the corresponding vector equilibrium problem is weakly admissible as defined in  \cite[Assumption 2.1]{HK}. In particular (see \cite[Corollary 2.7]{HK} and the sentence that follows it), this guarantees that
\[ {\mathcal{J}}_{\varphi^*} = \inf \{{\mathcal{J}}_{\varphi^*}(\vec{\mu}): \vec{\mu} \in \mathfrak{M}(\sigma)\} > -\infty.
\]

\medskip

Set
\[\,\, \quad \widetilde{\mathfrak{\mathfrak{M}}}(\sigma) = \{\vec{\mu} \in {\mathfrak{M}}(\sigma): \mathcal{I}(\mu_1) < \infty, \mathcal{I}(\mu_2) < \infty\}, \] \[\mathfrak M^*(\sigma) = \{\vec{\mu} \in {\mathfrak{M}}(\sigma): \mu_1,\mu_2\,\, \mbox{verify}\,\, \eqref{C1}\}.\]
A vector measure $\vec{\lambda} \in \widetilde{\mathcal{\mathfrak{M}}}(\sigma)$  is said to be extremal if
\[ -\infty < {\mathcal{J}}_{\varphi^*}(\vec{\lambda}) = {\mathcal{J}}_{\varphi^*}  < +\infty.
\]
In case that $\vec{\mu} \in \mathfrak M^*(\sigma)$, it is easy to check that
\begin{equation} \label{funcionales} {\mathcal{J}}_{\varphi^*}(\vec{\mu}) = 2\left(I(\mu_1) - I(\mu_1,\mu_2) + I(\mu_2) +\int \varphi\, d\mu_1\right) := J_{\varphi}(\vec{\mu}).
\end{equation}

\medskip

The next theorem complements, in the present setting, results from \cite{HK}.
\begin{teo} \label{lem1} Let $\varphi$ satisfy \eqref{C2} and let $\sigma$, $\supp (\sigma) = {\R}_-$, $|\sigma| > 1$, be a positive Borel measure such that $U^{\sigma|_K}$ is continuous on $\C$ for every compact subset $K \subset {\R}_-$. The following statements are equivalent and have the same unique solution:
\begin{itemize}
\item[$(A''')$] There  exists $\vec{\lambda} \in \widetilde{\mathfrak{M}}(\sigma)$ which is extremal.
\item[$(B''')$] There  exists $\vec{\lambda} \in \widetilde{\mathfrak{M}}(\sigma)$ such that for all $\vec{\nu} \in \widetilde{\mathfrak{M}}(\sigma)$
    \[\int {\mathcal{W}}^{\vec{\lambda}} \cdot d(\vec{\nu}-\vec{\lambda}) := \int {\mathcal{W}}_1^{\vec{\lambda}} d( {\nu}_1-{\lambda}_1) +   \int {\mathcal{W}}_2^{\vec{\lambda}} d( {\nu}_2-{\lambda}_2)\geq 0 .\]
\item[$(C''')$] There  exist  $\vec{\lambda} = ( \lambda_1,  \lambda_2) \in \widetilde{\mathfrak{M}}(\sigma)$ and constants ${\mathfrak{w}}_1 = {\mathfrak{w}}_1(\sigma,\varphi),{\mathfrak{w}}_2 = {\mathfrak{w}}_2(\sigma,\varphi)$  such that
\begin{itemize}
\item[$(i)$]
\[ {\mathcal{W}}_1^{\vec{\lambda}}(x) = 2 {{\mathcal{U}}}^{ \lambda_1}(x) - {{\mathcal{U}}}^{ \lambda_2}(x) + \varphi(x)
\left\{
\begin{array}{ll}
= {\mathfrak{w}}_1, & x \in \supp( \lambda_1), \\
\geq {\mathfrak{w}}_1, & x \in {\mathbb{R}}_+,
\end{array}
\right.
\]
\item[$(ii)$]
\[ {\mathcal{W}}_2^{\vec{\lambda}}(x) = 2 {{\mathcal{U}}}^{ \lambda_2}(x) - {{\mathcal{U}}}^{ \lambda_1}(x)
\left\{
\begin{array}{ll}
\leq {\mathfrak{w}}_2, & x \in \supp( \lambda_2), \\
\geq  {\mathfrak{w}}_2, & x \in \supp(\sigma -  \lambda_2).
\end{array}
\right.
\]
\end{itemize}
\end{itemize}
The constants ${\mathfrak{w}}_1,{\mathfrak{w}}_2$ are uniquely determined.  ${\mathcal{U}}^{\lambda_1}$ and ${\mathcal{U}}^{\lambda_2}$ are continuous on $\mathbb{C}$.
\end{teo}

{\bf Proof.} The proof is similar to that of Theorem \ref{lem1scalar} so we will be brief. As shown in \cite[Theorem 2.6]{HK}, the functional ${\mathcal{J}}_{\varphi^*}$ is lower semicontinuous and strictly convex on ${\mathfrak{M}}(\sigma)$, from which the existence of a unique solution to $(A''')$ is guaranteed, see \cite[Corollary 2.7]{HK}. By the way in which the functional is defined, the extremal measure must belong to $\widetilde{\mathfrak{M}}(\sigma)$.

\medskip

The equivalence of $(A''')$ and $(B''')$ comes from the identity
\[ {\mathcal{J}}_{\varphi^*}(\vec{\nu}_{\varepsilon}) - {\mathcal{J}}_{\varphi^*}(\vec{\lambda}) = \varepsilon^2 {\mathcal{J}}_0(\vec{\nu} - \vec{\lambda}) + 2\varepsilon \int {{\mathcal{W}}}^{\vec{\lambda}} \cdot d(\vec{\nu} - \vec{\lambda}),
\]
valid for any $\vec{\lambda}, \vec{\nu} \in \widetilde{\mathfrak{M}}(\sigma)$ and $0 \leq \varepsilon \leq 1$, where $\vec{\nu}_{\varepsilon} = \varepsilon \vec{\nu} + (1-\varepsilon)\vec{\lambda}$ and ${\mathcal{J}}_0(\vec{\nu} - \vec{\lambda})$ is the functional applied to $\vec{\nu} - \vec{\lambda}$ with $\varphi^* \equiv 0$. To prove $(B''')$ implies $(A''')$ one also uses that ${\mathcal{J}}_0(\vec{\nu} - \vec{\lambda}) \geq 0$ with equality   only if $\vec{\nu} = \vec{\lambda}$ (see \cite[Proposition 3.5]{HK} and \cite[Theorem 2.5]{CKL}).

\medskip

If $\vec{\lambda} = ( \lambda_1, \lambda_2)^t$  verifies $(C''')$ and $\vec{\nu} = (\nu_1,\nu_2)^t \in \widetilde{\mathfrak{M}}(\sigma) $.
From $(C'''-i)$, we have
\[ \int {{\mathcal{W}}}_1^{\vec{\lambda}} \, d(\nu_1 -  \lambda_1) = \int {{\mathcal{W}}}_1^{\vec{\lambda}} \, d\nu_1 - \int {{\mathcal{W}}}_1^{\vec{\lambda}} \, d \lambda_1  \geq {\mathfrak{w}}_1 - {\mathfrak{w}}_1 = 0.
\]
On the other hand, $ |\lambda_2|= |\nu_2| = 1$; therefore,
\[  \int {{\mathcal{W}}}_2^{\vec{\lambda}} \, d(\nu_2 -  \lambda_2) = \int ({{\mathcal{W}}}_2^{\vec{\lambda}} - {\mathfrak{w}}_2) \, d(\nu_2 -  \lambda_2).
\]
To show that this integral is also $\geq 0$ one uses the same arguments as in proving $(C')$ implies $(B')$ defining now
\[ E_+ = \{ t \in \R_- : {{\mathcal{W}}}_2^{\vec{\lambda}}(t) - {\mathfrak{w}}_2 > 0\}, \qquad E_- = \{ t \in R_-: {{\mathcal{W}}}_2^{\vec{\lambda}}(t) - {\mathfrak{w}}_2 < 0\}.
\]
Putting these relations together, we obtain
\[ \int {{\mathcal{W}}}^{\vec{\lambda}} \cdot  d(\vec{\nu} - \vec{\lambda}) \geq 0, \qquad \nu \in \widetilde{\mathfrak{M}}(\sigma).
\]
So, $(C''')$ implies $(B''')$.

\medskip

Assume that  $\vec{\lambda} = ( \lambda_1, \lambda_2)^t$  solves $(B''')$. Set
\[ {\mathfrak{w}}_1 := \frac{1}{2} \int {\mathcal{W}}_1^{\vec{\lambda}} \, d\lambda_1.
\]
Let us prove that
\begin{equation} \label{quasi1} {\mathcal{W}}_1^{\vec{\lambda}}(x) \geq {\mathfrak{w}}_1 \quad \mbox{quasi-everywhere on}\,\, \R_+\,,
\end{equation}
where ``quasi-everywhere'' means except on a set of capacity zero. If this was not so, there would exist  a compact subset   $K_1 \subset {\R}_+$, $\mbox{cap}(K_1) > 0,$ such that ${\mathcal{W}}_1^{\vec{\lambda}}(x) < {\mathfrak{w}}_1$, $x \in K_1$. Taking $\nu_1 \in {\mathcal{M}}_2^+({\R}_+)$, $\supp(\nu_1) \subset K_1,$ and $\nu_2 = \lambda_2$, we obtain
\[ \int {\mathcal{W}}^{\vec{\lambda}} \cdot  d(\vec{\nu}-\vec{\lambda}) =  \int {\mathcal{W}}_1^{\vec{\lambda}} \, d(\nu_1 -\lambda_1) < 2{\mathfrak{w}}_1 -2{\mathfrak{w}}_1 = 0,
\]
which contradicts $(B''')$. Now, we prove that
\[ {\mathcal{W}}_1^{\vec{\lambda}}(x) \leq {\mathfrak{w}}_1, \quad x \in \supp(\lambda_1).
\]
To the contrary, assume that there exists $x_0 \in \supp(\lambda_1)$ such that ${\mathcal{W}}_1^{\vec{\lambda}}(x_0) > {\mathfrak{w}}_1$. By the lower semi-continuity of  ${\mathcal{W}}_1^{\vec{\lambda}}$ on $\R_+$ (${\mathcal{U}}^{\lambda_2}$ is continuous by Lemma \ref{cont} and $\varphi$ by assumption) it follows that there exists $\delta >0$ such that ${\mathcal{W}}_1^{\vec{\lambda}}(x) > {\mathfrak{w}}_1$, $|x-x_0| \leq \delta$. Take $K_2 = \supp(\lambda_1) \cap \{x:|x-x_0|\leq \delta\}$. Then $\lambda_1(K_2) > 0$ and
\[ 2{\mathfrak{w}}_1 =   \int_{\supp(\lambda_1) \setminus K_2} {\mathcal{W}}_1^{\vec{\lambda}} \, d\lambda_1 +  \int_{K_2} {\mathcal{W}}_1^{\vec{\lambda}}\,  d\lambda_1 >  {\mathfrak{w}}_1 (\lambda_1(\supp(\lambda_1) \setminus K_2) + \lambda_1(K_2)) = 2{\mathfrak{w}}_1,
\]
which is also a contradiction. From \eqref{quasi1}, reasoning as in \cite[Theorem 5.4.1]{NS}, it follows that ${\mathcal{W}}_1^{\vec{\lambda}}\geq {\mathfrak{w}}_1$ on all $\R_+$. Hence,  $(C'''-i)$ is obtained. We have also obtained that ${\mathcal{U}}^{\lambda_1}$ is continuous on all $\C$ because  on $\supp(\lambda_1)$ it is equal to the continuous function $\frac{1}{2}\left({\mathfrak{w}}_2 - \varphi +  {\mathcal{U}}^{\lambda_2}\right)$.

\medskip

For the proof of  $(C'''-ii)$   take
\[  {\mathfrak{w}}_2 := \sup\{{\mathfrak{w}} \in {\R}: {\mathcal{W}}_2^{\vec{\lambda}} \geq {\mathfrak{w}} \quad (\sigma - \lambda_2)\,\,\mbox{a.e.}\}.
\]
If there exists $x_0 \in \supp(\lambda_2)$ such that ${\mathcal{W}}_2^{\vec{\lambda}}(x_0) > {\mathfrak{w}}_2$ proceeding as in the scalar case one can construct a signed measure $\eta$ of total mass $1$ supported on a compact subset of $\R_-$ such that $\vec{\nu} := (\lambda_1, \lambda_2 + \eta)^t \in \widetilde{\mathfrak{M}}(\sigma)$ and
\[ \int {\mathcal{W}}^{\vec{\lambda}} \cdot  d(\vec{\nu}-\vec{\lambda}) =  \int {\mathcal{W}}_2^{\vec{\lambda}} \, d\eta <   0,
\]
in contradiction with $(B''')$. From the continuity of $ {\mathcal{W}}_2^{\vec{\lambda}}$ on $\C$, the inequality in the second part of $(C''' -ii)$ holds for all $x \in \supp(\sigma - \lambda_2)$. Therefore, $(C''')$ has been proved.

\medskip

From the uniqueness of $\vec{\lambda}$ and the fact that $\supp(\sigma - \lambda_2) \cap \supp(\lambda_2)   \neq \emptyset$ it readily follows that  ${\mathfrak{w}}_1, {\mathfrak{w}}_2$ are uniquely determined.   \hfill $\Box$

\medskip

\begin{cor} \label{varios} With the assumptions of Theorem \ref{lem1}, let $\vec{\lambda}$ be extremal. Then, $\supp(\lambda_2)$ is connected and $0 \in \supp(\lambda_2)$. If $x\varphi'(x)$ is an increasing function on ${\R}_+$ then $\supp(\lambda_1)$ is connected. If $\varphi$ is increasing on $\R_+$ then $0 \in \supp(\lambda_1)$. If
\begin{equation}
\label{conditionInf}
\lim_{x \to +\infty} \left( \varphi(x) - 4\log x\right) = +\infty,
\end{equation}
then $\supp(\lambda_1)$ is a compact set, $\supp(\lambda_2) = \R_-$, and $\lambda_1, \lambda_2$ verify \eqref{cond4}.
\end{cor}

{\bf Proof.} Notice that for any finite measure $\mu$ on the real line $\left(\mathcal{U}^{\lambda}(x)\right)' = \left({U}^{\lambda}(x)\right)'$ and thus $\left(x\left(\mathcal{U}^{\lambda}(x)\right)'\right)' = \left(x\left({U}^{\lambda}(x)\right)'\right)'$ for all $x \in \R \setminus \supp(\lambda)$. Arguing as in Corollary \ref{miscelanea} (f)-(g) one proves that $\supp(\lambda_2)$ is connected and $0 \in \supp(\lambda_2)$. Similarly, one proves that $\supp(\lambda_1)$ is connected and $0 \in \supp(\lambda_1)$ when $x\varphi'$ and $\varphi$ are increasing, respectively.

\medskip

The first relation in $(C'''-i)$ of Theorem \ref{lem1} can be rewritten as follows
\[2 \int \log \frac{\sqrt{1+x^2}\sqrt{1+y^2}}{|x-y|} d\lambda_1(y) - \int \log \frac{ \sqrt{1+y^2}}{|x-y|} d\lambda_2(y) +
 \]
 \[\varphi(x) - 2\log(1+x^2) = \mathfrak{w}_1, \qquad x \in \supp(\lambda_1).\]
If $x \geq 1$, we have $\sqrt{1+y^2}/|x-y| \leq 1, y \in \R_-$ and taking \eqref{desig2} into consideration we obtain from the previous equality
\[\varphi(x) - 2\log(1+x^2) \leq \mathfrak{w}_1, \qquad x \in \supp(\lambda_1), \qquad x \geq 1.\]
Consequently, $\supp(\lambda_1)$ must be a compact set when\eqref{conditionInf} takes place. Condition \eqref{cond4} immediately follows for $\lambda_1$.

\medskip

Now, assume that $\supp(\lambda_2)$ is also compact. Then, $\lambda_2$ verifies \eqref{cond4} and
\[ \lim_{x \to \infty} {\mathcal{W}}_2^{\vec{\lambda}}(x) = \int \log (1+y^2) d\lambda_2(y) - \frac{1}{2} \int \log (1+y^2) d\lambda_1(y).
\]
In particular, taking the limit as $x \to -\infty$ along $\R_-$ from the second part of $(C'''-ii)$ we have that
$\int \log (1+y^2) d\lambda_2(y) - \frac{1}{2} \int \log (1+y^2) d\lambda_1(y) \geq {\mathfrak{w}}_2$. According to the first part of $(C'''-ii)$, ${\mathcal{W}}_2^{\vec{\lambda}}(x) \leq {\mathfrak{w}}_2$ on $\supp(\lambda_2)$. However, ${\mathcal{W}}_2^{\vec{\lambda}}$ is subharmonic in $\overline{\C} \setminus \supp (\lambda_2)$ and continuous on $\C$. By the maximum principle for subharmonic
function this means that ${\mathcal{W}}_2^{\vec{\lambda}} \equiv {\mathfrak{w}}_2$ on all $\C$ which is false. Therefore, $\supp(\lambda_2) = \R_-$ as claimed.

\medskip

In order to prove that $\lambda_2$ verifies \eqref{cond4} use $(C'''-ii)$ and argue as in Theorem \ref{lem1scalar} for proving that $\lambda$ satisfies \eqref{cond4}.
\hfill $\Box$

\medskip
{\bf Proof of Theorem \ref{teo:equil}}. Under the present assumptions, from the last assertions of Corollary \ref{varios} we know that  $\vec{\lambda} \in  {\mathfrak{M}}^*(\sigma) \subset \widetilde{\mathfrak{M}}(\sigma)$. The combined statements of Theorem \ref{lem1} and Corollary \ref{varios} give all but the last assertion of Theorem \ref{teo:equil}. Take into account that
\[ 2{\mathcal{U}}^{\lambda_1} - {\mathcal{U}}^{\lambda_2} + \varphi = 2U^{\lambda_1} - U^{\lambda_2} + \varphi + C_1, \qquad 2{\mathcal{U}}^{\lambda_2} - {\mathcal{U}}^{\lambda_1}   = 2U^{\lambda_2} - U^{\lambda_1} +  C_2,.
\]
where
\begin{align*}
 C_1 &= \int \log(1+y^2) d\lambda_1(y) - \frac{1}{2} \int \log(1+y^2) d\lambda_2(y),\\
 C_2 & = \int \log(1+y^2) d\lambda_2(y) - \frac{1}{2} \int \log(1+y^2) d\lambda_1(y).
\end{align*}
Thus
\[ w_1(\sigma,\varphi) = {\mathfrak{w}}_1(\sigma,\varphi) - C_1, \qquad w_2(\sigma,\varphi) =  {\mathfrak{w}}_2(\sigma,\varphi) - C_2.
\]
If $\int \log(1+y^2) d\sigma(x) = +\infty$, combining the arguments employed in the proof of $(c)$ and $(h)$ in Corollary \ref{miscelanea} it follows that $w_2(\sigma,\varphi) = 0$.
\hfill $\Box$


\section{Proof of Theorem \ref{teo:weak}}

{\bf Proof.} The sequences of zero counting measures $\left(\nu_{Q_{n}}\right), \left(\nu_{Q_{n,2}}\right), n \in \mathbb{Z}_+$ belong to $\mathcal{M}_1^+(\R_+)$ and $\mathcal{M}_1^+(\R_-)$, respectively. By Helly's selection theorem, there exists a sequence of indices $\Lambda \subset \mathbb{Z}_+$ and positive measures $\lambda_1^*, \lambda_2^*,$ $|\lambda_1^*| \leq 1, |\lambda_2^*| \leq 1$ such that
\begin{equation} \label{eq:13}
\lim_{n \in \Lambda} \nu_{Q_{n} } = \lambda_1^*, \qquad \lim_{n \in \Lambda} \nu_{Q_{n,2} } = \lambda_2^*.
\end{equation}
in the vague topology of measures. That is, for any continuous functions $f, g$ on $\R_+$ and $\R_-$, respectively, with compact support
\begin{equation} \label{weak}
\lim_{n \in \Lambda} \int f d \nu_{Q_{n} }= \int f d\lambda_1^*, \qquad \lim_{n \in \Lambda} \int g d \nu_{Q_{n,2} }= \int g d\lambda_2^*
\end{equation}
It easily follows that \eqref{weak} also holds for any $f \in \mathcal{C}_0(\R_+), g \in \mathcal{C}_0(\R_-)$ (the class of continuous functions on the indicated sets with limit equal to $0$ at infinity).

\medskip

In principle, it may occur that $|\lambda_1^*| < 1$ or $|\lambda_2^*| < 1$, but we will show that under our assumptions this is not the case. Moreover, we will show that $(2\lambda_1^*,\lambda_2^*) \in \mathfrak{M}^*(\sigma)$ and solves problem (C) in Theorem \ref{teo:equil}. After this is done, from uniqueness it follows that all convergent subsequences verifying \eqref{eq:13} have the same limit and the corresponding measures are precisely $\lambda_1/2$ and $\lambda_2$ where $(\lambda_1,\lambda_2)$ is the solution of Theorem \ref{teo:equil}. Then, since the limit measures in \eqref{eq:13} have mass one from \cite[Theorems 6.21, 6.22]{D} it follows that \eqref{weak} takes place for all bounded continuous functions $f, g$ on $\R_+, \R_-$, respectively, which amounts to \eqref{lessweak}.

\medskip

We begin by showing that $\lambda_2^* \leq \sigma$. Indeed, between two consecutive mass points of the discrete measure $\sigma_{2,n}$ there may be at most one zero of $Q_{n,2} $. Choose $-\infty < T_1 < T_2 \leq 0$, then from \eqref{limsigma} it follows that
\[ \limsup_n \int_{[T_1,T_2]} d \nu_{Q_{n,2}} \leq  \lim_{n} \frac{1}{n} \int_{[T_1,T_2]} d \left(\sum_{k \geq 1} \delta_{\xi_{k,n}}\right) = \int_{[T_1,T_2]} d \sigma. \]
On the other hand, since ${U}^{\sigma|_K}$ is continuous on $\mathbb{C}$ for every compact subset $K$ of $\R_-$ it follows that  $\sigma$ has no mass points; therefore, $\limsup_n \nu_{Q_{n,2}}(\{T\}) = 0 =\sigma(\{T\}) $ for each $T \in \mathbb{R}_-$.
These facts and the second part of \eqref{weak} imply that $\lambda_2^*  \leq \sigma$; whence, $\mathcal{U}^{\lambda_2^*}$ is continuous on $\C$ by Lemma \ref{cont}. Additionally, $\lambda_2^*$ satisfies \eqref{cond5} since $\sigma$ verifies it (see (iii)).

\medskip

Our next goal is to deduce the variational relations. We start with $\R_+$. To this aim we use the theorem on page 124 in \cite{GR1}.  From \eqref{eq:3^*}, it follows that
\[ \int  \frac{|Q_{n}(x)|^2}{|Q_{n,2}(x)|} {C_n x^{\alpha} s_1'(d_n x)\frac{d x}{x^{\alpha}}}  \leq \int  \frac{|Q(x)|^2}{|Q_{n,2}(x)|} {C_n x^{\alpha}s_1'(d_n x)\frac{d x}{x^\alpha}},\] with \[ C_n = \prod_{Q_{n,2}(x_{n,k}) =0} \sqrt{1 + x_{n,k}^2},
\]
for any monic polynomials $Q,\,\,\, \deg Q = 2n$. So $Q_{n}$ is the monic polynomial of degree~$2n$ that minimizes the $L_2$ norm with respect to the varying weight
\[ \frac{C_n x^{\alpha} s_1'(d_n x) }{|Q_{n,2}(x)|} \frac{d x}{x^{\alpha}}
\]
Since $\alpha < 1$ the measure $dx/x^{\alpha}$ is locally integrable on $\R_+$.

\medskip

We have
\[   g_n(x) := \frac{1}{n} \log  \frac{|Q_{n,2}(x)|}{C_n}  =  - \int \log \frac{\sqrt{1 + y^2}}{|x - y|} d\,\nu_{Q_{n,2}}(y),
\]
and $\log \frac{\sqrt{1 + y^2}}{|x - y|} \in \mathcal{C}_0(\R_-)$ for every $x >0$. From
\eqref{weak} we have
\begin{equation} \label{a} \lim_{n \in \Lambda} \frac{1}{2n} \log \left(\frac{|Q_{n,2}(x)|}{C_n}\right)^{1/2} = -\frac{1}{4} \int \log \frac{\sqrt{1 + y^2}}{|x - y|} d\lambda_2^*(y) = - \frac{1}{4}{\mathcal{U}}^{\lambda_2^*}(x)
\end{equation}
pointwise on $(0,+\infty)$.  On the other hand, if $0 < x < x' < + \infty$
\[ \left|\int \log \frac{\sqrt{1 + y^2}}{|x - y|} d\,\nu_{Q_{n,2}}(y) - \int \log \frac{\sqrt{1 + y^2}}{|x' - y|} d\,\nu_{Q_{n,2}}(y)\right| =  \int \log \frac{x' - y}{x - y} d\,\nu_{Q_{n,2}}(y) =
\]
\[ \int \log \left(1 + \frac{x' - x}{x - y}\right) d\,\nu_{Q_{n,2}}(y) < (x' -x) \int  \frac{ d\,\nu_{Q_{n,2}}(y)}{x - y} \leq \frac{x' -x}{x},
\]
which means that the family of functions $\left(g_n \right), n \in \mathbb{N}, $ is equicontinuous on compact subsets of $(0,+\infty)$. Therefore, \eqref{a} takes place uniformly on each compact subset of $(0,+\infty)$. Let us show that indeed \eqref{a} holds true uniformly on each compact subset of $\R_+$. It remains to show that this is true, for example, on the interval $[0,1/2]$.

\medskip

Take $\delta \in (0,1/2)$ and $x \in [0,1/2]$. Then
\[ \left| \int \log \frac{\sqrt{1+y^2}}{|x-y|} d \lambda_2^*(y) - \int \log \frac{\sqrt{1+y^2}}{|x-y|} d \nu_{Q_{n,2}}(y) \right| \leq
\]
\[ \left| \int_{|y| \geq \delta} \log \frac{\sqrt{1+y^2}}{|x-y|} d \lambda_2^*(y) - \int_{|y| \geq \delta} \log \frac{\sqrt{1+y^2}}{|x-y|} d \nu_{Q_{n,2}}(y) \right| +
\]
\[\left| \int_{|y| \leq \delta} \log \frac{\sqrt{1+y^2}}{|x-y|} d \lambda_2^*(y)\right| + \left| \int_{|y| \leq \delta} \log \frac{\sqrt{1+y^2}}{|x-y|} d \nu_{Q_{n,2}}(y) \right| \leq
\]
\[ \left| \int_{|y| \geq \delta} \log \frac{\sqrt{1+y^2}}{|x-y|} d \lambda_2^*(y) - \int_{|y| \geq \delta} \log \frac{\sqrt{1+y^2}}{|x-y|} d \nu_{Q_{n,2}}(y) \right| +
\]
\[ \int_{|y| \leq \delta} \log \frac{\sqrt{1+y^2}}{|y|} d \lambda_2^*(y)  + \left| \int_{|y| \leq \delta} \log \frac{ 1}{|y|} d \nu_{Q_{n,2}}(y) \right|   + \log \sqrt{1 + \delta^2}.
\]
Fix $\varepsilon > 0$. Since $\mathcal{U}^{\lambda_2^*}$ is continuous on $\mathbb{C}$ (in particular at $x=0$) we have that $\log \frac{\sqrt{1+y^2}}{|y|}$ is integrable with respect to $\lambda_2^*$ and $0$ is not a mass point of $\lambda_2^*$; consequently, for all $\delta$ sufficiently small it follows that
\[  \int_{|y| \leq \delta} \log \frac{\sqrt{1+y^2}}{|y|} d \lambda_2^*(y)  < \varepsilon.
\]
The last term on the last line is obviously $< \varepsilon$ for all sufficiently small $\delta$. Let us show that the same is true for the middle term.

\medskip

Between two mass point of $\sigma_{2,n}$ there is at most one zero on $Q_{n,2}$; therefore,
\[  \left| \int_{|y| \leq \delta} \log \frac{ 1}{|y|} d \nu_{Q_{n,2}}(y) \right| =  \frac{1}{n}\left| \log \prod_{|x_{n,k}|\leq \delta}|x_{n,k}|\right| \leq \left| \log \left(\prod_{|\xi_{n,k}|\leq \delta}|\xi_{n,k}|\right)^{1/n}\right|. \]
Let $\rho = \min\{\rho(x): x\in [-\delta,0]\} (> 0)$ where $\rho(x)$ is the function which appears in condition (i) in Section 2. According to (i)
\begin{equation}
\label{acotxi}
|\xi_{k,n}| = |\xi_{k,n} - \xi_{k-1,n}| + \cdots + |\xi_{1,n}|\geq k\rho/n.
\end{equation}
Let $\ell_n$ be the number of $\xi_{k,n}$ in $[-\delta,0]$. From \eqref{limsigma} $\lim_{n\to \infty} \ell_n/n  = \sigma([-\delta,0])$. Condition (i) also implies that $\ell_n \leq n\delta/\rho$; consequently $\lim_n \ell_n^{1/n} = 1$. Using Stirling's formula and \eqref{acotxi}
\[1 > \left(\prod_{|\xi_{n,k}|\leq \delta}|\xi_{n,k}|\right)^{1/n}\geq  \left(\frac{\rho}{n}\frac{2\rho}{n}\cdots\frac{\ell_n\rho}{n}\right)^{1/n} = \left(\frac{\rho}{n}\right)^{\ell_n/n}\left(\ell_n!\right)^{1/n} \geq \]
\[\left(\frac{\rho}{e}\right)^{\ell_n/n} \left(\frac{\ell_n}{n}\right)^{\ell_n/n}\ell_n^{1/(2n)}\mathcal{O}(1)^{1/n}.\]
Consequently,
\begin{equation}
  \label{limQn}
  \begin{split}
  \limsup_{n\to \infty}\left| \int_{|y| \leq \delta} \log \frac{ 1}{|y|} d \nu_{Q_{n,2}}(y) \right| \leq &\left|\log \left(\lim_{n\to \infty} \left(\frac{\rho}{e}\right)^{\ell_n/n} \left(\frac{\ell_n}{n}\right)^{\ell_n/n}\ell_n^{1/(2n)}\mathcal{O}(1)^{1/n}\right)\right| =\\
&\left|\sigma[-\delta,0] \log \frac{\rho \sigma[-\delta,0]}{e}\right|, \end{split}\end{equation}
which tends to zero as $\delta \to 0$.

\medskip

Therefore, we can choose and fix $\delta \in (0,1/2)$ such that
\[ \int_{|y| \leq \delta} \log \frac{\sqrt{1+y^2}}{|y|} d \lambda_2^*(y)  + \left| \int_{|y| \leq \delta} \log \frac{ 1}{|y|} d \nu_{Q_{n,2}}(y) \right|   + \log \sqrt{1 + \delta^2} < 3\varepsilon.
\]
For $\delta$ fixed, it is easy to show that
\[ \lim_{ n \in \Lambda} \int_{|y| \geq \delta} \log \frac{\sqrt{1+y^2}}{|x-y|} d \nu_{Q_{n,2}}(y) =   \int_{|y| \geq \delta} \log \frac{\sqrt{1+y^2}}{|x-y|} d \lambda_2^*(y)
\]
uniformly with respect to $x \in [0,1/2]$. Putting all this together we find that for any $\varepsilon > 0$ there exists $n_0$ such that if $n \geq n_0, n \in \Lambda,$ then
\[ \left| \int \log \frac{\sqrt{1+y^2}}{|x-y|} d \lambda_2^*(y) - \int \log \frac{\sqrt{1+y^2}}{|x-y|} d \nu_{Q_{n,2}}(y) \right| \leq 4 \varepsilon
\]
independent of $x \in [0,1/2]$. Thus \eqref{a} takes place uniformly on each compact subset of $\R_+$ as we wanted to prove.

\medskip

Set
\[ {f}_n(x) := \frac{1}{2n} \log \left(\frac{|Q_{n,2}(x)|}{C_n x^{\alpha} s_1'(d_n x)}\right)^{1/2}.
\]
What was proved in the previous sentence and \eqref{b} imply
\[
\lim_{n \in \Lambda}  \frac{1}{2n} \log \left(\frac{|Q_{n,2}(x)|}{C_n x^{\alpha} s_1'(d_n x)}\right)^{1/2} = \frac{1}{4}(\varphi(x) - {\mathcal{U}}^{\lambda_2^*}(x)),
\]
uniformly on each compact subset of $\R_+$. In particular,  for any closed interval $\Delta \subset \R_+$
\begin{equation} \label{delta} \lim_{n \in \Lambda} \min_{x \in \Delta} {f}_n(x) = \min_{x \in \Delta} \frac{1}{4}(\varphi(x) - {\mathcal{U}}^{\lambda_2^*}(x)).
\end{equation}
For $x \geq 1, y \leq 0,$ we have that $\log \sqrt{1 + y^2}/(x-y) \leq 0 $; therefore, from \eqref{cond1} and \eqref{cond3} it follows that
\begin{equation} \label{c} \liminf_{x \to +\infty} \frac{\varphi(x) - {\mathcal{U}}^{\lambda_2^*}(x)}{4 \log x} > 1, \qquad \liminf_{n\in \Lambda, x \to +\infty} \frac{ {f}_n(x)}{\log x} > 1.
\end{equation}

\medskip

Relations \eqref{delta} and \eqref{c} certify that a) and b) on page 124 of  \cite{GR1} are fulfilled. Therefore, using  the lemma on page 121 and the theorem on page 124 in \cite{GR1} it follows that $\lambda_1^*$ is the unique probability measure on $\R_+$ which solves the extremal problem
\begin{equation} \label{extremal1*} U^{\lambda_1^*}(x) + \frac{1}{4}(\varphi(x) - {\mathcal{U}}^{\lambda_2^*}(x))  \left\{
\begin{array}{ll}
 = w_1^*, & x \in \supp(\lambda_1^*), \\
 \geq w_1^*, & x \in \R_+,
\end{array}
\right.
\end{equation}
for some constant $w_1^*$,  and (recall that $\deg Q_{n} = 2n$)
\begin{equation} \label{eq:15}
\lim_{n \in \Lambda} \left(\int  \frac{|Q_n (x)|^2}{|Q_{n,2} (x)|} {C_n s_1'(d_n x)dx}  \right)^{1/4n} = e^{-w_1^*}.
\end{equation}
The arguments employed on \cite[page 127]{GR1} to prove the main theorem allow to conclude that for each $\varepsilon > 0$ there exists $R > 0$ such that
\begin{equation} \label{eq:15*}
\liminf_{n \in \Lambda} \left(\int_0^R  \frac{|Q_n (x)|^2}{|Q_{n,2} (x)|} {C_n s_1'(d_n x)dx}  \right)^{1/4n} \geq e^{-w_1^*-\varepsilon}.
\end{equation}
The first part of \eqref{c} guarantees that $\supp (\lambda_1^*)$ is a compact subset of $[0,+\infty)$. This is shown in \cite{GR2} (see also \cite[Theorem 1.3.1]{ST}, or even Corollary \ref{miscelanea}(a) applied to measures supported on $\R_+$). Notice that \eqref{extremal1*} and the continuity of $\varphi$ and ${\mathcal{U}}^{\lambda_2^*}$ on $\R_+$ imply that $\mathcal{U}^{\lambda_1^*}$ is continuous on $\supp(\lambda_1^*)$ and thus on all $\C$. Using the compactness  of $\supp (\lambda_1^*)$, we have
\[I(\lambda_1^*) < +\infty, \qquad \int \log (1+y^2) d \lambda_1^*(y) < \infty.
\]

\medskip

Now, let us obtain the variational relations on $\R_-$. The varying discrete measure with respect to which $Q_{n,2}$ is orthogonal, see \eqref{eq:4^*} and \eqref{sigma2scaled}, may be regarded as
\[
\sum_{k=1}^{\infty} \frac{\beta_k \eta_{n,k}}{|\xi_{k,n}|} \frac{ D_n}{|Q_n (\xi_{k,n})|} \delta_{\xi_{k,n}}(t),
   \,\,\,   \eta_{n,k} = \int_{\mathbb{R}_+} \frac{|Q_n (x)|^2}{|Q_{n,2} (x)| } \frac{C_n s_1'(d_n x) dx}{1 - (x/\xi_{k,n})},\] \[ D_n = \prod_{Q_n(y_{n,k}) = 0} \sqrt{1 + y^2_{n,k}}.
\]

\medskip

Since $\sum_{k=1}^\infty \beta_k/t_k < + \infty$ and $ \lim_n d_n^{1/n}= 1$, we have
\[\lim_{n \to \infty}\left(\sum_{k=1}^{\infty} \frac{\beta_k }{|\xi_{k,n}|}\right)^{1/n} = 1.\]
Using (\ref{eq:15})
\begin{equation} \label{eq:17}  \limsup_{n \in \Lambda} \eta_{n,k}^{1/n}   \leq e^{-4w_1^*}.
\end{equation}
On the other hand, from  \eqref{eq:15*} for any $\varepsilon > 0$ we can choose $R > 0$ such that
\begin{equation} \label{eq:17*}   \liminf_{n \in \Lambda} \eta_{n,k}^{1/n} \geq  \liminf_{n \in \Lambda} \left(\int_0^R \frac{|Q_n(x)|^2}{|Q_{n,2}(x)|} \frac{C_n s_1'(d_n x) dx}{1 - (Rd_n/t_1) }\right)^{1/n} \geq e^{-4w_1^*-4\varepsilon}.
\end{equation}
From \eqref{eq:17} and \eqref{eq:17*} it follows that
\begin{equation} \label{eq:17**} \lim_{n \in \Lambda} \eta_{n,k}^{1/n} = e^{-4w_1^*},
\end{equation}
uniformly on $k$.

\medskip

Since $\log\frac{\sqrt{1 + y^2}}{|x-y|} \in \mathcal{C}_0(\R_+)$ for every $x <0$,  arguing as we did for the sequence of polynomials $\left(Q_{n,2}\right)$, we have
\begin{equation} \label{eq:18}\lim_{n \in \Lambda} \left(\frac{|Q_{n}(x)|}{D_n}\right)^{1/n} = e^{- 2{{\mathcal{U}}}^{\lambda_1^*}(x)},
\end{equation}
uniformly on each compact subset of $(-\infty,0)$.
Set $\phi(x) := 4w_1^* - \mathcal{U}^{2\lambda_1^*}(x)$. Using (\ref{eq:17**}) and (\ref{eq:18}), we obtain
\begin{equation} \label{eq:19}
\lim_{n \in \Lambda}  \left(\frac{\eta_{n,k} D_n}{|Q_n (\xi_{k,n})|}\right)^{1/n} - e^{-\phi(\xi_{k,n})} = 0
\end{equation}
uniformly on each compact $K \subset (-\infty,0)$  and $k$ such that $\xi_{k,n} \in K$.

\bigskip

Let $\lambda \in \mathcal{M}^*(\sigma)$ be the extremal solution of Corollary \ref{potescalar} with $\sigma$  as in Theorem \ref{teo:weak} and $\phi(x) := 4w_1^* - \mathcal{U}^{2\lambda_1^*}(x)$. In Theorem \ref{teo:weak} we have assumed that $0 \not \in \supp (\sigma \setminus \lambda_2)$ so we will assume here that $0 \not\in \supp (\sigma \setminus \lambda)$. We will show that $\lambda_2^* = \lambda$ using modified versions of some results which appear in \cite[Lemmas 5.3, 5.5, and 3.2]{DrSa}. In \cite{DrSa} the corresponding $\lambda$ had compact support while in our case the support is $\R_-$. More exactly, applying Corollaries \ref{potescalar} and \ref{miscelanea}(c),(h), it follows that there  exist  $ {\lambda}  \in  {\mathcal{M}}^*(\sigma)$ and a constant ${ {w}}  = { {w}} (\sigma,\phi) = 4 w_1^*- \int \log(1+y^2) d\lambda_1^*(y)$  such that
\begin{equation}
\label{extremal2*}
2 {{ {U}}}^{ \lambda }(x) + \phi(x)
\left\{
\begin{array}{ll}
\leq w , & x \in \supp( \lambda ) = \R_-, \\
= w , & x \in \supp(\sigma -  \lambda ),
\end{array}
\right.
\end{equation}
and $\supp(\sigma -  \lambda )$ is unbounded.

\medskip

Set
\[
\|Q_{n,2}\|_{2,n} = \left(\sum_{k=1}^{\infty} |Q_{n,2}(\xi_{k,n})|^2\frac{\beta_k \eta_{n,k}}{|\xi_{k,n}|} \frac{ D_n}{|Q_n (\xi_{k,n})|} \right)^{1/2}.
\]
Let us show that
\begin{equation}
\label{Qn2}
\limsup_{n \in \Lambda} \|Q_{n,2}\|_{2, n}^{1/n } \leq e^{-w(\sigma,\phi)/2}.
\end{equation}
For this, we follow the approach in \cite[Lemma 5.3]{DrSa}

\medskip

Fix $\varepsilon > 0$. Set $w = w(\sigma,\phi)$. Choose $A \supset \supp(\sigma - \lambda)$ to be the union of finitely many closed intervals such that $2 {{ {U}}}^{ \lambda }(x) - {{ \mathcal{U}}}^{ 2\lambda_1^* }(x) > w - \varepsilon, x \in A$, and $0 < \lambda(A) < 1$. The existence of such a set is guaranteed because, according to \eqref{unif2} and Lemma~\ref{cont},
\begin{equation}
\label{wsf}\lim_{x \to \infty} 2 {{ {U}}}^{ \lambda }(x) + \phi(x) = w(\sigma,\phi) = 4w_1^* - \int \log(1+y^2) d\lambda_1^*(y),
\end{equation}
when $x\to \infty$ in any direction; in particular as $x \to -\infty$. Moreover, because $\R_- \setminus \supp(\sigma - \lambda) \neq \emptyset$  since $|\lambda| = 1 < |\sigma|$. Since $0 \not\in \supp(\sigma \setminus \lambda)$ we can take $A$ so that $0\in \R_- \setminus A$.

\medskip

Let $ \widetilde{\lambda} = \lambda|_{\R_- \setminus A}$, and $\widetilde{\sigma}_n = \frac{1}{n} \sum_{k \geq 1} \delta_{\xi_{k,n}}|_{\R_- \setminus A}$. $\R_- \setminus A$ is a compact set and from \eqref{limsigma}, we obtain $\lim_{n \in \Lambda} \widetilde{\sigma}_n = \widetilde{\lambda}  $ in the vague topology. In particular, $\lim_{n\in \Lambda} \frac{m_n}{n} = \lambda(\R_- \setminus A) < 1$, where $m_n$ is the number of points $\xi_{k,n}$ which lie in $\R_- \setminus A$. Therefore, there exists $n_0$ such that $m_n < n$ for $n\geq n_0, n \in \Lambda.$

\medskip

Let $P_n$ be a monic polynomial of degree $n$ whose zeros consist of the $m_n$ points $\xi_{k,n} \in \R_- \setminus A$ and $n - m_n$ points in $A$ chosen so that $\lim_{n \in \Lambda} \nu_{P_n} = \lambda$ in the vague topology. It is sufficient to discretize $\lambda$ on $A$. Since $\lambda \in \mathcal{M}^*(\sigma)$ and $\log (1+y^2)$ is positive and decreasing in $\R_-$ one can also ensure that
\begin{equation}
\label{Pn}
  \lim_{n\in \Lambda} \int \log(1+y^2) d \nu_{P_n}(y) = \int \log(1+y^2) d \lambda(y)
\end{equation}

\medskip

For $n \geq n_0, n \in \Lambda$ we have
\[\|Q_{n,2}\|_{2,n}^{2/n} \leq \|P_n\|_{2,n}^{2/n} \leq \left(\sum_{\xi_{k,n} \in A} |P_n(\xi_{k,n})|^2\frac{\beta_k \eta_{n,k}}{|\xi_{k,n}|} \frac{ D_n}{|Q_n (\xi_{k,n})|} \right)^{1/n} \leq \]
\[ \left(\sum_{k=1}^{\infty} \frac{\beta_k}{|\xi_{k,n}|}\right)^{1/n}\exp\left\{-  \left( 2U^{\nu_{P_n}}(\xi_n) - 2\mathcal{U}^{\nu_{Q_n}}(\xi_n) + \frac{1}{n} \log \eta_n \right)\right\},
\]
where $\xi_n$ is a point $\xi_{k,n} \in A$ for which
\[  2U^{\nu_{P_n}}(\xi_n) - 2 \mathcal{U}^{\nu_{Q_n}}(\xi_n) + \frac{1}{n} \log \eta_n = \min_{\xi_{k,n} \in A} \left(2U^{\nu_{P_n}}(\xi_{k,n}) - 2\mathcal{U}^{\nu_{Q_n}}(\xi_{k,n}) + \frac{1}{n} \log \eta_{n,k}\right)\]
and $\eta_n$ is the $\eta_{n,k}$ corresponding to that point.

\medskip

Let $\xi \in A$ be any limit point of the sequence $(\xi_n), n \in \Lambda$; that is, $\lim_{n \in \Lambda'} \xi_n = \xi (\neq 0)$ with $\Lambda'\subset \Lambda$. Then, using \eqref{eq:19}, \eqref{Pn}, and the principal of descent
\[  \liminf_{n \in \Lambda'} \left(2U^{\nu_{P_n}}(\xi_n) - 2\mathcal{U}^{\nu_{Q_n}}(\xi_n) + \frac{1}{n} \log \eta_n \right)\geq 2U^{\lambda}(\xi) + \phi(\xi) \geq w(\sigma,\phi) - \varepsilon.
\]
Consequently,
\[\limsup_{n\in \Lambda} \|Q_{n,2}\|_{2,n}^{2/n} \leq e^{-w(\sigma,\phi) + \varepsilon}.\]
Letting $\varepsilon \to 0$ we obtain \eqref{Qn2}.

\medskip

Now, using the scheme employed in \cite[Lemmas 3.2, 5.5]{DrSa} we prove that
\begin{equation}
\label{Qn2-}
\liminf_{n \in \Lambda} \|Q_{n,2}\|_{2, n}^{2/n } \geq e^{-F_{\lambda_2^*}},
\end{equation}
where
\[
\label{F*}
F_{\lambda_2^*} = \max \{C \in \R: 2 U^{\lambda_2^*}(x) + \phi(x) \geq C\,\,\, \mbox{holds}\,\,\,(\sigma -\lambda_2^*)\,\, \mbox{a.e.} \}.
\]

\medskip

Let $x_0 \in \supp (\sigma \setminus \lambda_2^*)\setminus \{0\}$. Fix $0 < \varepsilon < 1/2,$ sufficiently small so that $[x_0 - \varepsilon, x_0 +\varepsilon] \subset (-\infty,0)$. Set $\Delta_\varepsilon = (x_0 - \varepsilon, x_0 +\varepsilon)$. Now, choose $0 <\delta < \varepsilon$ and set $\Delta_\delta := (x_0-\delta,x_0 +\delta)$. Choose $M > 0$ such that $-M < x_0 - \varepsilon -1$. Define
\[ Q_{n,2}^{(1)}(x) := \prod_{y_{n,k} \in \Delta_\varepsilon}(x - y_{n,k}), \qquad Q_{n,2}^{(2)}(x) := \prod_{y_{n,k} \in [-M,0] \setminus \Delta_\varepsilon}(x - y_{n,k}), \]\[Q_{n,2}^{(3)} := Q_{n,2}/(Q_{n,2}^{(1)}Q_{n,2}^{(2)}). \]
Since $x_0 \,\in\, \supp\,(\sigma -\lambda_2^*)$ we have that $q:= (\sigma -\lambda_2^*)\,(\Delta_\delta)\, >\, 0$. Let $\ell_n$ be the number of zeros of $Q_{n,2}$ in $\Delta_\delta$ and $m_n$ be the number of $\xi_{k,n}$ in $\Delta_\delta$. Then $\lim_{n\to \infty} (m_n -\ell_n)/n = q$. Since the intervals $((\xi_{k,n} + \xi_{k-1,n})/2,(\xi_{k,n} + \xi_{k+1,n})/2)$ around the mass point $\xi_{k,n}$ are disjoint for all sufficiently large $n$ there exists at least one interval containing no zeros of $Q_{n,2}$ whose corresponding mass point is in~$\Delta_\delta$. Denote this mass point by $\xi_n^*$ and its adjacent mass points by $\xi_n^{(1)}$ and $\xi_n^{(2)}$.
Now, using again that between two mass points of $\sigma_{2,n}$ there is at most one zero of $Q_{n,2}$, one obtains
\[ |Q_{n,2}^{(1)}(\xi_n^*)|^{1/n} \geq \left(\frac{|\xi_{n}^* - \xi_{n}^{(1)}||\xi_{n}^* - \xi_{n}^{(1)}|}{4}\right)^{1/n}\left(\prod_{\xi_n^*\neq \xi_{k,n}\in \Delta_\varepsilon}|\xi_n^* - \xi_{k,n}|\right)^{1/n} \geq \]
\[(1/4)^{1/n}\left(\prod_{\xi_n^*\neq \xi_{k,n}\in \Delta_\varepsilon}|\xi_n^* - \xi_{k,n}|\right)^{2/n}.\]
Let $p_n$ be the number of $\xi_{k,n} > \xi_n^*$ in $\Delta_\varepsilon$, $q_n$ be the number of $\xi_{k,n} < \xi_n^*$ in $\Delta_\varepsilon$. Using \eqref{limsigma}, we have $\lim_{n\to \infty} (p_n +q_n)/n = \sigma(\Delta_\varepsilon)$. Let $\rho := \inf\{\rho(x): x \in \Delta_\varepsilon\}$. The previous inequalities and (i) of Section 2 imply that
\[ |Q_{n,2}^{(1)}(\xi_n^*)|^{1/n} \geq (1/4)^{1/n} \left(\frac{\rho}{n}\right)^{2p_n/n} (p_n!)^{2/n}\left(\frac{\rho}{n}\right)^{2q_n/n} (q_n!)^{2/n} \geq
\]
\[ (1/4)^{1/n} \left(\frac{\rho}{n}\right)^{2(p_n+q_n)/n} \left((r_n -1)! \right)^{2/n}\]
where $r_n$ denotes the integer part of $(p_n+q_n)/2$.
From here, using Stirling's formula,   it is easy to deduce that
\begin{equation}
\label{cotainf1}
\liminf_{n \to \infty}|Q_{n,2}^{(1)}(\xi_n^*)|^{1/n} \geq \left(\frac{\rho \sigma(\Delta_{\varepsilon})}{2e}\right)^{2\sigma(\Delta_{\varepsilon}) }.
\end{equation}
Notice that the right hand tends to 1 as $\varepsilon \to 0$.

\medskip

We have
\begin{equation}
\label{cotainf}
\|Q_{n,2}\|_{2,n}^{2/n} = \left(\sum_{k=1}^{\infty} |Q_{n,2}(\xi_{k,n})|^2\frac{\beta_k \eta_{n,k}}{|\xi_{k,n}|} \frac{ D_n}{|Q_n (\xi_{k,n})|} \right)^{1/n} \geq
\end{equation}
\[ \left(|Q_{n,2}(\xi_{n}^*)|^2\frac{\beta_n^* \eta_{n}^*}{|\xi_{ n}^*|} \frac{ D_n}{|Q_n (\xi_{n}^*)|} \right)^{1/n} \geq \left(|Q_{n,2}^{(1)}(\xi_{n}^*)Q_{n,2}^{(2)}(\xi_{n}^*)|^2\frac{\beta_n^* \eta_{n}^*}{|\xi_{ n}^*|} \frac{ D_n}{|Q_n (\xi_{n}^*)|} \right)^{1/n},
\]
where $\beta_n^*, \eta_{n}^*$ are the values of $\beta_k$, and $\eta_{n,k}$, respectively, corresponding to $\xi_{k,n} = \xi_n^*$. In the last inequality we skip $Q_{n,2}^{(3)}$ because all its zeros are at distance greater than $1$ from $\xi_n^*$. Let us find a lower bound for
\[\left(|Q_{n,2}^{(2)}(\xi_{n}^*)|^2\frac{\beta_n^* \eta_{n}^*}{|\xi_{ n}^*|} \frac{ D_n}{|Q_n (\xi_{n}^*)|} \right)^{1/n}.\]

\medskip

Since $\nu_{Q_{n,2}^{(2)}}$ converges vaguely to $\lambda_2^*|_{[-M,0] \setminus \Delta_\varepsilon}, n \in \lambda,$  $U^{\nu_{Q_{n,2}^{(2)}}}$ converges uniformly on $\Delta_{\delta}$ to $U^{\lambda_2^*|_{[-M,0] \setminus \Delta_\varepsilon}}, n\in \Lambda$ , and $U^{\lambda_2^*|_{[-M,0]}}$ is continuous on $\R_-$ (in particular at $x_0$, recall that $\xi_n^* \in \Delta_\delta)$, given $\varepsilon$ we can find $\delta, 0 <\delta < \varepsilon,$ such that
\begin{equation}
\label{cotainf2}
\liminf_{n\in \Lambda} |Q_{n,2}^{(2)}(\xi_{n}^*)|^{2/n} \geq e^{- 2U^{\lambda_2^*|_{[-M,0]}}(x_0) - 2\varepsilon}.
\end{equation}
Also, because of the continuity of $\phi$ and \eqref{eq:19}, $\delta$ may be chosen so that $|\phi(x) - \phi(x_0)| < \varepsilon,  x\in \delta_\delta,$ and for all sufficiently large $n \in \Lambda$ and $k$ with $\xi_{k,n}\in \Delta_\delta $
\begin{equation}
\label{cotainf3}
\left|  \left(\frac{\eta_{n,k} D_n}{|Q_n (\xi_{k,n})|}\right)^{1/n} -  e^{-{\phi(\xi_{k,n})}}\right| < \varepsilon.
\end{equation}
so that \eqref{cotainf3} holds, in particular, for $\xi_n^*$ and $\eta_n^*$. Since $\xi_n^* \in \Delta_\delta$, we have $$\lim_{n\to\infty}|\xi_n^*|^{1/n} = 1.$$   On the other hand, from (i) if $\xi_{k,n} \in (x_0-\delta,x_0+\delta)$ then $k < n|x_0 - \delta|/\rho^*$ where $ \rho^* = \inf\{{\rho}(x): x \in [x_0 - \delta,0]\} > 0$ which combined with (ii)
implies that $\liminf_{n\to\infty} |\beta_n^*|^{1/n} \geq 1$.

\medskip

Using \eqref{cotainf1}-\eqref{cotainf3}, it follows that for all sufficiently small $\varepsilon > 0$ and $M >0$ sufficiently large
\begin{equation}
\label{cotainf4}
\liminf_{n \in \Lambda} \|Q_{n,2}\|_{2,n}^{2/n}  \geq \left(\frac{\rho \sigma(\Delta_{\varepsilon})}{2e}\right)^{2\sigma(\Delta_{\varepsilon}) }e^{- 2U^{\lambda_2^*|_{[-M,0]}}(x_0) - \phi(x_0)- 4\varepsilon}
\end{equation}
Now, \eqref{Qn2} and \eqref{cotainf4} imply that
\begin{equation}
\label{cotainf5}
 2U^{\lambda_2^*|_{[-M,0]}}(x_0) + \phi(x_0)+ 4\varepsilon - 2\sigma(\Delta_{\varepsilon}) \log \left(\frac{\rho \sigma(\Delta_{\varepsilon})}{2e}\right) \geq w(\sigma,\phi) > -\infty.
\end{equation}
Suppose that $\int \log (1+y^2) d\lambda_2^*(y) = \infty$. In this case it is easy to prove that $U^{\lambda_2^*|_{[-M,0]}}(x_0)$ tends to $-\infty$ as $M\to +\infty$ which contradicts \eqref{cotainf5}. Consequently, $\int \log (1+y^2) d\lambda_2^*(y) < \infty$. In this case $U^{\lambda_2^*}$ is well defined on all $\C$, and is continuous on $\R_-$; moreover,
\[ \lim_{M \to \infty} U^{\lambda_2^*|_{[-M,0]}}(x) = U^{\lambda_2^*}(x)\]
uniformly on any compact subset of $\C$. Making $M\to \infty$ and $\varepsilon \to 0$ from \eqref{cotainf5} it follows that
\[2U^{\lambda_2^* }(x_0) + \phi(x_0) \geq w(\sigma,\phi) = F_{\lambda}.\]
Now this occurs for every $x_0 \in \supp(\sigma \setminus \lambda_2^*) \setminus \{0\}$ and by continuity also at $0$ should this be an accumulation point of $\supp(\sigma \setminus \lambda_2^*)$. Consequently, $F_{\lambda_2^*} \geq F_\lambda$. From $(D'')$ of Corollary \ref{potescalar} we conclude that $F_{\lambda}= w(\sigma,\phi) = F_{\lambda_2^*} $; therefore,
\begin{equation}
\label{Qn2+-}
\lim_{n \in \Lambda} \|Q_{n,2}\|_{2, n}^{2/n} = e^{-F_{\lambda_2^*}},
\end{equation}
and $\lambda = \lambda_2^*$ according to the unicity statement in that part of Corollary \ref{potescalar}.

\medskip

Using  that $\int (1+y^2) d\lambda_1^*(y) < +\infty$ and $\int (1+y^2) d\lambda_2^*(y) < +\infty$, we can rewrite \eqref{extremal1*} and \eqref{extremal2*} in terms of $U^{\lambda_1^*}, U^{\lambda_2^*}$ and $\phi$  as follows:
\[2U^{2\lambda_1^*}(x) - {{U}}^{\lambda_2^*}(x) + \varphi(x) \left\{
\begin{array}{ll}
 = 4w_1^* + \frac{1}{2} \int \log (1+y^2) d\lambda_2^*(y), & x \in \supp(\lambda_1^*), \\
 \geq 4w_1^* + \frac{1}{2} \int \log (1+y^2) d\lambda_2^*(y), & x \in \R_+,
\end{array}
\right.\]
\[  2 {{ {U}}}^{ \lambda_2^* }(x) -  {{ {U}}}^{ 2\lambda_1^* }(x)
\left\{
\begin{array}{ll}
\leq 0 , & x \in \supp( \lambda_2^* ) = \R_-, \\
= 0 , & x \in \supp(\sigma -  \lambda_2^* ),
\end{array}
\right. \]
Therefore, the   pair $(2\lambda_1^*,\lambda_2^*)$ satisfies \eqref{EqCo1p} and \eqref{EqCo2p} in part $(C)$ of Theorem \ref{teo:equil}. This means that $(2\lambda_1^*,\lambda_2^*) = (\lambda_1,\lambda_2)$ is the extremal solution of Theorem \ref{teo:equil} and the extremal constants are $w_1 = 4w_1^* + \frac{1}{2} \int \log (1+y^2) d\lambda_2(y), w_2 = 0$. In particular, $|\lambda_1^*| = |\lambda_2^*| = 1$ and, as explained in the beginning of the proof, \eqref{lessweak} follows from \eqref{eq:13}. \hfill $\Box$

\medskip

\begin{rem}
\label{rem:l}
From \eqref{eq:15} and \eqref{Qn2+-} we also have
\begin{equation}
\label{nth3}
\lim_{n } \left(\int  \frac{|Q_n (x)|^2}{|Q_{n,2} (x)|} {C_n s_1'(d_n x)dx}  \right)^{1/n} = e^{-4w_1^*}, \qquad \lim_{n} \|Q_{n,2}\|_{2, n}^{2/n} = e^{-F_{\lambda_2^*}},
\end{equation}
where $F_{\lambda_2^*} = 4w_1^* - \frac{1}{2} \int \log (1+y^2) d\lambda_1(y)$ (see \eqref{wsf}). Direct computation gives
\[\|Q_{n,2}\|^{2/n} = (D_n C_n)^{1/n} \left| \int  \frac{Q_{n,2}^2(t)}{Q_{n}(t)}\int  \frac{Q_{n}^2(x)}{Q_{n,2}(x)} \frac{\sigma'_1(d_n x)d x}{x-t}d \sigma_{2,n}(t)\right|^{1/n}. \]
Therefore, using \eqref{nth3} we could establish that
\begin{equation} \label{nth1} \lim_n \left| \int \frac{Q_{n}^2(x)}{Q_{n,2}(x)} {\sigma'_1(d_n x)d x}\right|^{1/n} = e^{- w_1},
\end{equation}
and
\begin{equation} \label{nth2}  \lim_n \left| \int  \frac{Q_{n,2}^2(t)}{Q_{n}(t)}\int  \frac{Q_{n}^2(x)}{Q_{n,2}(x)} \frac{\sigma'_1(d_n x)d x}{x-t}d \sigma_{2,n}(t)\right|^{1/n} = e^{-w_1},
\end{equation}
where $w_1$ is the corresponding equilibrium constant from \eqref{EqCo1p} (here $w_2 =0$), if we could prove that
\[ \lim_n C_n^{1/n} = e^{\frac{1}{2}\int \log(1+y^2) d\lambda_2(y)}, \qquad \lim_n D_n^{1/n} = e^{\frac{1}{2}\int \log(1+y^2) d\lambda_1(y)}.\]
In order to do this, it is necessary to obtain some bound on the rate of growth of the largest zeros of the polynomials $Q_n$ and $Q_{n,2}$.
\end{rem}


\end{document}